\documentclass[12pt,leqno,a4paper]{amsart}
\usepackage{amssymb}
\usepackage{enumerate}
\usepackage{txfonts}
\usepackage{mathrsfs}
\usepackage{hyperref,cite}
\usepackage{xcolor}
\usepackage{framed}

\newcommand{\FF}{\mathbb{F}}
\newcommand{\ZZ}{\mathbb{Z}}

\newcommand{\ba}{\mathbf{a}}
\newcommand{\bb}{\mathbf{b}}
\newcommand{\bG}{{\mathbf{G}}}
\newcommand{\bH}{\mathbf{H}}
\newcommand{\bL}{{\mathbf{L}}}
\newcommand{\bP}{\mathbf{P}}
\newcommand{\bT}{\mathbf{T}}
\newcommand{\bZ}{\mathbf{Z}}

\newcommand{\cA}{\mathcal{A}}
\newcommand{\cB}{\mathcal{B}}

\newcommand{\cF}{\mathcal{F}}
\newcommand{\cG}{\mathcal{G}}
\newcommand{\cH}{\mathcal{H}}
\newcommand{\cI}{\mathcal{I}}
\newcommand{\cL}{\mathcal{L}}

\newcommand{\cW}{\mathcal{W}}

\newcommand{\cY}{\mathcal{Y}}
\newcommand{\cZ}{\mathcal{Z}}

\newcommand{\fb}{\mathfrak{b}}
\newcommand{\fc}{\mathfrak{c}}
\newcommand{\fS}{\mathfrak{S}}
\newcommand{\fZ}{\mathfrak{Z}}

\newcommand{\ty}[1]{\mathsf{#1}}

\newcommand{\cC}{\mathscr{C}}
\newcommand{\cE}{\mathscr{E}}

\newcommand{\Alp}{\operatorname{Alp}}
\newcommand{\Aut}{\operatorname{Aut}}
\newcommand{\bl}{\operatorname{bl}}

\newcommand{\de}{\operatorname{def}}
\newcommand{\Diag}{\operatorname{Diag}}
\newcommand{\dz}{\operatorname{dz}}
\newcommand{\IBr}{\operatorname{IBr}}
\newcommand{\Ind}{\operatorname{Ind}}
\newcommand{\Irr}{\operatorname{Irr}}
\newcommand{\Lin}{\operatorname{Lin}}
\newcommand{\Out}{\operatorname{Out}}
\newcommand{\Res}{\operatorname{Res}}
\newcommand{\tr}{\operatorname{tr}}

\newcommand{\SC}{{\operatorname{sc}}}
\newcommand{\rdz}{\operatorname{rdz}}
\newcommand{\mrO}{\operatorname{O}}
\newcommand{\Z}{\operatorname{Z}}
\newcommand{\N}{\operatorname{N}}
\newcommand{\C}{\operatorname{C}}
\newcommand{\R}{\operatorname{R}}

\newcommand{\GL}{\operatorname{GL}}
\newcommand{\SL}{\operatorname{SL}}
\newcommand{\PSL}{\operatorname{PSL}}
\newcommand{\GU}{\operatorname{GU}}
\newcommand{\SU}{\operatorname{SU}}
\newcommand{\PSU}{\operatorname{PSU}}
\newcommand{\Sp}{\operatorname{Sp}}

\newcommand{\tvhi}{\widetilde{\vhi}}

\newcommand{\tchi}{\widetilde{\chi}}
\newcommand{\tpsi}{\widetilde{\psi}}
\newcommand{\hchi}{\widehat{\chi}}
\newcommand{\ttheta}{\widetilde{\theta}}

\newcommand{\tB}{{\widetilde{B}}}
\newcommand{\tC}{\widetilde{C}}

\newcommand{\tG}{\widetilde{G}}

\newcommand{\tN}{\widetilde{N}}
\newcommand{\tR}{\widetilde{R}}
\newcommand{\tZ}{\widetilde{Z}}

\newcommand{\tbZ}{\widetilde{\mathbf{Z}}}
\newcommand{\tbG}{{\widetilde{\mathbf{G}}}}
\newcommand{\tbL}{{\widetilde{\mathbf{L}}}}
\newcommand{\tbT}{\widetilde{\mathbf{T}}}
\newcommand{\barFF}{{\overline{\FF}}}

\newcommand{\RLG}{{\R_\bL^\bG}}

\def\embed{\hookrightarrow}

\let\al=\alpha
\let\eps=\epsilon
\let\ga=\gamma
\let\Ga=\Gamma
\let\ka=\kappa
\let\la=\lambda
\let\vOm=\varOmega
\let\vhi=\varphi
\let\ze=\zeta
\let\ti=\times

\theoremstyle{theorem}

\newtheorem{thm}{Theorem}[section]
\newtheorem{lem}[thm]{Lemma}
\newtheorem{prop}[thm]{Proposition}
\newtheorem{cor}[thm]{Corollary}
\newtheorem{amp}[thm]{Assumption}
\newtheorem{conj}[thm]{Conjecture}
\newtheorem{condition}[thm]{Condition}
\newtheorem{question}[thm]{Question}

\theoremstyle{definition}
\newtheorem{defn}[thm]{Definition}
\newtheorem{notation}[thm]{Notation}
\newtheorem{rmk}[thm]{Remark}

\numberwithin{equation}{section}

\begin{document}

%%%%%%%%%%%%%%%%%%%%%%%%%%%%%%%%%%%%%%%%%%%%%%%%%%%%%%%%%%%%
\title[Generic weights for finite reductive groups]{Generic weights for finite reductive groups}
%%%%%%%%%%%%%%%%%%%%%%%%%%%%%%%%%%%%%%%%%%%%%%%%%%%%%%%%%%%%

\author{Zhicheng Feng}
\address[Z. Feng]{Shenzhen International Center for Mathematics and
 Department of Mathematics, Southern University of Science and Technology,
 Shenzhen 518055, China}
\makeatletter
\email{fengzc@sustech.edu.cn}
\makeatother

\author{Gunter Malle}
\address[G. Malle]{FB Mathematik, RPTU Kaiserslautern, Postfach 3049,
 67653 Kaisers\-lautern, Germany}
\makeatletter
\email{malle@mathematik.uni-kl.de}
\makeatother

\author{Jiping Zhang}
\address[J. Zhang]{School of Mathematical Sciences, Peking University,
 Beijing 100871, China} 
\makeatletter
\email{jzhang@pku.edu.cn}
\makeatother

\thanks{The first and third authors gratefully acknowledge financial support
 by National grants (12431001, 12350710787 \& 2020YFE0204200).
 The second author gratefully acknowledges financial support by the DFG
 --Project-ID 286237555--TRR 195.}

\begin{abstract}
This paper is motivated by the study of Alperin's weight conjecture in the
representation theory of finite groups. We generalize the notion of
$e$-cuspidality in the $e$-Harish-Chandra theory of finite reductive groups,
and define generic weights in non-defining characteristic.
We show that the generic weights play an analogous role as the weights defined
by Alperin in the investigation of the inductive Alperin weight condition for
simple groups of Lie type at most good primes. We hope that our approach will
constitute a major step towards a proof of Alperin's weight conjecture.
\end{abstract}



\date{\today}

\maketitle

%\tableofcontents 

%\pagestyle{myheadings}
%\markboth{for personal use only}{preliminary}

%%%%%%%%%%%%%%%%%%%%%%%%%%%%%%%%%%%%%%%%%%%%%%%%%%%%%%%%%%%%%%%%%%%%%%%%%
\section{Introduction}

Harish-Chandra theory is a significant tool in Lie theory, such as in the
representation theory of Lie groups, Lie algebras and finite reductive groups.
For finite groups with a BN-pair, Harish-Chandra theory provides a way to
construct irreducible characters in non-defining characteristic in terms of the
representation theory of Hecke algebras. The generalized, so-called
$e$-Harish-Chandra theory is built by using Deligne--Lusztig induction instead
of Harish-Chandra induction, and it plays a fundamental role in modular
representation theory of finite reductive groups; for example, in the
classification of blocks of finite reductive groups (cf.\ 
\cite{CE94,CE99,KM13,KM15}).

Let $\bG$ be a connected reductive linear algebraic group with a Frobenius
endomorphism $F\colon\bG\to\bG$ endowing $\bG$ with an $\FF_q$-structure. 
We let $\ell$ be a prime not dividing $q$. Let $e_\ell(q)$ denote the
multiplicative order of~$q$ modulo~$\ell$ or~4 according as $\ell\ge 3$ or
$\ell=2$. In this paper, we give a generalization of $e$-cuspidality and
$e$-Jordan-cuspidality defined by Brou\'e--Malle--Michel \cite{BMM93} and
Kessar--Malle respectively \cite{KM15}. We classify the unipotent
$e$-generalized-cuspidal characters of $\bG^F$ for all odd primes~$\ell$ with
$e=e_\ell(q)$. For groups of type $\ty A$ and odd primes, the unipotent
$e$-generalized-cuspidal characters are labeled by hook partitions
(Proposition~\ref{lem:e-GC-type-A}), which have been widely highlighted in the
representation theory of both symmetric groups and groups of type $\ty A$.
It seems reasonable to expect our generalization here to have potential
applications in the modular representation theory of finite reductive groups.
For example, we expect a new partition for the irreducible characters of
relative Weyl groups.

The Alperin weight conjecture, announced by Alperin \cite{Al87} in 1986, is one
of the central problems in the modular representation theory of finite groups. 
It remains open to the present day, and perhaps the most promising approach is
to reduce the problem to simple groups. In 2011, Navarro and Tiep \cite{NT11}
achieved a reduction for the Alperin weight conjecture; they proved that if
every finite non-abelian simple group satisfies the so-called \emph{inductive
AW condition}, then the Alperin's weight conjecture holds for every finite
group. The verification of the inductive condition has been achieved for
alternating groups, sporadic groups, simple group of Lie type at their defining
characteristic, as well as simple groups of type $\ty A$. For the recent
developments around the inductive investigation of Alperin weight conjecture,
we refer to the survey paper \cite{FZ22} by the first and third authors.
Even though several families of simple groups have been proved to satisfy the
inductive AW condition, it is still a great challenge to complete the
verification, even for a given prime.

For the McKay conjecture, the second author \cite{Ma07} used the normalizers of
Sylow $e$-tori in place of normalizers of Sylow subgroups, and this approach has
already proved successful in the proof of the McKay conjecture; see the recent
paper \cite{CS24} of Cabanes--Sp\"ath and the references therein. 
A similar approach was also used by Rossi to study Dade's conjecture in
\cite{Ro23,Ro24}. For the
inductive Alperin weight condition, we often need to consider a plethora of
radical subgroups and analyze the representations of their normalizers.
This seems to be very difficult for several families of quasi-simple groups of
Lie type, especially for types $\ty D$ and $\ty E$. In this paper, we present
an approach to connect the inductive AW condition and the generalized
Harish-Chandra theory.
Using the $e$-Jordan-generalized-cuspidal pairs, we define new objects,
called \emph{generic weights}, for finite reductive groups in non-defining
characteristic. 
We show that in the verification of the inductive Alperin weight condition for
simple groups of Lie type and most good primes, the weights defined by Alperin
can be replaced by our generic weights. The advantage is that we may consider
less local subgroups, and use more generic generalized Harish-Chandra theory.
We hope that our approach will constitute a major step towards an eventual
proof of Alperin's weight conjecture, as with the similar approach in the proof
of the McKay conjecture.
Additionally, in the spirit of Alperin's weight conjecture and
$e$-Harish-Chandra theory, we introduce a problem
(Question~\ref{ques:partition-local}) concerning the correspondence of
characters at the level of relative Weyl groups, that partitions the
irreducible characters in terms of defect zero characters of some smaller
relative Weyl groups. We prove it for all unipotent blocks, as well as for
quasi-isolated blocks of exceptional groups.
\medskip

The paper is organized as follows.
After introducing some notation in Section \S\ref{sec:Preli}, we propose a
generalization of $e$-cuspidality, define the generic weights, and reduce the
determination of generic weights to quasi-isolated blocks inductively in
Section \S\ref{sec:gene-eJC}. In Section \S\ref{sec:wei-par}, we partition
weights in terms of the center of radical subgroups. Section \S\ref{sec:type-A}
is devoted to the study of the generic weights of groups of type $\ty A$, while
Section \S\ref{sec:IBAW} establishes the relation between weights and generic
weights, and give criteria for the inductive condition of the Alperin weight
conjecture, non-blockwise version or blockwise version.
In Section~\S\ref{sec:open-problem}, we propose a problem for correspondence of
characters on the level of relative Weyl groups, and demonstrate its validity
for all unipotent blocks, as well as for quasi-isolated blocks of exceptional
groups.
\bigskip

\noindent{\bf Acknowledgement:} The authors thank Damiano Rossi for his
 pertinent comments on an earlier version.

%%%%%%%%%%%%%%%%%%%%%%%%%%%%%%%%%%%%%%%%%%%%%%%%%%%%%%%%%%%%%%%%%%%%%%%%%
\section{Preliminaries}   \label{sec:Preli}

%%%%%%%%%%%%%%%%%%%%%%%%%%%%%%%%%%%%%
\subsection{General notation}
If a group $G$ acts on a set $X$, we let $G_x$ denote the stabilizer of
$x\in X$ in $G$, and analogously we denote the setwise stabilizer of
$X'\subseteq X$ in $G$ by $G_{X'}$.
If $H\le G$, then we denote by $X/\!\sim_H$ the set of $H$-orbits on $X$.
Moreover, if a group $G$ acts on two sets $X$, $Y$ and $x\in X$, $y\in Y$, we
denote by $G_{x,y}$ the stabilizer of $y$ in $G_x$. For a positive integer $n$,
we denote the symmetric group on $n$ symbols by $\fS_n$.

Suppose that $G$ is a finite group and $H\le G$.
We denote the restriction of $\chi\in\Irr(G)$ to $H$ by $\Res_H^G(\chi)$, and
$\Ind^G_H(\theta)$ denotes the character induced from $\theta\in\Irr(H)$ to $G$.
As usual, for $\theta\in\Irr(H)$ the set of irreducible constituents of
$\Ind_H^G(\theta)$ is denoted by $\Irr(G\mid\theta)$, while $\Irr(H\mid\chi)$
denotes the set of irreducible constituents of $\Res_H^G(\chi)$ for
$\chi\in\Irr(G)$. For a subset $\cH\subseteq\Irr(H)$, we define
\[\Irr(G\mid\cH)=\bigcup_{\theta\in\cH}\Irr(G\mid\theta)\]
and for a subset $\cG\subseteq\Irr(G)$, we define
\[\Irr(H\mid\cG)=\bigcup_{\chi\in\cG}\Irr(H\mid\chi).\]
Additionally, for $N \unlhd G$, we sometimes identify the characters of $G/N$
with the characters of $G$ whose kernel contains $N$.

Let $\ell$ be a prime number. Throughout, all modular representations considered
are with respect to $\ell$. For $\chi\in\Irr(G)$, the $\ell$-block of~$G$
containing $\chi$ is denoted by $\bl(\chi)$, which is also denoted by
$\bl_G(\chi)$ where we add a subscript to indicate the ambient group $G$.
If $b$ is a union of blocks of~$G$, then we write
$\Irr(b)=\bigcup_{B\in b}\Irr(B)$. For a block $b$ of a subgroup $H\le G$ we
denote by $b^G$ the induced block of $G$, when it is defined.

Denote by $\dz(G)$ the set of irreducible characters of $G$ of ($\ell$-)defect
zero. For $N\unlhd G$ and $\theta\in\Irr(N)$ we set
\[\rdz(G\mid\theta):=
  \{\,\chi\in\Irr(G\mid\theta)\mid\chi(1)_\ell/\theta(1)_\ell=|G/N|_\ell \,\}.\]
If moreover $\theta\in\dz(N)$, then $\rdz(G\mid\theta)\subseteq\dz(G)$, and
then we also write $\dz(G\mid\theta)$ for $\rdz(G\mid\theta)$.

Let $\Lin(G)$ denote the set of linear characters of $G$, which can be seen
as a multiplicative group. Then $\Lin(G)$ acts on $\Irr(G)$ by multiplication.
The Hall $\ell'$-subgroup of $\Lin(G)$ is denoted $\Lin_{\ell'}(G)$.

Denote by $\mrO_\ell(G)$ the largest normal $\ell$-subgroup of $G$. Similarly,
$\mrO_{\ell'}(G)$ denotes the largest normal $\ell'$-subgroup of $G$.
If $A$ is an abelian group, then we also write $A_\ell$ for its Sylow
$\ell$-subgroup and write $A_{\ell'}$ for its Hall $\ell'$-subgroup;
note that $A_\ell=\mrO_\ell(A)$ and $A_{\ell'}=\mrO_{\ell'}(A)$.

%%%%%%%%%%%%%%%%%%%%%%%%%%%%%%%%%%%%%
\subsection{Radical subgroups and weights}   \label{subsec:preli-weights}

Let $G$ be a finite group. We denote by $\Re^0(G)$ the set of radical
$\ell$-subgroups of $G$, and write $\Re(G)=\Re^0(G)/\!\sim_G$.

A \emph{weight} of $G$ is a pair $(R,\vhi)$, where $R$ is a (possibly trivial)
$\ell$-subgroup of $G$ and $\vhi\in\dz(\N_G(R)/R)$.
For an $\ell$-subgroup $R$ of $G$, if there exists a weight $(R,\vhi)$ of $G$,
then we say that $R$ is a \emph{weight ($\ell$-)subgroup} of $G$.
Denote by $\Re_w^0(G)$ the set of weight $\ell$-subgroups of $G$.
Then $\Re_w^0(G)\subseteq \Re^0(G)$. Let $\Re_w(G)=\Re_w^0(G)/\!\sim_G$.
Let $\Alp^0(G)$ denote the set of weights of $G$. 
The $G$-orbit of $(R,\vhi)$ is denoted by $\overline{(R,\vhi)}$ and we define
$\Alp(G)=\Alp^0(G)/\!\sim_G$. Sometimes we also write $\overline{(R,\vhi)}$
simply as $(R, \vhi)$ when no confusion can arise.
For $\nu\in\Lin_{\ell'}(\Z(G))$, we denote by $\Alp^0(G\mid \nu)$ the set of
weights $(R,\vhi)$ of $G$ with $\vhi\in\Irr(\N_G(R)\mid\nu)$ and write
$\Alp(G\mid\nu)=\Alp^0(G\mid \nu)/\!\sim_G$.

The group $\Lin_{\ell'}(G)$ acts on $\Alp^0(G)$ by $\mu.(R,\vhi)=(R,\mu'\vhi)$
where $\mu'$ is the restriction of $\mu\in\Lin_{\ell'}(G)$ to $\N_G(R)$
(sometimes we also write $\mu$ for $\mu'$ when no confusion can arise); see
\cite[Lemma~2.4]{BS22}.
This induces an action of $\Lin_{\ell'}(G)$ on $\Alp(G)$.

\begin{lem}   \label{lem:cenprod-rad-corr}
 Let $G$ be a finite group with $G=HZ$, where $H\le G$, $Z\le \Z(G)$.
 Then $S\mapsto S Z_\ell$ gives a bijection from $\Re^0(H)$ to $\Re^0(G)$ with
 inverse $R\mapsto R\cap H$. In addition, this induces bijections
 $\Re(H)\to\Re(G)$ and $\Re_w(H)\to\Re_w(G)$.
\end{lem}

\begin{proof}
Let $Z_0=Z\cap H$.
Then $(Z_0)_\ell$ is contained in any radical $\ell$-subgroup of $G$.
Denote by $\pi:G\to G/Z_0$ the canonical epimorphism.
Then $R\mapsto RZ_0/Z_0$ gives a bijection $\Re(G)\to\Re(G/Z_0)$ with inverse
$\overline R\mapsto \mrO_\ell(\pi^{-1}(\overline R))$.
Now $G/Z_0\cong H/Z_0\ti Z/Z_0$.
So $\overline S\mapsto \overline S\ti (Z/Z_0)_\ell$ gives a bijection
$\Re^0(H/Z_0)\to\Re^0(G/Z_0)$ with inverse
$\overline R\mapsto \overline R\cap (H/Z_0)$.
Thus the assertion holds.
\end{proof}

Each weight may be assigned to a unique block. Let $B$ be a block of $G$. 
A weight $(R,\vhi)$ of $G$ is said to be a $B$-weight if
$\bl_{\N_G(R)}(\vhi)^G=B$. We denote the set of $B$-weights by $\Alp^0(B)$.
Let $\Alp(B)=\Alp^0(B)/\!\sim_G$. If $b$ is a union of blocks of $G$, then we
define $\Alp(b)=\bigcup_{B\in b}\Alp(B)$.

In \cite{BS22}, Brough and Sp\"ath defined a relationship ``covering" between
weights of a finite group and its normal subgroups.
Let $G$ be a normal subgroup of a finite group $\tG$.
If $(R,\vhi)$ is a weight of $G$, then we write $\Alp^0(\tG\mid(R,\vhi))$ for
the set of those $(\tR,\tvhi)\in\Alp^0(\tG)$ covering $(R,\vhi)$.
If $(\tR,\tvhi)$ is a weight of $\tG$, then we write $\Alp^0(G\mid(\tR,\tvhi))$
for the set of those $(R,\vhi)\in\Alp^0(G)$ covered by $(\tR,\tvhi)$.

For $(\tR,\tvhi)\in\Alp^0(\tG)$ and $(R,\vhi)\in\Alp^0(G)$, we say that
$\overline{(\tR,\tvhi)}$ \emph{covers} $\overline{(R,\vhi)}$ if $(\tR,\tvhi)$
covers $(R^g,\vhi^g)$ for some $g\in\tG$. If $(R,\vhi)$ is a weight of $G$ we
write $\Alp(\tG\mid\overline{(R,\vhi)})$ for the set of those
$\overline{(\tR,\tvhi)}\in\Alp(\tG)$ covering $\overline{(R,\vhi)}$.
If $(\tR,\tvhi)$ is a weight of $\tG$, then we write
$\Alp(G\mid\overline{(\tR,\tvhi)})$ for the set of those
$\overline{(R,\vhi)}\in\Alp(G)$ covered by $\overline{(\tR,\tvhi)}$.
For a subset $\cA\subseteq\Alp(G)$ we define
\[\Alp(\tG\mid\cA)
  =\bigcup_{\overline{(R,\vhi)}\in\cA} \Alp(\tG\mid\overline{(R,\vhi)})\]
and for a subset $\widetilde\cA\subseteq\Alp(\tG)$, we define
\[\Alp(G\mid\widetilde\cA)=\bigcup_{\overline{(\tR,\tvhi)}\in\widetilde\cA}
  \Alp(G\mid\overline{(\tR,\tvhi)}). \]

%%%%%%%%%%%%%%%%%%%%%%%%%%%%%%%%%%%%%%%%%%%%%%%%%%%%%%%%%%%%%%%%%%%%%%%%%
%%%%%%%%%%%%%%%%%%%%%%%%%%%%%%%%%%%%%%%%%%%%%%%%%%%%%%%%%%%%%%%%%%%%%%%%%
\section{A generalization of cuspidality and generic weights}   \label{sec:gene-eJC}

Let $\bG$ be a connected reductive group over $\barFF_p$ for a prime~$p$
and let $F\colon\bG\to \bG$ be a Frobenius endomorphism defining an
$\FF_q$-structure on $\bG$, where~$q$ is a power of~$p$.
Write $\cZ(\bG):=\Z(\bG)/\Z^\circ(\bG)$. Denote by $\cZ(\bG)_F$ the largest
quotient of $\cZ(\bG)$ on which $F$ acts trivially.
Let $\bG^*$ be the Langlands dual of $\bG$, whose root datum can be obtained
from that of $\bG$ by exchanging character group and cocharacter group, as well
as roots and coroots.
We denote the corresponding Frobenius endomorphism of $\bG^*$ also by $F$ for
simplicity (see \cite[\S1.5]{GM20}).
Let $\ell$ be a prime different from~$p$ throughout.

%%%%%%%%%%%%%%%%%%%%%%%%%%%%%%%%%%%%%
\subsection{Embeddings between reductive groups}

\begin{defn}\label{defn-embedding}
Let $\bG$, $\tbG$ be connected reductive groups over $\barFF_p$ with
Frobenius endomorphisms $F\colon\bG\to\bG$, $\widetilde F\colon\tbG\to\tbG$.
Suppose that $i\colon\bG\to\tbG$ is a homomorphism of algebraic groups such
that $i\circ F=\widetilde F\circ i$.
\begin{enumerate}[\rm(a)]
 \item We say $i$ is a \emph{weakly regular embedding} if $i$ is an isomorphism
  of $\bG$ with a closed subgroup of $\tbG$ and $[\tbG,\tbG]=[i(\bG),i(\bG)]$.
 \item We say $i$ is an \emph{$\ell$-regular embedding} if $i$ is a
  weakly regular embedding and $\ell\nmid |\cZ(\tbG)_{\widetilde F}|$.	
 \item Following Lusztig \cite[\S7]{Lu88}, we say $i$ is a \emph{regular
  embedding} if $i$ is a weakly regular embedding and $\Z(\tbG)$ is connected
 (i.e., $\cZ(\tbG)=1$).
\end{enumerate} 
\end{defn}

If $i$ is a weakly regular embedding, we identify~$\bG$ with~$i(\bG)$ and
denote~$\widetilde F$ briefly by~$F$ since~$\widetilde F$ can be viewed as an
extension of~$F$.

\begin{lem}   \label{lem:diga-ell'}
 \begin{enumerate}[\rm(a)]
	\item Let $\bG\embed\tbG$ be a weakly regular embedding. Then
     $\tbG^F/\bG^F\Z(\tbG)^F$ is isomorphic to a subgroup $X$ of $\cZ(\bG)_F$
     such that $\cZ(\bG)_F/X\cong \cZ(\tbG)_F$.	
	\item Let $\bG\embed\tbG$ be an $\ell$-regular embedding. Then
     $(\tbG^F/\bG^F\Z(\tbG)^F)_\ell\cong(\cZ(\bG)_F)_\ell$.
 \end{enumerate}
\end{lem}
	
\begin{proof}
We prove (a), from which (b) follows directly.	
Let $\tbG\embed \tbG'$ be a regular embedding. Then $\bG\embed\tbG\embed\tbG'$
is also a regular embedding. Write $\tbZ:=\Z(\tbG)$ and $\tbZ':=\Z(\tbG')$.
Then $\tbZ^F=\tbZ'^F\cap\tbG^F$. From this,
\[\tbG^F\tbZ'^F/\bG^F\tbZ'^F\cong
  \tbG^F/(\tbG^F\cap\bG^F\tbZ'^F)=\tbG^F/\bG^F\tbZ^F.\]
By \cite[Rem.~1.7.6]{GM20}, one has $\tbG'^F/\bG^F\tbZ'^F\cong\cZ(\bG)_F$ and
$\tbG'^F/\tbG^F\tbZ'^F\cong\cZ(\tbG)_F$.
Then $\tbG^F/\bG^F\tbZ^F$ is isomorphic to a subgroup of $\cZ(\bG)_F$ with
$\cZ(\bG)_F/(\tbG^F/\bG^F\tbZ^F)\cong \cZ(\tbG)_F$.	
\end{proof}

\begin{lem}   \label{lem:regular-embedding}
 Let $\bG\embed\bG_1$ and $\bG\embed \bG_2$ be weakly regular embeddings.
 Then there exists a connected reductive group~$\tbG$ and regular embeddings
 $\bG\embed\tbG$, $\bG_1\embed\tbG$ and $\bG_2\embed\tbG$.
\end{lem}

\begin{proof}
Let $\bG_i\embed \tbG_i$, $i=1,2$, be regular embeddings. Then
$\bG\embed \tbG_i$ are regular embeddings.
By a result of Asai (cf.\ \cite[Lemma~7.1]{Lu88}), there exists a connected
reductive group~$\tbG$ with regular embeddings $\tbG_i\embed \tbG$, $i=1,2$.
Then the assertion holds for this~$\tbG$.
\end{proof}

\begin{lem}   \label{relative-regular-Levi}
 Let $\bG\embed\tbG$ be an $\ell$-regular embedding. Then for any
 $F$-stable Levi subgroup $\bL$ of $\bG$, $\bL\embed \tbL:=\bL\Z^\circ(\tbG)$
 is an $\ell$-regular embedding.
\end{lem}

\begin{proof}
Note that $\tbL$ is a Levi subgroup of $\tbG$.
Now, by \cite[Prop.~2.4]{GH91}, for any Levi subgroup $\tbL$ of $\tbG$ we have
that $\cZ(\tbL)_F$ is a factor group of $\cZ(\tbG)_F$. Since $\bG\embed\tbG$ is
$\ell$-regular, $\cZ(\tbG)_F$ has order prime to $\ell$, and thus the same is
true for $\cZ(\tbL)_F$, whence $\bL\embed \tbL$ is $\ell$-regular.
\end{proof}

%%%%%%%%%%%%%%%%%%%%%%%%%%%%%%%%%%%%%
\subsection{Blocks of finite reductive groups}
For our notation and basic facts about representation theory of finite
reductive groups we refer the reader to \cite{GM20}.
To distribute the irreducible characters of $\bG^F$ into $\ell$-blocks, we
define for each semisimple $\ell'$-element $s$ of ${\bG^*}^F$, the set
$\cE_\ell(\bG^F,s)$, which is the union of Lusztig series $\cE(\bG^F,st)$ where
$t$ runs through the semisimple $\ell$-elements of ${\bG^*}^F$ commuting
with~$s$.

By a theorem of Brou\'e--Michel \cite[Thm.~9.12]{CE04}, $\cE_\ell(\bG^F,s)$ is
a union of $\ell$-blocks of $\bG^F$ for every semisimple $\ell'$-element $s$ of
${\bG^*}^F$. We denote by $\cE(\bG^F,\ell')$ the union of Lusztig series
$\cE(\bG^F,s)$ where $s$ runs through the semisimple $\ell'$-element of
${\bG^*}^F$.

For a positive integer $e$, we denote by $\phi_e$ the $e$-th cyclotomic
polynomial. We will make use of the terminology of Sylow $e$-theory (see for
instance \cite[\S3.5]{GM20} or \cite[\S25]{MT11}). For $\bT$ an $F$-stable
maximal torus, $\bT_{\phi_e}$ denotes its Sylow $e$-torus.

Let $E\subseteq \ZZ_{\ge 1}$. Recall that an \emph{$E$-torus} of $\bG$ is an
$F$-stable torus whose polynomial order is a product of cyclotomic polynomials
in $\{\phi_e\mid e\in E\}$ and an \emph{$E$-split Levi subgroup} of $\bG$ is
the centralizer of an $E$-torus of $\bG$.
We say that an irreducible character $\chi\in\Irr(\bG^F)$ is
\emph{$E$-cuspidal} if $^*\R^\bG_{\bL\subseteq \bP}(\chi)=0$ for all proper
$E$-split Levi subgroups $\bL$ of $\bG$ and any parabolic subgroup $\bP$ of
$\bG$ containing $\bL$ as Levi complement. If $\bL\le\bG$ is $E$-split
and $\la\in\Irr(\bL^F)$ is $E$-cuspidal, then
$(\bL,\la)$ is called an \emph{$E$-cuspidal pair of $\bG$}.
If $E=\{ e\}$, then we also say \emph{$e$-cuspidal} for $E$-cuspidal.

Let~$e$ be a positive integer and $s\in{\bG^*}^F$ be semisimple. As in
\cite[Def.~2.1]{KM15}, we say $\chi\in\cE(\bG^F,s)$ is
\emph{$e$-Jordan-cuspidal} if
$\Z^\circ(\C^\circ_{\bG^*}(s))_{\phi_e}=\Z^\circ(\bG^*)_{\phi_e}$ and $\chi$
corresponds under Lusztig's Jordan decomposition (cf.\ \cite{Lu88}) to the
$\C_{\bG^*}(s)^F$-orbit of a unipotent $e$-cuspidal character of
$\C^\circ_{\bG^*}(s)^F$.
If $\bL$ is an $e$-split Levi subgroup of $\bG$ and $\la\in\Irr(\bL^F)$ is
$e$-Jordan-cuspidal, then $(\bL,\la)$ is called an \emph{$e$-Jordan-cuspidal
pair} of $\bG$.
The blocks of $\bG^F$ for good primes were classified in \cite{CE99,KM15} in
terms of $e$-Jordan-cuspidal pairs.

%%%%%%%%%%%%%%%%%%%%%%%%%%%%%%%%%%%%%
\subsection{$e$-Jordan-generalized-cuspidal pairs}
Now we propose a generalization of $e$-(Jordan-) cuspi\-dality. Set
$E_{e,\ell}=\{\,e\ell^i\mid i=0,1,2,\ldots \,\}$ or $\{1,2,4,8,\ldots \}$
according as $\ell\ge3$ or $\ell=2$.

\begin{defn}\label{defn:eJGC}
 Let~$e$ be a positive integer and let $E:=E_{e,\ell}$.
 \begin{enumerate}[\rm(a)]
  \item Let $\chi\in\Irr(\bG^F)$. We say $\chi$ is
   \emph{$(e,\ell)$-generalized-cuspidal} (\emph{$(e,\ell)$-GC}) if $\langle\chi,\R_{\bL\subseteq\bP}^{\bG}(\la)\rangle\ne0$ for some $E$-cuspidal pair $(\bL,\la)$ of $\bG$ with
   $\Z^\circ(\bL)_{\phi_e}=\Z^\circ(\bG)_{\phi_e}$ (i.e.,
   $\Z^\circ(\bL)_{\phi_e}\subseteq\Z(\bG)$) and some parabolic subgroup $\bP$
   of $\bG$ containing $\bL$ as a Levi subgroup.
  \item Let $s\in{\bG^*}^F$ be semisimple. We say $\chi\in\cE(\bG^F,s)$ is
   \emph{$(e,\ell)$-Jordan-generalized-cuspidal ($(e,\ell)$-JGC)} if
   \begin{itemize}
	\item $\Z^\circ(\C^\circ_{\bG^*}(s))_{\phi_e}=\Z^\circ(\bG^*)_{\phi_e}$, and
	\item $\chi$ corresponds under Jordan decomposition to the
     $\C_{\bG^*}(s)^F$-orbit of a unipotent $(e,\ell)$-GC character of
     $\C^\circ_{\bG^*}(s)^F$.
   \end{itemize}
  \item If $\bL$ is an $e$-split Levi subgroup of $\bG$ and $\la\in\Irr(\bL^F)$
   is $(e,\ell)$-GC (resp. $(e,\ell)$-JGC), then $(\bL,\la)$
   is called an \emph{$(e,\ell)$-GC pair} (resp. \emph{$(e,\ell)$-JGC pair})
   of~$\bG$.
 \end{enumerate}
\end{defn}	

\begin{lem}\label{lem:e-J-eJGC}
 The $e$-(Jordan-)cuspidal characters are $(e,\ell)$-(J)GC for any prime $\ell$.
\end{lem}

\begin{proof}
By definition, it suffices to show that $e$-cuspidal characters are
$(e,\ell)$-GC. Let $\chi$ be an $e$-cuspidal character of $\bG^F$.
Assume that $\chi$ occurs in $\RLG(\la)$ for some proper $E$-split Levi
subgroup $\bL$ of $\bG$ and some $E$-cuspidal $\la\in\Irr(\bL^F)$.
If $\Z^\circ(\bL)_{\phi_e}\not\subseteq\Z^\circ(\bG)_{\phi_e}$ then $\chi$ also
occurs in $\R_\bH^\bG(\la')$ with $\bH=\C_\bG(\Z^\circ(\bL)_{\phi_e})$, an
$e$-split proper Levi subgroup of $\bG$, and some $\la'\in\Irr(\bH^F)$, in
contradiction to $\chi$ being $e$-cuspidal. So we have
$\Z^\circ(\bL)_{\phi_e}=\Z^\circ(\bG)_{\phi_e}$, and then $\chi$ is
$(e,\ell)$-GC by definition.
\end{proof}

We will simply write $e$-GC (resp.\ $e$-JGC) for $(e,\ell)$-GC (resp.\ 
$(e,\ell)$-JGC) when $\ell$ is clear from the context. 

\begin{rmk}   \label{rmk-construction-e-GC}
	Let $(\bL_0,\la_0)$ be an $E$-cuspidal pair of $\bG$,
	$\bL=\C^\circ_\bG(\Z^\circ(\bL_0)_{\phi_e})$ and $\la$ be an irreducible
	constituent of $\R_{\bL_0}^\bL(\la_0)$. Then $(\bL,\la)$ is an $e$-GC
	pair of $\bG$. Moreover, all $e$-GC pairs of $\bG$ can be obtained in this
	way.	
\end{rmk}

\begin{lem}   \label{lem:Clifford-eJGC}
 Let $\bG\embed\tbG$ be a weakly regular embedding, and $\tbL\le\tbG$ be an
 $F$-stable Levi subgroup. Let $\widetilde\la\in\Irr(\tbL^F)$,
 $\bL:=\tbL\cap \bG$ and $\la\in\Irr(\bL^F\mid\widetilde\la)$.
 Then $(\tbL,\widetilde\la)$ is an $e$-JGC pair of $\tbG$ if and only if
 $(\bL,\la)$ is an $e$-JGC pair of $\bG$.
\end{lem}

\begin{proof}
Let $\tbG\embed\tbG_1$ be a regular embedding. Then the composition
$\bG\embed\tbG_1$ is also a regular embedding. It is shown in the proof of
\cite[Prop.~1.10]{CE99} that $e$-Jordan-cuspidality is preserved under any
regular embedding, hence $e$-Jordan-cuspidality is equivalent for $\bG$ and
$\tbG_1$, as well as for $\tbG$ and $\tbG_1$, hence also for $\bG$ and $\tbG$.
The very same argument also works for $e$-generalized-Jordan-cuspidality.
\end{proof}

Recall that $e_\ell(q)$ denotes the multiplicative order of~$q$ modulo~$\ell$
or~4 according as $\ell\ge 3$ or $\ell=2$. Then 
$E_{e_\ell(q),\ell}=\{d\in\ZZ_{\ge1}\mid \ell\ \text{ divides }\ \phi_d(q)\}$.
As formulated in \cite[1.11]{CE99} for $e$-cuspidality and $e$-Jordan
cuspidality it seems reasonable to expect the following (see
Corollary~\ref{cor:e-JGC} below for some evidence):

\begin{conj}   \label{conj:e-JGC}
 Let $\bG$ be connected reductive with Frobenius map $F$, let $\ell$ be a
 prime and $e:=e_\ell(q)$. Then $\chi\in\Irr(\bG^F)$ is $(e,\ell)$-JGC if and
 only if $\chi$ is $(e,\ell)$-GC.
\end{conj}

\begin{lem}   \label{lem:odd-eGC}
 Assume that $\ell$ is odd, good for $\bG$, and $\ell>3$ if $\bG^F$ has a
 component of type $^3\ty{D}_4$. Set $e:=e_\ell(q)$.
 If $\ell\nmid|\Z(\bG_\SC)^F|$ then $\chi\in\Irr(\bG^F)$ is $e$-GC if and only
 if $\chi$ is $e$-cuspidal.
\end{lem}

\begin{proof}
One direction is in Lemma~\ref{lem:e-J-eJGC}. Now by \cite[Thm.~22.2]{CE04},
under our assumptions every proper $E:=E_{e,\ell}$-split Levi subgroup of $\bG$
is contained in a proper $e$-split Levi subgroup of $\bG$, and thus
$\chi\in\Irr(\bG^F)$ is $E$-cuspidal if and only if it is $e$-cuspidal.
	
If $\bL$ is an $F$-stable Levi subgroup of $\bG$ such that
$\Z^\circ(\bL)_{\phi_e}\subseteq \Z^\circ(\bG)$, then
$\Z^\circ(\bL)_{\phi_E}\subseteq \Z^\circ(\bG)$ by \cite[Lemma~22.3]{CE04}.
From this, if $(\bL,\la)$ is an $E$-cuspidal pair of $\bG$ with
$\Z^\circ(\bL)_{\phi_e}\subseteq\Z^\circ(\bG)$, then
$\bL=\C_\bG(\Z^\circ(\bL)_{\phi_E})=\bG$. So $\chi$ is $e$-GC if and only if
$\chi$ is $E$-cuspidal, and thus if and only if $\chi$ is $e$-cuspidal.
\end{proof}	

\begin{lem}   \label{cor:abel-sylow}
 Suppose that $\bG$ is simple and $\bG^F$ has an abelian Sylow $\ell$-subgroup.
 Let $\bL$ be an $e$-split Levi subgroup of $\bG$, for $e:=e_\ell(q)$, and
 $\chi\in\cE(\bL^F,\ell')$. Then $\chi$ is $e$-cuspidal if and only if $\chi$
 is $e$-Jordan-cuspidal, if and only if $\chi$ is $e$-GC, if and only if $\chi$
 is $e$-JGC.
\end{lem}

\begin{proof}
By \cite[\S 2.1]{Ma14}, $e$ is the unique positive integer such that
$\ell\mid \phi_e(q)$ and $\phi_e$ divides the order polynomial of $(\bG,F)$,
and~$\ell$ is odd, good for $\bG$, $\ell>3$ if $(\bG,F)$ has type $^3\ty{D}_4$,
and does not divide the orders of $\cZ(\bG)^F$ and $\cZ(\bG^*)^F$.
Thus the $e$-tori and $E$-tori of $\bG^F$ coincide. By definition, the $e$-GC
characters are just the $e$-cuspidal characters, and hence the $e$-JGC
characters are the $e$-Jordan-cuspidal characters. Moreover, by
\cite[Thm.~4.2 and Rem.~5.2]{CE99}, $e$-Jordan-cuspidality and $e$-cuspidality
agree, which completes the proof.
\end{proof}

\begin{prop}   \label{lem:e-GC-type-A}
 Let $\bG=\SL_n(\overline \FF_q)$ (with $n\ge 2$) and $F\colon\bG\to\bG$ be a
 Frobenius endomorphism such that $\bG^F=\SL_n(\eps q)$ with $\eps\in\{\pm1\}$.
 Let $e:=e_\ell(q)$. Then $\bG^F$ possesses a unipotent $e$-GC character that
 is not $e$-cuspidal if and only if one of the following holds.
 \begin{enumerate}[\rm(a)]
  \item $\ell\mid (q-\eps)$ when $\ell$ is odd, respectively $4\mid(q-\eps)$
   when $\ell=2$, and $n=\ell^k$ for some integer~$k\ge1$, in which case the
   unipotent $e$-GC characters of $\bG^F$ are those parameterized by the
   hook partitions $(n), (n-1,1), (n-2,1^2), \ldots, (1^n)$ of $n$.
  \item $\ell=2$ and $4\mid(q+\eps)$, in which case all unipotent characters of
   $\bG^F$ are $e$-GC.
 \end{enumerate} 
\end{prop}	

\begin{proof}
Let $\tbG=\GL_n(\barFF_q)$ so that $\bG\embed\tbG$ is a regular embedding. Let
$\chi$ be a unipotent character of $\bG^F$, and let $\tchi$ be the unipotent
character of $\tbG^F$ with $\chi=\Res^{\tbG^F}_{\bG^F}(\tchi)$.
By Lemma~\ref{lem:Clifford-eJGC} and the proof of \cite[Prop.~1.10]{CE99},
$\chi$ is $e$-GC (resp.\ $e$-cuspidal) if and only if $\tchi$ is $e$-GC
(resp.\ $e$-cuspidal). For the classification of unipotent $e$-cuspidal
characters of $\tbG^F$, see \cite[\S4.3]{GM20}.

First let $\ell$ be odd. If $\ell\nmid(q-\eps)$, then by
Lemma~\ref{lem:odd-eGC} unipotent $e$-GC characters of $\bG^F$ are
$e$-cuspidal. Now we assume $\ell\mid(q-\eps)$, so $e=\frac{3-\eps}{2}$. In
this case, $\bG^F$ has no unipotent $e$-cuspidal character. In particular,
if $\tchi$ is $e$-GC, there exists a \emph{proper} $E$-split Levi
subgroup~$\tbL$ of~$\tbG$ with
$\Z(\tbL)_{\phi_e}\subseteq \Z(\tbG)$ and a unipotent $E$-cuspidal character
$\la$ of $\tbL^F$ such that $\tchi$ is an irreducible constituent of
$\R_\tbL^\tbG(\la)$. Now the order polynomial of $\Z(\tbG)$ is $\phi_e$, and
therefore $\Z(\tbL)_{\phi_e}\subseteq \Z(\tbG)$ implies that
$\tbL^F\cong \GL_m((\eps q)^{\ell^k})$ for some positive integers~$k$ and~$m$
with $n=m\ell^k$. Then $(\tbL,\la)$ is an $e\ell^k$-cuspidal pair of $\tbG$ and
thus $m=1$ and $\la=1_{\tbL^F}$. In particular, $n=\ell^k$.
So by \cite[\S4.3]{GM20} the irreducible constituents of $\R^\tbG_\tbL(\la)$
are parameterized by the partitions of $n$ that have an $n$-hook. This
completes the proof of necessity. The sufficiency is clear by construction.

Now let $\ell=2$. If $4\mid(q-\eps)$, then the above proof also applies.
Suppose that $4\mid(q+\eps)$. Then $e=2$ or 1 according as $\eps=1$ or $-1$. 
Any Sylow $\phi_e$-torus $\tbT$ of $\tbG$, with $|\tbT^F|=(q-\eps)^n$, is a
$d$-split Levi subgroup where $d=1$ or $2$ according as $\eps=1$ or $-1$, and
$\Z(\tbT)_{\phi_e}=\Z(\tbG)_{\phi_e}$.
So $(\tbT,1_{\tbT^F})$ is an $E$-cuspidal pair of $\tbG$ and every irreducible
constituent of $\R_{\tbT}^{\tbG}(1_{\tbT^F})$, and thus every unipotent
character of $\tbG^F$, is $e$-GC. This completes the proof.
\end{proof}

\begin{lem}   \label{lem:noncus-GC}
 Let $e:=e_\ell(q)$. If $\bG^F$ possesses a unipotent $e$-GC character that is
 not $e$-cuspidal, then one of the following holds.
 \begin{enumerate}[\rm(1)]
  \item $\ell$ is bad for $\bG$;
  \item $\bG$ has a component of type $^3\ty D_4$ and $\ell=3$;
  \item $\bG$ has a component of type $\ty{A}_n(\eps q^m)$ with $n+1=\ell^k>1$,
   and $\ell|(q^m-\eps)$ when $\ell\ge 3$, respectively $4\mid(q^m-\eps)$ when
   $\ell=2$; or
  \item $\bG$ has a component of type $\ty{A}_n(\eps q^m)$, $\ell=2$ and
   $4\mid(q^m+\eps)$.
 \end{enumerate}
\end{lem}

\begin{proof}
The inclusion $[\bG,\bG]\le\bG$ clearly is a weakly regular embedding, so by
Lemma~\ref{lem:Clifford-eJGC} we may assume $\bG$ is semisimple. Passing to a
simply connected covering, we reduce to the case that $\bG$ is simple and of simply connected type, in which
case the claim follows by combing Lemma~\ref{lem:odd-eGC} and
Proposition~\ref{lem:e-GC-type-A}.
\end{proof}

For completeness we classify the unipotent $e$-GC characters for bad primes
$\ell\ge3$.

\begin{prop}
 Let $\bG$ be simple with a Frobenius map $F$, $\ell\ge3$ a bad prime for~$\bG$,
 or $\ell=3$ and $\bG^F={}^3\ty{D}_4(q)$,
 and $e:=e_\ell(q)$. Then the $e$-GC unipotent characters of $\bG^F$ that are
 not $e$-cuspidal lie in the $e\ell^i$-Harish-Chandra series as described in
 Table~\ref{tab:e-GC}, up to Ennola duality. In the table, $\bT_i$, $\bT_i'$
 denote suitable $i$-tori of $\bG$.
\end{prop}

\begin{table}[htbp]
\caption{Unipotent $e$-GC characters for bad primes $\ell\ge3$}   \label{tab:e-GC}
$$\begin{array}{l|c|lcc}
 \bG^F& (\ell,e)& \text{$e\ell^i$-Harish-Chandra-series}\\
\hline
 \ty{G}_2(q),{}^3\ty{D}_4(q),\ty{F}_4(q)& (3,1)& (\bT_3,1)\\
% ^3\ty{D}_4(q)& (3,1)& (\bT_3,1)\\
% F_4(q)& (3.1)& (\bT_3,1) \\
 \ty{E}_6(q)& (3,1)& (\bT_3,1),(\bT_9,1),(\bT_3'.^3\ty{D}_4,{}^3\ty{D}_4[-1]) \\
 ^2\ty{E}_6(q)& (3,1)&(\bT_3\bT_6,1) \\
% \ty{E}_7(q)& (3.1)& (\bT_3,1)\\
 \ty{E}_8(q)& (3,1)& (\bT_3,1),(\bT_3'.^3\ty{D}_4,{}^3\ty{D}_4[-1])\\
 \ty{E}_8(q)& (5,1)& (\bT_5,1)\\
 \ty{E}_8(q)& (5,4)& (\bT_{20},1)
\end{array}$$
\end{table}

\begin{proof}
Since all classical groups are good for $\ell\ge3$ and we only consider
Frobenius endomorphisms, $\bG^F$ is one of $\ty{G}_2(q)$, $^3\ty{D}_4(q)$,
$\ty{F}_4(q)$, $^{(2)}\ty{E}_n(q)$ and $\ell=3$, or $\bG$ is of type $\ty{E}_8$
and $\ell=5$. If $\ell=3$ then $e\in\{1,2\}$ and up to Ennola duality we may
assume $e=1$. We refer to \cite[Tab.~3.3]{GM20} for a list of $d$-split Levi
subgroups of exceptional type groups. The only $E$-split Levi subgroups $\bL$
of $\bG$ of type $\ty G_2$ that are not 1-split are the Sylow 3-tori, whose only
($E$-cuspidal) unipotent character is the trivial character. The constituents
of $\RLG(1)$ are, by definition, the unipotent characters in the principal
3-Harish-Chandra series. For $\bG^F={}^3\ty{D}_4(q)$ the only 3-split but not
1-split Levi subgroups are the Sylow 3-tori and the centralizers of a 3-torus
of rank~1, with derived subgroup of type $\ty{A}_2$ but $\ty{A}_2$ does not
possess $E$-cuspidal unipotent characters, so again we only get the principal
3-Harish-Chandra series. The situation for $\ty{F}_4(q)$ is entirely similar.
For $\ty{E}_6(q)$ we obtain, in addition, the unipotent characters in the
3-Harish-Chandra series above the 1- and 3-cuspidal character $^3\ty{D}_4[-1]$
of a 3-split Levi subgroup of type $^3\ty{D}_4$ and the characters in the
principal 9-Harish-Chandra series. The arguments for the other groups of
type $\ty{E}_n$ are entirely analogous.
\end{proof}

For the bad prime $\ell=2$ we expect quite a few unipotent $e$-GC characters,
we will not go into this here.

\begin{prop}   \label{prop-JGC}
 Suppose that $\bG$ is simple of simply connected type, $\ell$ is odd and good
 for~$\bG$, does not divide $|\Z(\bG)^F|$, and $\ell>3$
 if~$\bG^F={}^3\ty{D}_4(q)$. Set $e:=e_\ell(q)$. Then $\chi\in\Irr(\bG^F)$ is
 $e$-JGC if and only if $\chi$ is $e$-Jordan-cuspidal.
\end{prop}

\begin{proof}
Thanks to Lemma~\ref{lem:e-J-eJGC} it suffice to show the necessity.
Let $\chi$ be an $e$-JGC character of~$\bG^F$. If $\chi$ is unipotent, then
the claim is a consequence of Lemma~\ref{lem:odd-eGC}. Now assume
$\chi\in\cE(\bG^F,s)$ with $1\ne s\in{\bG^*}^F$ semisimple.
Let $\bH=\C_{\bG^*}^\circ(s)$. By definition, $\Z^\circ(\bH)_{\phi_e}=1$ since
$\Z(\bG^*)=1$. Let $\psi$ be a unipotent $e$-GC character of $\bH$ which
corresponds to $\chi$ under Jordan decomposition. 
We are left to prove that $\psi$ is $e$-cuspidal.

Assume not. Then by Lemma~\ref{lem:noncus-GC}, $\bH$ has
a component of type $\ty{A}_n(\eps q^m)$ with $\ell|(q^m-\eps)$ and
$n+1=\ell^k>1$. Note that for $\bH$ to have a component of type $^3\ty{D}_4$
the group $\bG$ has to be of exceptional type, but then $\ell=3$ is bad
for~$\bG$, contrary to assumption.
Assume that $\bG$ is of exceptional type. Then $\ell\ge5$, and $\ell\ge7$ if
$\bG$ is of type $\ty{E}_8$. If $\bG=\ty{E}_6$, then any centralizer with
an $\ty{A}_4(\eps q^m)$-component has $\Z^\circ(\bH)_{\phi_1}\ne1$; note that
here $\eps q^m=q$, so $e=1$. Similarly, if $\bG=\ty{E}_7$, any centralizer $\bH$
with an $\ty{A}_{\ell-1}(q^m)$-component has
$\Z^\circ(\bH)_{\phi_1}\ne1$, for $\ell\in\{5,7\}$, and any centralizer $\bH$ in
$\bG=\ty{E}_8$ with an $\ty{A}_6(q^m)$-component has
$\Z^\circ(\bH)_{\phi_1}\ne1$. (These claims can be checked easily in
{\sf Chevie} \cite{MChev} using the command {\tt Twistings}). Thus, $\bG$ is of
classical type. But then all centralizers of semisimple elements with a
component $\ty{A}_n(\eps q^m)$ have $|\Z^\circ(\bH)^F|$ divisible by
$q^m-\eps$, so $\Z^\circ(\bH)_{\phi_e}\ne1$, contradiction.
\end{proof}

We provide the following evidence for the validity of
Conjecture~\ref{conj:e-JGC}.

\begin{cor}   \label{cor:e-JGC}
 Let $\bG$ be connected reductive with Frobenius map $F$ such that $[\bG,\bG]$
 is of simply connected type. Assume that $\ell$ is good for $\bG$ with
 $\ell\nmid 2|\Z([\bG,\bG])^F|$, and $\ell>3$ if $\bG^F$ has a component of
 type $^3\ty{D}_4$. Then Conjecture~\ref{conj:e-JGC} holds for all
 $\chi\in\cE(\bG^F,\ell')$.
\end{cor}	

\begin{proof}
Let $\chi\in\Irr(\bG^F)$ and set $e:=e_\ell(q)$.
By	\cite[Thm.~4.2 and Rem.~5.2]{CE99}, $\chi$ is $e$-Jordan-cuspidal if and
only if it is $e$-cuspidal, while by Lemma~\ref{lem:odd-eGC}, $\chi$ is
$e$-cuspidal if and only if it is $e$-GC. 
Therefore, to prove Conjecture~\ref{conj:e-JGC}, it suffices to show that
$\chi$ is $e$-JGC if and only if it is $e$-Jordan-cuspidal.
By \cite[Lemma~2.3]{KM15} and Lemma~\ref{lem:Clifford-eJGC}, we may assume that
$\bG$ is semisimple, as in the proof of Lemma~\ref{lem:noncus-GC}.
Then this assertion follows from Proposition~\ref{prop-JGC} immediately.
\end{proof}

%%%%%%%%%%%%%%%%%%%%%%%%%%%%%%%%%%%%%
\subsection{Generic weights}\label{subsec:generic-weights}

For an $\ell$-block $B$ of $\bG^F$, we denote by $\cL(B)$ the set of $e$-JGC
pairs $(\bL,\la)$ of $\bG^F$ such that $\la\in\cE(\bL^F,\ell')$ and there is
some $\chi\in\Irr(B)$ with
$\langle\chi,\R_{\bL\subseteq \bP}^\bG(\la)\rangle\ne0$ for any parabolic
subgroup $\bP$ of $\bG$ containing $\bL$ as a Levi subgroup.

\begin{defn}\label{defn:cW}
Let $B$ be an $\ell$-block of $\bG^F$ and let $\bT$ be an $e$-torus of $\bG$. 
\begin{enumerate}[\rm(a)]
 \item If $\bT=\Z^\circ(\C_\bG(\bT))_{\phi_e}$, then we define
  \[\cW^0(B,\bT):=\{\,\eta\in\rdz(\N_{\bG^F}(\bT)\mid\la)\mid
    \la\in\cE(\C_{\bG^F}(\bT),\ell')\text{ with } (\C_\bG(\bT),\la)\in\cL(B)\,\}.\]
 \item If $\bT\ne\Z^\circ(\C_\bG(\bT))_{\phi_e}$, then $\cW^0(B,\bT):=\emptyset$.
\end{enumerate}	
\end{defn}

Define $\cW^0(\bG^F,\bT)$ to be the union of the $\cW^0(B,\bT)$ where $B$ runs
through the blocks of $\bG^F$. For an $e$-split Levi subgroup $\bL$ of $\bG$,
we write $\cW^0(B,\bL)=\cW^0(B,\Z^\circ(\bL)_{\phi_e})$ and
$\cW^0(\bG^F,\bL)=\cW^0(\bG^F,\Z^\circ(\bL)_{\phi_e})$.
When $\bL=\bG$, sometimes we also abbreviate $\cW^0(B,\bG)$ (resp.\ 
$\cW^0(\bG^F,\bG)$) to $\cW(B,\bG)$ (resp.\ $\cW(\bG^F,\bG)$).

\begin{defn}\label{defn-pairs11}
 When $\bT\le\bG$ is an $e$-torus and $\eta\in\cW^0(\bG^F,\bT)$, we call
 $(\bT,\eta)$ a \emph{generic $(e,\ell)$-weight} of $\bG^F$. We say that a
 generic $(e,\ell)$-weight $(\bT,\eta)$ \emph{belongs to a block $B$} if
 $\eta\in\cW^0(B,\bT)$.
\end{defn}

Denote by $\cW^0(\bG^F)$ the set of generic $(e,\ell)$-weights of $\bG^F$,
and by $\cW^0(B)$ the set of generic $(e,\ell)$-weights of $\bG^F$ belonging
to~$B$. Set $\cW(\bG^F)=\cW^0(\bG^F)/\!\sim_{\bG^F}$ and
$\cW(B)=\cW^0(B)/\!\sim_{\bG^F}$. For $(\bT,\eta)\in \cW^0(\bG^F)$, we write
$\overline{(\bT,\eta)}$ for the $\bG^F$-conjugacy class of $(\bT,\eta)$.

We highlight the following maximal extendibility property, which is also
involved in the study of the inductive conditions of the McKay conjecture and
the Alperin--McKay conjecture.

\begin{amp}\label{ext-char-e-cuspidal}
 Let $e$ be a positive integer and $\bL$ be an $e$-split Levi subgroup
 of~$\bG$. Then $\ze$ extends to $\N_{\bG^F}(\bL,\ze)$ for every
 $\ze\in \Irr(\bL^F)$.	
\end{amp}

If $\bG$ is a simple simply connected linear algebraic group, then
Assumption~\ref{ext-char-e-cuspidal} holds if $\bL$ is the centralizer of a
Sylow $e$-torus of $\bG$ by Sp\"ath~\cite{Sp09,Sp10a,Sp10b}.
But in general it is still open.

\begin{prop}
 If $\bG$ is simple and of simply connected type, then
 Assumption~\ref{ext-char-e-cuspidal} holds if $\bG$ is of classical type
 $\ty{A}$, $\ty{B}$ or $\ty{C}$, or of exceptional type $\ty{F}_4$.
\end{prop}	

\begin{proof}
If $\bG$ is of type $\ty{A}$, $\ty{C}$ or $\ty{F}_4$, then this assertion holds
by \cite[Thm.~1.2]{BS20}, \cite[Thm.~1.2]{Br22} and \cite[Cor.~4.20]{AHL21}
respectively. If $\bG$ is of type $\ty{B}$, then this assertion is indeed
proved in \cite[\S7]{FLZ22a}.
\end{proof}	

Usually $W_{\bG^F}(\bL,\ze):=\N_{\bG^F}(\bL,\ze)/\bL^F$ denotes the relative
Weyl group of a pair $(\bL,\ze)$ in $\bG$ where $\bL$ is a $F$-stable Levi
subgroup of $\bG$ and $\ze\in\Irr(\bL^F)$.

\begin{lem}   \label{lem-genwei-dz}
 Let $(\bL,\la)\in\cL(B)$. Under Assumption~\ref{ext-char-e-cuspidal} the sets
 $\rdz(\N_{\bG^F}(\bL)\mid\la)$ and $\dz(W_{\bG^F}(\bL,\ze))$ are in bijection.
 In particular, in Definition~\ref{defn:cW} (a), $\cW^0(B,\bT)$ is
 parameterized by defect zero characters of relative Weyl groups.
\end{lem}

\begin{proof}
This follows immediately from Gallagher's theorem.
\end{proof}

\begin{defn}
 Let $\bG\embed\tbG$ be a weakly regular embedding. 
 \begin{enumerate}[\rm(a)]
	\item Let $(\bT,\eta)\in\cW^0(\bG^F)$ and $(\tbT,\widetilde\eta)\in\cW^0(\tbG^F)$. We say that $(\tbT,\widetilde\eta)$ \emph{covers} $(\bT,\eta)$ if $\tbT=\bT\Z^\circ(\tbG)$ and $\widetilde\eta\in\Irr(\N_{\tbG^F}(\tbT)\mid\eta)$, and say that $\overline{(\tbT,\widetilde\eta)}$ \emph{covers} $\overline{(\bT,\eta)}$ if $(\tbT,\widetilde\eta)$ covers $(\bT^g,\eta^g)$ for some $g\in\tbG^F$.
	\item We write $\cW^0(\tbG^F\mid(\bT,\eta))$ for the set of those $(\tbT,\widetilde\eta)\in\cW^0(\tbG^F)$ covering $(\bT,\eta)\in\cW^0(\bG^F)$, while we write $\cW^0(\bG^F\mid(\tbT,\widetilde\eta))$ for the set of those $(\bT,\eta)\in\cW^0(\bG^F)$ covered by $(\tbT,\widetilde\eta)\in\cW^0(\tbG^F)$.
	\item Write $\cW(\tbG^F\mid\overline{(\bT,\eta)})$ for the set of those $\overline{(\tbT,\widetilde\eta)}\in\cW(\tbG^F)$ covering $\overline{(\bT,\eta)}\in\cW(\bG^F)$, while we write $\cW(\bG^F\mid\overline{(\tbT,\widetilde\eta)})$ for the set of those $\overline{(\bT,\eta)}\in\cW(\bG^F)$ covered by $\overline{(\tbT,\widetilde\eta)}\in\cW(\tbG^F)$.
 \end{enumerate}
\end{defn}

Now assume further that $\ell$ is good for $\bG$.  Let~$e:=e_\ell(q)$.
If $\bL$ is an $e$-split Levi subgroup of $\bG$ and $b$ is an $\ell$-block of
$\bL^F$, then by \cite[Thm.~2.5]{CE99} there exists an $\ell$-block of $\bG^F$,
denoted $\RLG(b)$, such that for any $\ze\in\cE(\bL^F,\ell')\cap\Irr(b)$
and any parabolic subgroup $\bP$ of $\bG$ containing $\bL$ as a Levi subgroup,
one has $\R^\bG_{\bL\subseteq\bP}(\ze)\in\ZZ\Irr(\RLG(b))$.
This implies that $\cW^0(B_1)\cap\cW^0(B_2)=\emptyset$ if $B_1$ and $B_2$ are
distinct $\ell$-blocks of $\bG^F$. By \cite[Thm.~3.4]{KM15}, this continues to
hold for bad primes if $\bG$ is a Levi subgroup of some \emph{simple} algebraic
group of simply connected type.

\begin{condition}   \label{condition-1}
	Let $\bG$ be connected reductive and $F\colon\bG\to\bG$ a Frobenius
	endomorphism with respect to an $\FF_q$-structure on $\bG$. Assume
	that~$\ell$ is odd, good for~$\bG$ and does not divide
	$|\cZ(\bG)^F|\, |\cZ(\bG^*)^F|$. Let~$e:=e_\ell(q)$.
\end{condition}

\begin{prop}   \label{prop:block-induce}
 Keep Condition~\ref{condition-1}.
 Let $(\bT,\eta)\in\cW^0(\bG^F)$ and $\bL:=\C_\bG(\bT)$.
 \begin{enumerate}[\rm(a)]
  \item We have $\bL=\C_\bG^\circ(\Z(\bL)_\ell^F)$,
   $\bL^F=\C_{\bG^F}(\Z(\bL)_\ell^F)$ and
   $\N_{\bG^F}(\bT)=\N_{\bG^F}(\Z(\bL)_\ell^F)$.
  \item Let $B$ be an $\ell$-block of $\bG^F$ such that
   $(\bT,\eta)\in\cW^0(B)$. Then $\bl(\eta)^{\bG^F}$ is defined and equals~$B$.
 \end{enumerate}
\end{prop}

\begin{proof}
Note that $\bT=\Z^\circ(\bL)_{\phi_e}$ by definition.
By \cite[Prop.~2.2]{CE94} we have $\bL=\C_\bG^\circ(\Z(\bL)_\ell^F)$ and
$\bL^F=\C_{\bG^F}(\Z(\bL)_\ell^F)$. Thus
$\N_{\bG^F}(\bT)=\N_{\bG^F}(\bL)=\N_{\bG^F}(\Z(\bL)^F_\ell)$ and (a) is shown.

Suppose that $\eta\in\Irr(\N_{\bG^F}(\bT)\mid\la)$ where
$\la\in\cE(\bL^F,\ell')$ is an $e$-JGC character. By definition,
$\RLG(\bl(\la))=B$.  According to \cite[Thm.~2.5]{CE99}, there is an
inclusion of connected subpairs
$(1,B)^0\lhd (\Z(\bL)_\ell^F,\bl(\la))^0$ in the sense of
\cite[Prop.~2.1]{CE99}. Thus $\bL^F=\C_{\bG^F}(\Z(\bL)_\ell^F)$ forces that
$(1,B)\lhd (\Z(\bL)_\ell^F,\bl(\la))$, i.e., $\bl(\la)^{\bG^F}=B$.
As $\N_{\bG^F}(\bT)=\N_{\bG^F}(\Z(\bL)_\ell^F)$, we have that
$\bl(\la)^{\N_{\bG^F}(\bT)}$ is defined and equals $\bl(\eta)$ (see, e.g.,
\cite[Chap.~5, Thm.~5.15]{NT89}). So, by transitivity of block induction,
$\bl(\eta)^{\bG^F}$ is defined and equals~$B$.
\end{proof}	

\begin{lem}   \label{lem:def-zero-pairs}
 Keep Condition~\ref{condition-1}. Let $B$ be an $\ell$-block of $\bG^F$ of
 central defect. Then $\Irr(B)\cap\cE(\bG^F,\ell')$ consists of a unique
 character, denoted $\chi$, and $\cW^0(B)=\{(\Z(\bG)_{\phi_e},\chi)\}$.
 Moreover, $\chi$ is $e$-Jordan-cuspidal.
\end{lem}	

\begin{proof}
Let $G:=\bG^F$.	Since $B$ has central defect, the dominated block of
$G/\Z(G)_\ell$ has defect zero, and hence contains a unique character, and thus
all other characters in $\Irr(B)$ are non-trivial on $\Z(G)_\ell$ and so lie in
series $\cE(G,t)$ with $t$ not an $\ell'$-element.
Therefore, $\Irr(B)\cap\cE(G,\ell')=:\{\chi\}$ is a singleton.
	
Next, we claim that if $(\bL,\la)\in\cL(B)$, then $\bL=\bG$ and $\la=\chi$.
By the proof of Proposition~\ref{prop:block-induce}, if $\bL$ is any $e$-split
Levi subgroup of $\bG$ and $\la\in\cE(\bL^F,\ell')$ with
$\RLG(\bl(\la))=B$, then $\Z(\bL)_\ell^F$ is a subgroup of a defect
group of $B$. This forces that $\Z(\bL)_\ell^F$ is central in $G$, and thus
$\bL=\C_\bG^\circ(\Z(\bL)_\ell^F)=\bG$ and $\la=\chi$.  So the claim holds.

By \cite[Thm.~3.6]{KM15}, there exists an $e$-Jordan-cuspidal pair $(\bL,\la)$
of $\bG$ with $\la\in\cE(\bL^F,\ell')$ and $\RLG(\bl(\la))=B$,
and by Lemma~\ref{lem:e-J-eJGC}, $\la$ is also $e$-JGC.
So by the previous paragraph $\bL=\bG$ and $\chi=\la$ is $e$-JGC. From this
$\cW^0(B)=\{(\Z(\bG)_{\phi_e},\chi)\}$ by definition.
\end{proof}

\begin{cor}   \label{cor-par-gen-wei}
 Suppose that $\bG$ is simple of simply connected type, $\ell$ is odd and good
 for~$\bG$, does not divide the order of $\Z(\bG)^F$ and $\ell>3$ if
 $\bG^F={}^3\ty{D}_4(q)$. Let $B$ be an $\ell$-block of $\bG^F$, and let
 $e:=e_\ell(q)$. Then the set $\cW(B,\bG)$ is non-empty if and only if $B$ is
 of defect zero. In particular, $\cW(\bG^F,\bG)=\dz(\bG^F)$.
\end{cor}

\begin{proof}
The	sufficiency follows by Lemma~\ref{lem:def-zero-pairs}. Now we prove the
necessity and assume that $\cW(B,\bG)\ne\emptyset$, which implies that
there exists an $e$-JGC character $\chi\in\Irr(B)\cap\cE(\bG^F,\ell')$
of~$\bG^F$. By Proposition~\ref{prop-JGC}, $\chi$ is also $e$-Jordan-cuspidal.
Now by \cite[Thm.~A]{KM15}, $\chi$ is of quasi-central defect in the
sense of \cite[Def.~2.4]{KM13}. Since $[\bG,\bG]=\bG$, we have that $\chi$ is
of central defect, and thus $B$ is of defect zero since
$\ell\nmid |\Z(\bG)^F|$.	
\end{proof}	
	
%%%%%%%%%%%%%%%%%%%%%%%%%%%%%%%%%%%%%
\subsection{Reduction to quasi-isolated blocks} 

Let $s\in {\bG^*}^F$ be a semisimple $\ell'$-element and $B$ be an $\ell$-block
of $\bG^F$ with $\Irr(B)\subseteq\cE_\ell(\bG^F,s)$. Let $\bL^*$ be an
$F$-stable Levi subgroup of $\bG^*$ with
$\C^\circ_{\bG^*}(s)\C_{{\bG^*}^F}(s)\subseteq\bL^*$ and let $C$ be the block
of $\bL^F$ that is the Bonnaf\'e--Dat--Rouquier correspondent
(cf.\ \cite[Thm.~7.7]{BDR17}) of $B$. Let $e:=e_\ell(q)$.

\begin{lem}   \label{lem:eJGC-inj}
 The map
 $$\Xi:(\bL'',\la'')\mapsto(\bL',\la'):=
   \big(\C_{\bG}(\Z^\circ(\bL'')_{\phi_e}),\pm\R^{\bL'}_{\bL''}(\la'')\big)$$
 is an injection $\cL(C)\to\cL(B)$.
 Moreover, every element in $\cL(B)$ has a $\bG^F$-conjugate in $\Xi(\cL(C))$.
\end{lem}

\begin{proof}
Let $(\bL',\la')\in\cL(B)$ where $\cL(B)$ is defined as in
\S\ref{subsec:generic-weights}. Then there is a semisimple $\ell'$-element
$s'\in\bL'^{*F}$ such that $\la'\in\cE(\bL'^F,s')$, and by definition we obtain
$\Irr(B)\cap\cE_\ell(\bG^F,s')\ne\emptyset$. Hence $s$ and $s'$ are conjugate
in ${\bG^*}^F$, and so up to conjugation we may assume that $s=s'\in{\bL'}^*$.
Thus $\Z^\circ({\bL'}^*)\subseteq \C_{\bG^*}^\circ(s)\subseteq \bL^*$.
Let ${\bL''}^*:={\bL'}^*\cap{\bL}^*=\C_{\bL^*}(\Z^\circ({\bL'}^*)_{\phi_e})$.
Then ${\bL''}^*$ is an $e$-split Levi subgroup of $\bL^*$ with $s\in{\bL''}^*$.
As $\C^\circ_{\bG^*}(s)\C_{{\bG^*}^F}(s)\subseteq\bL^*$, we obtain
$\C^\circ_{\bL'^*}(s)\C_{{\bL'^*}^F}(s)\subseteq\bL''^*$, and hence
$\C^\circ_{\bL'^*}(s)=\C^\circ_{\bL''^*}(s)$ and $\Z^\circ({\bL'}^*)\subseteq
\Z^\circ({\bL''}^*)\subseteq\Z^\circ(\C^\circ_{{\bL'}^*}(s))$.
By $(\bL',\la')\in\cL(B)$, we have
$\Z^\circ({\bL'}^*)_{\phi_e}=\Z^\circ(\C^\circ_{{\bL'}^*}(s))_{\phi_e}$.
Therefore, $\Z^\circ({\bL'}^*)_{\phi_e}=\Z^\circ({\bL''}^*)_{\phi_e}$, and thus
${\bL'}^*=\C_{\bG^*}(\Z^\circ({\bL''}^*)_{\phi_e})$. Let $\bL''$ be an
$e$-split Levi subgroup of $\bL$ in duality with ${\bL''}^*$. Then the
previous argument shows that, after conjugation, we may assume $\bL''\le\bL'$.

By \cite[Thm.~3.3.22]{GM20}, there exists a unique $\la''\in\cE({\bL''}^F,s)$
such that $\la'=\pm\R^{\bL'}_{\bL''}(\la'')$. Then by \cite[Thm.~4.7.1]{GM20},
$\la''$ corresponds via Jordan decomposition to the same unipotent character of
$\C^\circ_{\bL'^*}(s)^F=\C^\circ_{\bL''^*}(s)^F$ as $\la'$, whence $\la''$ is
an $e$-JGC character of~$\bL''^F$. Hence $(\bL'',\la'')\in\cL(C)$ by
construction.  Moreover, by \cite[Prop.~13.8]{CE04},
$\Z^\circ(\bL')_{\phi_e}=\Z^\circ(\bL'')_{\phi_e}$, which implies
that $\bL''=\bL'\cap\bL$.

Conversely, we let $(\bL'',\la'')\in\cL(C)$ and
$\bL':=\C_{\bG}(\Z^\circ(\bL'')_{\phi_e})$ so that $\bL''=\bL'\cap\bL$.
Then $\Z^\circ(\bL')_{\phi_e}=\Z^\circ(\bL'')_{\phi_e}$, and from
\cite[Prop.~13.8]{CE04} we get
$\Z^\circ({\bL'}^*)_{\phi_e}=\Z^\circ({\bL''}^*)_{\phi_e}$.
Similarly as above, we may assume that $s\in\bL''^*$. Moreover,
$${\bL''}^*=\C_{\bL^*}(\Z^\circ({\bL''}^*)_{\phi_e})
  =\C_{\bL^*}(\Z^\circ({\bL'}^*)_{\phi_e})={\bL'}^*\cap{\bL}^*.$$
Therefore, $\C^\circ_{\bL'^*}(s)\C_{{\bL'^*}^F}(s)\subseteq\bL''^*$ and
$\la':=\pm\R^{\bL'}_{\bL''}(\la'')\in\Irr(\bL'^F)$.
By $(\bL'',\la'')\in\cL(C)$, we have
$\Z^\circ({\bL''}^*)_{\phi_e}=\Z^\circ(\C^\circ_{{\bL''}^*}(s))_{\phi_e}$ and
so $\Z^\circ({\bL'}^*)_{\phi_e}=\Z^\circ(\C^\circ_{{\bL'}^*}(s))_{\phi_e}$.
Thus similarly as above, we have $(\bL',\la')\in\cL(B)$.
\end{proof}

\begin{cor}   \label{cor:eJGC-bij}
 Keep Condition~\ref{condition-1}. The map $\Xi$ from Lemma~\ref{lem:eJGC-inj}
 induces a bijection
 \[\cL(C)/{\sim_{\bL^F}}\to\cL(B)/{\sim_{\bG^F}}.\]
\end{cor}

\begin{proof}
Note that $\ell$ is also good for $\bL$, and by \cite[Prop.~13.12(ii)]{CE04},
$\cZ(\bL)^F$ is of $\ell'$-order. As in the proof of
Proposition~\ref{prop:block-induce}, $\bL'^F=\C_{\bG^F}(\Z(\bL'')^F_\ell)$
and $\bL''^F=\C_{\bL^F}(\Z(\bL'')^F_\ell)$.
Therefore, since the Brauer categories of splendid Rickard equivalent blocks
are equivalent, by \cite[Rem.~4.7]{FLZ22b}, the Bonnaf\'e--Dat--Rouquier
splendid Rickard equivalence induces a bijection between the $\bL^F$-conjugacy
classes of $C$-Brauer pairs of $\bL^F$ and the $\bG^F$-conjugacy classes of
$B$-Brauer pairs of~$\bG^F$.  So this assertion follows by
Lemma~\ref{lem:eJGC-inj}.
\end{proof}

\begin{thm}
 Keep Condition~\ref{condition-1}. Then
 \[(\bT,\eta)\mapsto(\bT,\pm\R^{\N_\bG(\bT)}_{\N_\bL(\bT)}(\eta))\]
 induces a bijection $\cW(C)\to\cW(B)$.
\end{thm}

\begin{proof}
Let $e:=e_\ell(q)$ and let $\bT$ be an $e$-torus of $\bL$.
The normalizer $\N_{\bG^F}(\bT)=\N_{\bG^F}(\bT^F_\ell)$ is a local subgroup of
$\bG^F$ as in (the proof of) Proposition~\ref{prop:block-induce}.
Let $c$ be the union of blocks of $\N_{\bL^F}(\bT)$ whose induced block to
$\bL^F$ is $C$ and let $b$ be the union of blocks of $\N_{\bG^F}(\bT)$ whose
induced block to $\bG^F$ is $B$. By \cite[Thm.~3.10]{Rh22},
$\pm\R^{\N_\bG(\bT)}_{\N_\bL(\bT)}\colon\Irr(c)\to\Irr(b)$ is a bijection
(induced by a Morita equivalence).
Furthermore, by Lemma~\ref{lem:eJGC-inj} and the construction of this Morita
equivalence in the proof of \cite[Thm.~3.10]{Rh22}, $(\bT,\eta)$ is a generic
$(e,\ell)$-weight of $\bL^F$ if and only if
$(\bT,\pm\R^{\N_\bG(\bT)}_{\N_\bL(\bT)}(\eta))$ is a generic $(e,\ell)$-weight
of $\bG^F$. Thus the assignment
$(\bT,\eta)\mapsto(\bT,\pm\R^{\N_\bG(\bT)}_{\N_\bL(\bT)}(\eta))$ is
well-defined between generic weights. The bijectivity between $\cW(C)$ and
$\cW(B)$ follows directly from Corollary~\ref{cor:eJGC-bij} and
\cite[Rem.~4.9]{FLZ22b}.
\end{proof}

Therefore, we can reduce the determination of generic weights to quasi-isolated
blocks inductively.

%%%%%%%%%%%%%%%%%%%%%%%%%%%%%%%%%%%%%%%%%%%%%%%%%%%%%%%%%%%%%%%%%%%%%%%%%
%%%%%%%%%%%%%%%%%%%%%%%%%%%%%%%%%%%%%%%%%%%%%%%%%%%%%%%%%%%%%%%%%%%%%%%%%
\section{Weights of finite reductive groups}   \label{sec:wei-par}

In this section, we partition the weights of a finite reductive group into
several families in terms of the centers of radical subgroups. This will be
used in the following sections to compare weights with generic weights.
Throughout this section, $\bG$ denotes a connected reductive group with a
Frobenius endomorphism $F\colon\bG\to\bG$ endowing $\bG$ with an
$\FF_q$-structure. We let $\ell$ be a prime not dividing $q$.

\begin{notation}[{\cite[2.3]{CE94}}]   \label{notation:ab}
	%Assume that $\ell$ is good for $\bG$.
	Let $\bG_{\ba}$ be the central product in $\bG$ of $\Z^\circ(\bG)$ and all
	the rationally irreducible components of $[\bG,\bG]$ of type
	$(\ty{A}_n,\eps q^m)$ with $\ell$ dividing $q^m-\eps$.
	Let $\bG_\bb$ be the central product of the rationally irreducible
	components of $[\bG,\bG]$ which are not included in $\bG_{\ba}$.	
\end{notation}

In the situation of Notation~\ref{notation:ab}, $\bG=\bG_{\ba}\bG_\bb$ is a
central product, and $\Z(\bG_\bb)^F$ and $\bG^F/\bG_{\ba}^F \bG_\bb^F$ are
abelian $\ell'$-groups; in addition, if $R$ is an $\ell$-subgroup of $\bG^F$
such that $\Z(\C_{\bG^F}(R))_\ell\subseteq \Z(\bG)\bG_{\ba}$, then
$R\subseteq \bG_{\ba}$ (see \cite[p.~156]{CE94}).

%%%%%%%%%%%%%%%%%%%%%%%%%%%%%%%%%%%%%
\subsection{Groups of Lie type with abelian Sylow subgroups}
We recall the description of $\ell$-weights for groups of Lie type with abelian
Sylow~$\ell$-subgroups in \cite{Ma14}.

In this subsection, $\bG$ is moreover assumed to be simple.
If $\bG^F$ has abelian Sylow $\ell$-subgroups, then by \cite[\S 2.1]{Ma14},
there is a unique positive integer $e$ such that $\ell\mid \phi_e(q)$ and
$\phi_e$ divides the order polynomial of $(\bG,F)$, in which case~$\ell$ is
odd, good for $\bG$, $\ell>3$ if $(\bG,F)$ has type $^3\ty{D}_4$, and
does not divide the orders of $\cZ(\bG)^F$ and $\cZ(\bG^*)^F$.
We let $e:=e_\ell(q)$.

The radical $\ell$-subgroups of $\bG^F$ can be classified in terms of $e$-split
Levi subgroups.

\begin{prop}[{\cite[Cor.~3.2]{Ma14}}]
 Assume that $\bG^F$ has abelian Sylow $\ell$-subgroups. Then
 $R\mapsto\C_\bG(R)$ gives a bijection,with inverse
 $\bL\mapsto\Z(\bL)^F_\ell$, between the set of radical $\ell$-subgroups $R$ of
 $\bG^F$ and the set of $e$-split Levi subgroups~$\bL$ of~$\bG$ with
 $\Z(\bL)^F_\ell=\mrO_\ell(\bL^F)$.
\end{prop}

If $R$ and $\bL$ correspond to each other as above, then
$\N_\bG(R)=\N_\bG(\bL)$. Let $\bL$ be an $e$-split Levi subgroup of $\bG$ and
$\ze\in\dz(\bL^F/\mrO_\ell(\bL^F))$ be of defect zero. Then
$\ze\in\cE(\bL^F,\ell')$ and $\ze$ is $e$-cuspidal
by \cite[Prop.~3.4]{Ma14}.

\begin{thm}[{\cite[\S3]{Ma14}}]   \label{thm:wei-abel-Sylow}
 Assume that $\bG^F$ has abelian Sylow $\ell$-subgroups.
 Suppose that $(\bL,\la)$ is an $e$-cuspidal pair of $\bG$ with
 $\la\in\cE(\bL^F,\ell')$ and $R=\Z(\bL)^F_\ell$. Let
 $B=\RLG(\bl(\la))$.
 \begin{enumerate}[\rm(a)]
  \item $R$ is a defect group of $B$ and $\bl(\la)^{\bG^F}=B$. 
  \item Up to conjugation, the $B$-weights are $(R,\vhi)$ with
   $\vhi\in\rdz(\N_{\bG^F}(\bL)\mid\la)$.
 \end{enumerate}
\end{thm}

In Theorem~\ref{thm:wei-abel-Sylow}, if we assume further that
Assumption~\ref{ext-char-e-cuspidal} holds for $\bG$, then the conjugacy
classes of $B$-weights are in bijection with the irreducible characters
of~$W_{\bG^F}(\bL,\ze)$.

The result in \cite[\S3]{Ma14} is only stated and proved for primes $\ell\ge5$.
Now note that Sylow 2-subgroups of $\bG^F$ are never abelian, and that for
$\ell=3$, if Sylow 3-subgroups are abelian, then $Z=D$ in the setting of
\cite[Rem.~5.2]{CE99}, whence all ingredients in the proof of \cite[\S3]{Ma14}
taken from \cite{CE99} continue to hold, by \cite[Rem.~5.2]{CE99}.

Note that if $\bG^F$ is an abstract quasi-simple group and has abelian Sylow
2- or 3-subgroups, then the inductive BAW condition holds for every 
$\ell$-block of $\bG^F$, for any $\ell$; see \cite[Cor.~6.6]{Sp13} and
\cite[\S5]{FLZ23} (and the references therein). The construction of the weights
of $\bG^F$ in those cases can be found in those papers. 

\begin{cor}\label{cor-genwei-abel-def}
 Assume that $\bG^F$ has abelian Sylow $\ell$-subgroups. Then
 \[(R,\vhi)\mapsto (\Z^\circ(\C_\bG(R))_{\phi_e},\vhi)\]
 induces a canonical bijection $\vOm\colon\Alp(\bG^F)\to\cW(\bG^F)$ such
 that $\vOm(\Alp(B))=\cW(B)$ for every $\ell$-block $B$ of $\bG^F$. 
\end{cor}
	
\begin{proof}
According to \cite[Thm.]{CE99}, $(\bL,\la)\mapsto \RLG(\bl(\la))$
induces a bijection between the $\bG^F$-conjugacy classes of $e$-cuspidal pairs
$(\bL,\la)$ of $\bG$ with $\la\in\cE(\bL^F,\ell')$ and the $\ell$-blocks of
$\bG^F$. By Theorem~\ref{thm:wei-abel-Sylow}, it suffices to show that for
$(\bG,F)$ the $e$-cuspidal pairs and the $e$-JGC pairs coincide, which follows
from Lemma~\ref{cor:abel-sylow}.
\end{proof}

Some of the above results can be generalized to blocks with
abelian defect groups.

%%%%%%%%%%%%%%%%%%%%%%%%%%%%%%%%%%%%%
\subsection{Blocks with abelian defect groups}   \label{subsec:Abelian-defect-groups}

Let $\bH$ be a simple algebraic group of simply connected type with a
Frobenius endomorphism $F:\bH\to\bH$ endowing $\bH$ with an $\FF_q$-structure.
Let $\bG$ be an $F$-stable Levi subgroup of $\bH$.
Let $\ell$ be a prime not dividing $q$ such that $\ell$ is odd and good for
$\bG$. Assume that$\ell>3$ if $\bG^F={}^3\ty{D}_4(q)$. Let $e:=e_\ell(q)$.

Let $B$ be an $\ell$-block of $\bG^F$. Under our conditions
$\Irr(B)\cap\cE(\bG^F,\ell')$ is a basic set for $B$ (see
\cite[Thm.~14.4]{CE04}). Then by \cite[Thm.~A(e)]{KM15}, up to conjugacy, there
exists a unique $e$-Jordan-cuspidal pair $(\bL,\ze)$ of $\bG$ with
$\ze\in\cE(\bL^F,\ell')$ and $B=\RLG(\bl(\ze))$.
Thus, if $ \Irr(B)\cap\cE(\bG^F,\ell')$ satisfies a generalized
$e$-Harish-Chandra theory in the sense of \cite[Thm.~1.4]{KM13} then
$$|\IBr(B)|=|\Irr(B)\cap\cE(\bG^F,\ell')|=|\Irr(W_{\bG^F}(\bL,\ze))|.$$
Note that the generalized $e$-Harish-Chandra theory is known to hold in many
situations, for example whenever Lusztig's Jordan decomposition is known to
commute with Lusztig induction.

\begin{prop}   \label{prop:generic-abeldef}
 Keep the above hypotheses and assume further that $B$ has abelian defect
 groups. Then we have:
 \begin{enumerate}[\rm(a)]
  \item Up to conjugation, $(\bL,\ze)$ is the only $e$-JGC pair with
   $\ze\in\cE(\bL^F,\ell')$ and $B=\RLG(\bl(\ze))$. 
  \item The relative Weyl group $W_{\bG^F}(\bL,\ze)$ is of $\ell'$-order.
  \item If Assumption~\ref{ext-char-e-cuspidal} holds, then 
   $|\cW(B)=|\Irr(W_{\bG^F}(\bL,\ze))|$. 
 \end{enumerate}	
\end{prop}

\begin{proof}
Let $s\in\bG^{*F}$ be a semisimple $\ell'$-element such that
$\Irr(B)\subseteq \cE_\ell(\bG^F,s)$. Let $(\bL,\ze)$ let an $e$-JGC pair with
$\ze\in\cE(\bL^F,\ell')$ and $B=\RLG(\bl(\ze))$. Up to conjugation we may
assume that $s\in{\bL^*}^F$.
Recall from the proof Proposition~\ref{prop:block-induce} that
$\bl(\ze)^{\bG^F}=B$, whence $\bl(\ze)$ also has abelian defect groups.
By \cite[Prop.~5.1]{CE99}, $\bl(\ze)$ corresponds to a unipotent block of
${\bL(s)^\circ}^F$ such that they possess isomorphic defect groups. Here
$\bL(s)$ is a closed subgroup of $\bL$ such that $\bL(s)^\circ$ is in
duality with $\C^\circ_{\bL^*}(s)$ and $\bL(s)/\bL(s)^\circ$ is isomorphic to
$\C_{\bL^*}(s)/\C^\circ_{\bL^*}(s)$.
We will prove (a) by showing that there does not exist a non $e$-cuspidal
$e$-GC unipotent character of ${\bL(s)^\circ}^F$, which implies that $\ze$ must
be $e$-Jordan-cuspidal. According to Lemma~\ref{lem:noncus-GC}, it suffices to
prove that $\bL(s)^\circ$ does not
have a component of type $\ty A_{n}(\eps q^m)$ with $n+1=\ell^k$ and
$\ell\mid(q^m-\eps)$. If $\bL(s)^\circ$ has such a component, then by
\cite[Thm.~21.14]{CE04}, a Sylow $\ell$-subgroup of this component is a
subgroup of some defect group of~$\bl(\ze)$. However, it follows from
\cite[Prop.~2.2]{Ma14} that the Sylow $\ell$-subgroups of this component are
non-abelian, and this is a contradiction. Thus (a) holds.

Under our conditions the relative Weyl group $W_{\bG^F}(\bL,\ze)$ is of
$\ell'$-order by \cite[Lemma 4.16]{CE99}. Now
\[\cW(B)=\{ (\bT,\eta)\mid\eta\in\rdz(\N_{\bG^F}(\bT)\mid\ze)\}\]
where $\bT=\Z^\circ(\bL)_{\phi_e}$. Under Assumption~\ref{ext-char-e-cuspidal}
for $\bG$, $\cW(B)$ is in bijection with the irreducible characters of
$W_{\bG^F}(\bL,\ze)$ by Gallagher's theorem. This completes the proof.
\end{proof}

Therefore, by Proposition~\ref{prop:generic-abeldef} we obtain that
$|\IBr(B)|=|\cW(B)|$, modulo the generalized $e$-Harish-Chandra theory and the
maximal extendibility property (Assumption~\ref{ext-char-e-cuspidal}).
In the spirit of Alperin's weight conjecture, can we generalize
Corollary~\ref{cor-genwei-abel-def} and establish a correspondence between
$\cW(B)$ and $\Alp(B)$?
To do this, we will consider the general case, and establish an equivariant
bijection in Theorem~\ref{main-thm-weights}.

%%%%%%%%%%%%%%%%%%%%%%%%%%%%%%%%%%%%%%%%%%%%%%%%%%%%%%%%%%%%
\subsection{Centers of radical subgroups}
Now we consider the general case.

\begin{defn}\label{defn:weights-partition}
	Let~$B$ be an $\ell$-block of $\bG^F$ and~$\bT$ be an $e$-torus of~$\bG$.
	\begin{enumerate}[\rm(a)]
		\item We define $\Alp^0(B,\bT)$ to be the set of $B$-weights $(R,\vhi)$
		of~$\bG^F$ such that $\bT=\Z^\circ(\C^\circ_\bG(\Z(R)))_{\phi_e}$.
		\item We write
		\[\Alp^0(\bG^F,\bT):=\coprod_B\Alp^0(B,\bT)\]
		where~$B$ runs through the $\ell$-blocks of~$\bG^F$.
		\item Denote by $\Alp(\bG^F,\bT)$ (resp. $\Alp(B,\bT)$) the set of
		$\N_{\bG^F}(\bT)$-conjugacy classes of weights in $\Alp^0(\bG^F,\bT)$
		(resp.\ $\Alp^0(B,\bT)$).
	\end{enumerate}
\end{defn}

We remark that the notation $\Alp$ in Definition \ref{defn:weights-partition}
depends not only on $\ell$, but also on the choice of $e$. 

If $*\in\{0,\emptyset\}$ and $\bL$ is an $e$-split Levi subgroup of $\bG$, then
we also write
\[\Alp^*(\bullet,\bL)=\Alp^*(\bullet,\Z^\circ(\bL)_{\phi_e}).\]

\begin{rmk}   \label{rmk:cen-def}
 Let $B$ be an $\ell$-block of $\bG^F$ of central defect. Then
 \[\Alp(B)=\Alp(B,\Z^\circ(\bG)_{\phi_e})=\Alp(B,\bG).\]
\end{rmk}	

\begin{lem}
 Assume further that $\ell$ is good for $\bG$ and $e:=e_\ell(q)$. Then
 \[\Alp^0(B):=\coprod_{\bT}\Alp^0(B,\bT),\quad\textrm{and}\quad
   \Alp^0(\bG^F):=\coprod_{\bT}\Alp^0(\bG^F,\bT)\]
 where~$\bT$ runs through the $e$-tori of~$\bG$.
\end{lem}

\begin{proof}
This follows from the fact that $\C^\circ_\bG(\Z(R))$ is a Levi subgroup of
$\bG^F$ for every $\ell$-subgroup $R$ of $\bG$, using
\cite[Prop.~2.1(ii)]{CE94}.
\end{proof}

\begin{prop}   \label{prop:corr-rad}
 Let $\bG$ be connected reductive and let $\bG\le \bG_i$, $i=1,2$ be
 $\ell$-regular embeddings. Then there exists a bijection
 $\Re^0(\bG_1^F)\to\Re^0(\bG_2^F)$, $R_1\mapsto R_2$, such that
 \begin{enumerate}[\rm(a)]
  \item $R_1\cap \bG^F=R_2\cap \bG^F$;
  \item $\C_\bG(\Z(R_1))=\C_\bG(\Z(R_2))$ and
   $\Z(R_1)\cap \bG^F=\Z(R_2)\cap \bG^F$;
  \item for any $e\ge1$, $\Z^\circ(\C^\circ_{\bG_1}(\Z(R_1)))_{\phi_e}\subseteq\Z(\bG_1)$ if and
   only if $\Z^\circ(\C^\circ_{\bG_2}(\Z(R_2)))_{\phi_e}\subseteq\Z(\bG_2)$; and
  \item $R_1\in\Re_w^0(\bG_1^F)$ if and only if $R_2\in\Re_w^0(\bG_2^F)$.
 \end{enumerate}
 In addition, this induces bijections $\Re(\bG_1^F)\to\Re(\bG_2^F)$ and
 $\Re_w(\bG_1^F)\to\Re_w(\bG_2^F)$.
\end{prop}

\begin{proof}
By Lemma~\ref{lem:regular-embedding}, there exists a connected reductive group
$\tbG$ such that $\bG\le\tbG$ and $\bG_i\le \tbG$ ($i=1,2$) are regular
embeddings. Let $\tZ=\Z(\tbG)^F$. Now $\bG\le \bG_i$ are $\ell$-regular
embeddings, that is, $\ell\nmid |\cZ(\bG_i)_F|$, whence $\tbG^F/\bG_i^F\tZ$
are $\ell'$-groups by \cite[Rem.~1.7.6]{GM20}.
According to Lemma \ref{lem:cenprod-rad-corr}, for $i=1,2$,
$R_i\mapsto R_i\tZ_\ell$ defines a bijection $\Re^0(\bG_i^F)\to\Re^0(\tbG^F)$
with inverse $\tR\mapsto \tR\cap\bG_i^F$. For $R_1\in\Re^0(\bG_1^F)$, we let
$\tR:=R_1 \tZ_\ell$ and $R_2:=\tR\cap\bG_2^F$. This gives a bijection
$\Re^0(\bG_1^F)\to\Re^0(\bG_2^F)$, $R_1\mapsto R_2$.
	
Therefore, (d) follows by Lemma~\ref{lem:cenprod-rad-corr} and (a) follows from
the fact $R_i\cap \bG^F=\tR\cap \bG^F$ for $i=1,2$. Now we consider (b).
Note that $\Z(\tR)=\Z(R_i)\tZ_\ell$, so
$\C_\tbG(\Z(R_1))=\C_\tbG(\Z(\tR))=\C_\tbG(\Z(R_2))$ and
$\Z(R_i)=\Z(\tR)\cap \bG_i^F$.
From this, $\Z(R_i)\cap \bG^F=\Z(\tR)\cap \bG^F$ and (b) holds.

Finally, as $\bG_i=[\bG_i,\bG_i]\Z(\bG_i)=\bG\Z(\bG_i)$ for $i=1,2$ we have
$\C^\circ_{\bG_i}(\Z(R_i))=\C^\circ_\bG(\Z(R_i))\Z(\bG_i)$, so
$\Z^\circ(\C^\circ_{\bG_i}(\Z(R_i)))_{\phi_e}=
\Z^\circ(\C^\circ_\bG(\Z(R_i)))_{\phi_e}\Z^\circ(\bG_i)_{\phi_e}$. Thus
$\Z^\circ(\C^\circ_{\bG_i}(\Z(R_i)))_{\phi_e}\subseteq\Z(\bG_i)$ if and only if
$\Z^\circ(\C^\circ_\bG(\Z(R_i)))_{\phi_e}\subseteq\Z(\bG)$, and we obtain~(c)
since $\C_\bG(\Z(R_1))=\C_\bG(\Z(R_2))$ by (b).
\end{proof}
	
\begin{defn}
 Let $\bG$ be connected reductive with a Frobenius endomorphism
 $F\colon\bG\to\bG$. Suppose that $\bG\embed\tbG$ is an $\ell$-regular
 embedding.
 \begin{enumerate}[\rm(a)]
	\item We let $\Alp_0^0(\bG^F)$ be the set of weights of $\bG^F$ covered
         by weights in $\Alp^0(\tbG^F,\tbG)$.
	\item For an $\ell$-block $B$ of $\bG^F$, we define
		$\Alp^0_0(B):=\Alp^0_0(\bG^F)\cap \Alp^0(B)$.
	\item Define $\Alp_0(\bG^F):=\Alp_0^0(\bG^F)/\!\sim_{\bG^F}$ and
		$\Alp_0(B):=\Alp_0^0(B)/\!\sim_{\bG^F}$.
 \end{enumerate}
\end{defn}
By Proposition~\ref{prop:corr-rad}, $\Alp_0^0(\bG^F)$,
$\Alp^0_0(B)$, $\Alp_0(\bG^F)$, $\Alp_0(B)$ are independent of the choice
of~$\tbG$.

\begin{lem}   \label{lem:weights-ell'center}
 Suppose that $|\cZ(\bG)_F|$ is prime to~$\ell$. Then
 $\Alp^0_0(\bG^F)=\Alp^0(\bG^F,\bG)$.
\end{lem}

\begin{proof}
This follows from the fact that under our assumption $\bG$ can be regarded as
an $\ell$-regular embedding of itself.
\end{proof}

\begin{prop}   \label{prop:center-Levi}
 Keep Condition~\ref{condition-1}. Let~$R$ be a radical $\ell$-subgroup
 of~$\bG^F$ and $\bH:=\C^\circ_\bG(\Z(R))$.  Let~$E:=E_{e,\ell}$. Then:
 \begin{enumerate}[\rm(a)]
  \item $\bH$ is an $E$-split Levi subgroup of $\bG$ with
   $\Z(R)=\Z(\bH)^F_\ell$ and $\bH^F=\C_{\bG^F}(\Z(R))$.
  \item $\N_{\bH^F}(R)\unlhd \N_{\bG^F}(R)$. In particular, $R$ is a radical
   $\ell$-subgroup of~$\bH^F$.
  \item Let $\bL:=\C_\bG(\Z^\circ(\bH)_{\phi_e})$. Then
   $\Z^\circ(\bL)_{\phi_e}=\Z^\circ(\bH)_{\phi_e}$ and
   $\N_{\bL^F}(R)\unlhd \N_{\bG^F}(R)$.
   In particular, $R$ is a radical $\ell$-subgroup of~$\bL^F$.
 \end{enumerate}
\end{prop}	

\begin{proof}
By \cite[Prop.~2.1(iii)]{CE94}, $\bH$ is an $F$-stable Levi subgroup of $\bG$
and $\bH^F=\C_{\bG^F}(\Z(R))$. From
\[\N_{\bH^F}(R)=\{\,g\in\N_{\bG^F}(R)\mid [g,\Z(R)]=1 \,\}\]
we conclude that $\N_{\bH^F}(R)\unlhd \N_{\bG^F}(R)$.
So $R$ is a radical subgroup of $\bH^F$. This proves (b).
Therefore, $\mrO_\ell(\bH^F)\subseteq R$ by \cite[Lemma~2.3]{NT11}. In
particular, $\Z(\bH^F)_\ell\subseteq R$ which implies that
$\Z(\bH^F)_\ell\subseteq \Z(R)$ since $R\subseteq \bH^F$.
As $\Z(R)\subseteq \Z(\bH)$, we conclude that $\Z(R)=\Z(\bH)^F_\ell$.
Thus (a) holds by \cite[Prop.~13.19]{CE04}.

Now we consider (c). As $\bL=\C_\bG(\Z^\circ(\bH)_{\phi_e})$ we have
$\Z^\circ(\bH)_{\phi_e}\subseteq \Z^\circ(\bL)_{\phi_e}$.
Also note that $\bH=\C_\bG(\Z^\circ(\bH))$ and $\bH\subseteq\bL$, then
$\Z^\circ(\bL)\subseteq\Z^\circ(\bH)$ and thus
$\Z^\circ(\bL)_{\phi_e}=\Z^\circ(\bH)_{\phi_e}$.
For the $e$-split Levi subgroup $\bL$ of $\bG$, by \cite[Prop.~13.19]{CE04}, we
have $\bL=\C_\bG^\circ(\Z(\bL)^F_\ell)$. By construction,
$\N_\bG(R)\subseteq \N_\bG(\bL)$ and it follows that
$\Z(\bL)^F_\ell\unlhd \N_{\bG^F}(R)$. By \cite[Prop.~2.2]{CE94},
$\bL^F=\C_{\bG^F}(\Z(\bL)^F_\ell)$ which implies
\[\N_{\bL^F}(R)=\{\, g\in\N_{\bG^F}(R)\mid [g,\Z(\bL)^F_\ell]=1\,\}.\]
From this $\N_{\bL^F}(R)\unlhd \N_{\bG^F}(R)$ and $R$ is a radical
$\ell$-subgroup of $\bL^F$. This completes the proof.
\end{proof}

In Definition~\ref{defn:weights-partition}, if we also assume
Condition~\ref{condition-1}, then it follows from
Proposition~\ref{prop:center-Levi}(c) that $\Alp^0(B,\bT)$
(or $\Alp^0(\bG^F,\bT)$) is non-empty only when
$\bT=\Z^\circ(\C_\bG(\bT))_{\phi_e}$.

\begin{cor}   \label{cor:center-Levi}
 Keep Condition~\ref{condition-1}. Let $\bL$ be an $e$-split Levi subgroup of
 $\bG$ and~$R$ be a radical $\ell$-subgroup of~$\bL^F$ such that
 $\Z^\circ(\C^\circ_\bL(\Z(R)))_{\phi_e}\subseteq \Z(\bL).$
 Then $\bL=\C_\bG(\Z^\circ(\C^\circ_\bG(\Z(R)))_{\phi_e})$ and
 $\N_{\bL^F}(R)\unlhd \N_{\bG^F}(R)$.	
\end{cor}

\begin{proof}
According to \cite[Prop.~13.12]{CE04}, $\ell$ is good for~$\bL$ and does not
divide $|\cZ(\bL)^F|\,|\cZ(\bL^*)^F|$.
Let $\bH:=\C_\bL^\circ(\Z(R))$. From Proposition~\ref{prop:center-Levi} it
follows that $\bH$ is an $E$-split Levi subgroup of~$\bL$ and
$\Z(R)=\Z(\bH)_\ell^F$. Since $\bL$ is an $e$-split Levi subgroup of $\bG$ and
$\Z^\circ(\bH)_{\phi_e}\subseteq \Z(\bL)$ by assumption, we have
$\Z^\circ(\bH)_{\phi_e}=\Z^\circ(\bL)_{\phi_e}$ and thus 
$\bL=\C_\bG(\Z^\circ(\bH)_{\phi_e})$.
By \cite[Prop.~13.19]{CE04}, we have $\bL=\C_\bG^\circ(\Z(\bL)^F_\ell)$.
Since $\Z(\bL)_\ell^F\subseteq \Z(R)$, we have
$\C_\bG^\circ(\Z(R))\subseteq\C_\bG^\circ(\Z(\bL)_\ell^F)=\bL$, hence
$\C_\bG^\circ(\Z(R))=\C_\bL^\circ(\Z(R))$. So $\bH=\C_\bG^\circ(\Z(R))$ and
$\bL=\C_\bG(\Z^\circ(\C^\circ_\bG(\Z(R)))_{\phi_e})$.
Therefore, $\N_\bG(R)\subseteq \N_\bG(\bL)$ and thus
$\Z(\bL)^F_\ell\unlhd \N_{\bG^F}(R)$. So the final assertion follows by the
arguments in the proof of Proposition~\ref{prop:center-Levi}(c).
\end{proof}	

\begin{lem}   \label{lem:cen-weights}
 Keep Condition~\ref{condition-1}. Let $\bL$ be an $e$-split Levi subgroup
 of~$\bG$.
 \begin{enumerate}[\rm(a)]
  \item Let $(R,\vhi)\in\Alp^0(\bG^F,\bL)$. Then
   $(R,\vhi_0)\in\Alp^0(\bL^F,\bL)$ for all
   $\vhi_0\in\Irr(\N_{\bL^F}(R)\mid\vhi)$.
  \item Let $(R,\vhi_0)\in\Alp^0(\bL^F,\bL)$. Then
   $(R,\vhi)\in\Alp^0(\bG^F,\bL)$ for all 
   $\vhi\in\rdz(\N_{\bG^F}(R)\mid\vhi_0)$.
 \end{enumerate}
\end{lem}

\begin{proof}
Part (a) follows by Proposition~\ref{prop:center-Levi}(c). To prove (b), let
$(R,\vhi_0)\in\Alp^0(\bL^F,\bL)$. By Corollary~\ref{cor:center-Levi},
$\bL=\C_\bG(\Z^\circ(\C^\circ_\bG(\Z(R)))_{\phi_e})$ and
$\N_{\bL^F}(R)\unlhd \N_{\bG^F}(R)$. So, if
$\vhi\in\rdz(\N_{\bG^F}(R)\mid\vhi_0)$, then $(R,\vhi)\in\Alp^0(\bG^F,\bL)$.
\end{proof}
 
In Lemma~\ref{lem:cen-weights} it follows that if $(R,\vhi_0)$ runs through
the weights in $\Alp^0(\bL^F,\bL)$ and $\vhi$ runs through
$\rdz(\N_{\bG^F}(R)\mid\vhi_0)$, then $(R,\vhi)$ runs through the weights in
$\Alp^0(\bG^F,\bL)$.

\begin{prop}   \label{prop:partition-weights}
 Keep Condition~\ref{condition-1}, and assume further that $\ell$ does not
 divide $|\Z(\bG_\SC)^F|$ and $\ell>3$ if the rational type of $(\bG,F)$
 includes type $^3\ty{D}_4$.
 If $(R,\vhi)\in\Alp^0(\bG^F,\bG)$, then $\Z(R)\subseteq \Z(\bG)$.
\end{prop}

\begin{proof}
Let $(R,\vhi)\in\Alp^0(\bG^F,\bG)$. Then
$\Z^\circ(\C_{\bG}^\circ(\Z(R)))_{\phi_e}\subseteq \Z(\bG)$. 
According to Proposition~\ref{prop:center-Levi}, $\C_{\bG}^\circ(\Z(R))$ is an
$E$-split Levi subgroup of $\bG$, and so by \cite[Thm.~22.2]{CE04}, if
$\C_\bG^\circ(\Z(R))<\bG$, then
$\C_\bG(\Z^\circ(\C_\bG^\circ(\Z(R)))_{\phi_e})<\bG$, a contradiction.
Thus $\C_{\bG}^\circ(\Z(R))=\bG$, which implies that $\Z(R)\subseteq \Z(\bG)$.
\end{proof}

\begin{cor}   \label{cor-partition-wei}
 Suppose that $\bG$ is simple of simply connected type, $\ell$ is odd and good
 for~$\bG$, does not divide $|\Z(\bG)^F|$ and $\ell>3$ if
 $\bG^F={}^3\ty{D}_4(q)$. Let $B$ be an $\ell$-block of $\bG^F$ and
 $e:=e_\ell(q)$. Then the set $\Alp^0(B,\bG)$ is non-empty if and only if $B$
 is of defect zero. In particular,
 $\Alp^0(\bG^F,\bG)=\{(1,\chi)\mid\chi\in\dz(\bG^F)\}.$
\end{cor}

\begin{proof}
The sufficiency is clear by	Remark~\ref{rmk:cen-def}.
Conversely, if $(R,\vhi)\in\Alp^0(B,\bG)$, then by
Proposition~\ref{prop:partition-weights}, $\Z(R)=1$ as $\ell\nmid |\Z(\bG)^F|$.
Hence $R=1$ and $\vhi\in\Irr(B)\cap\dz(\bG^F)$, and thus $B$ is of defect zero.
\end{proof}

%%%%%%%%%%%%%%%%%%%%%%%%%%%%%%%%%%%%%%%%%%%%%%%%%%%%%%%%%%%%%%%%%%%%%%%%%
%%%%%%%%%%%%%%%%%%%%%%%%%%%%%%%%%%%%%%%%%%%%%%%%%%%%%%%%%%%%%%%%%%%%%%%%%
\section{Groups of type $\ty{A}$}\label{sec:type-A}

In this section, we let $\tbG=\GL_n(\barFF_q)$ and $\bG=\SL_n(\barFF_q)$.
For any positive integer $k$ we denote by $F_{p^k}\colon\tbG\to\tbG$ the field
automorphism $(a_{ij})\mapsto(a_{ij}^{p^k})$, and by $\ga\colon\tbG\to\tbG$ the
graph automorphism $(a_{ij})\mapsto (a_{ji})^{-1}=((a_{ij}^{-1}))^{\tr}$ where
$^{\tr}$ denotes the transpose of matrices.
Let $\eps\in\{\pm1\}$ and $F=\ga^{\frac{1-\eps}{2}}F_q$.
Let $\tG=\tbG^F=\GL_n(\eps q)$ and $G=\bG^F=\SL_n(\eps q)$.
Here, by convention $\GL_n(-q)=\GU_n(q)$ and $\SL_n(-q)=\SU_n(q)$.
Let $\cB=\langle F_p,\ga\rangle$ or $\langle F_p\rangle$ according as $n\ge 3$
or $n=2$. Then $\tG\rtimes\cB$ induces all automorphisms of $G$; explicitly,
$(\tG\rtimes\cB)/\Z(\tG)\cong\Aut(G)$.
Let $\ell$ be a prime not dividing $q$ and $e:=e_\ell(q)$. Throughout this
Section \S\ref{sec:type-A}, we assume that $4\mid (q-\eps)$ when $\ell=2$.

The main aim of this section is the proof of the following theorem.

\begin{thm}   \label{thm:mainthm-type-A}
 Let $\tG'$ be the subgroup of $\tG$ such that $\tG'/G=(\tG/G)_\ell$.	
 There is a $(\tG\rtimes\cB)$-equivariant bijection 
 \[\vOm\colon\cW(\bG^F,\bG)\to\Alp_0(\bG^F)/\!\sim_{\tG'}\]
 such that 
 \begin{enumerate}[\rm(a)]	
  \item $\vOm(\cW(B,\bG))=\Alp_0(B)/\!\sim_{\tG'}$ for every $\ell$-block
   $B$ of $\bG^F$, and
  \item	for every $\chi\in\cW(\bG^F,\bG)$, there is a weight $(R,\vhi)$ of $G$
   whose $\tG'$-orbit corresponds to $\chi$ via $\vOm$ satisfying that
   $(\tG\rtimes\cB)_{R,\vhi}\subseteq (\tG\rtimes\cB)_\chi$ and that
   \[((\tG\rtimes\cB)_\chi,G,\chi)
     \geqslant_{(g),b}((\tG\rtimes\cB)_{R,\vhi},\N_G(R),\vhi)\]
   is normal with respect to $\N_{\tG'}(R)_\vhi$.
 \end{enumerate}	
\end{thm}

Here, the relation $\geqslant_{(g),b}$ between character triples was introduced
by the first author \cite{Fe25}, a generalization of the block isomorphism
$\geqslant_b$ of character triples first introduced by Navarro and
Sp\"ath \cite{NS14}.  In \cite{Sp17}, Sp\"ath reformulated the inductive
conditions of some of the local-global conjectures, including the Alperin
weight conjecture, in terms of central isomorphisms $\geqslant_c$ and block
isomorphisms $\geqslant_b$ between character triples.
We refer to \cite[Def.~3.6]{Fe25} for the definition of $\geqslant_{(g),b}$,
to \cite[Def.~3.14]{Fe25} for the notion of normality, and to
\cite[Def.~3.8]{Fe25} for the definition of $\geqslant_c$ and $\geqslant_b$.

%%%%%%%%%%%%%%%%%%%%%%%%%%%%%%%%%%%%%
\subsection{Characters and blocks of general linear and unitary groups}

Denote by $\Irr(\FF_q[x])$ the set of all non-constant monic irreducible
polynomials over $\FF_q$.
For $\Delta(x)=x^m+a_{m-1}x^{m-1}+\cdots+a_0$ in $\FF_{q^2}[x]$, we define
$\widetilde\Delta(x)=x^ma_0^{-q}\Delta^q(x^{-1})$, where $\Delta^q(x)$ means
the polynomial in $x$ whose coefficients are the $q$-th powers of the
corresponding coefficients of $\Delta(x)$. Now, we denote by
\begin{align*}
\cF_0 &= \left\{~ \Delta\mid \Delta\in\Irr(\FF_q[x]),\Delta\neq x ~\right\},\\
\cF_1 &= \left\{~ \Delta\mid \Delta\in\Irr(\FF_{q^2}[x]),
  \Delta\neq x,\Delta=\widetilde\Delta ~\right\},\\
\cF_2 &= \left\{~ \Delta\widetilde\Delta\mid
  \Delta\in\Irr(\FF_{q^2}[x]),\Delta\neq x,\Delta\neq\widetilde\Delta ~\right\}.
\end{align*}
Following \cite[\S 1]{FS82}, we let $\cF=\cF_0$ if $\eps=1$, and
$\cF=\cF_1\cup\cF_2$ if $\eps=-1$.
We denote by $\deg(\Ga)$ the degree of any polynomial~$\Ga$.

The conjugacy classes of semisimple elements of $\tbG^F$ can be classified in
terms of the polynomials in $\cF$. For any semisimple element $s$ of $\tbG^F$,
we let $s=\prod_{\Ga\in\cF} s_\Ga$ be its primary decomposition and let
$m_\Ga(s)$ denote the multiplicity of $\Ga$ in $s_\Ga$. If $m_\Ga(s)$ is
non-zero, then $\Ga$ is said to be an \emph{elementary divisor} of $s$.
Then $\C_{\tbG^F}(s)\cong\prod_{\Ga\in\cF}\GL_{m_\Ga(s)}((\eps q)^{\deg(\Ga)})$.
The unipotent characters of $\C_{\tbG^F}(s)$ can be labeled by the
combinatorial objects $\mu=\prod_{\Ga\in\cF}\mu_\Ga$ with
$\mu_\Ga\vdash m_\Ga(s)$ and the unipotent character of $\C_{\tbG^F}(s)$
corresponding to $\mu$ is denoted
$\tchi^\mu:=\prod_{\Ga\in\cF}\tchi^{\mu_\Ga}$, where $\tchi^{\mu_\Ga}$ is the
unipotent character of $\GL_{m_\Ga(s)}((\eps q)^{\deg(\Ga)})$ labeled by
$\mu_\Ga$. Then Lusztig's Jordan decomposition can be constructed by 
\[\cE(\tbL^F,1)\to\cE(\tbG^F,s),\quad
  \tchi^\mu \mapsto \tchi^{s,\mu}:=\pm\R_{\tbL}^{\tbG}(\widehat s\tchi^\mu),\]
where $\tbL:=\C_{\tbG}(s)$ and $\widehat s$ denotes the image of $s$ under the
isomorphism (see e.g. \cite[(8.19)]{CE04})
\begin{equation}\label{equ:iso-cen-linear}
\Z(\tbL)^F\to\Lin(\tbL^F)
\addtocounter{thm}{1}\tag{\thethm}
\end{equation}
which can be chosen as in \cite[p.~177]{Br86}.

For $\Ga\in\cF$ we denote by $d_\Ga$ the multiplicative order of
$(\eps q)^{\deg(\Ga)}$ modulo $\ell$. Let $\cF'$ be the subset of $\cF$ of
polynomials whose roots are of $\ell'$-order.
Then by \cite{FS82} the $\ell$-blocks of $\tbG^F$ are in bijection with the
$\tbG^F$-conjugacy classes of pairs $(s,\ka)$ with a semisimple $\ell'$-element
$s\in\tbG^F$ and $\kappa=\prod_{\Ga\in\cF}\ka_\Ga$ where $\ka_\Ga$ is the
$d_\Ga$-core of a partition of $m_\Ga(s)$. 
Moreover, if $\tB$ is an $\ell$-block of $\tbG^F$ with label $(s,\ka)$, then
the set $\Irr(\tB)\cap\cE(\bG^F,s)$ consists of characters $\tchi^{s,\mu}$ such
that $\ka_\Ga$ is the $d_\Ga$-core of $\mu_\Ga$ for every $\Ga\in\cF'$.
%See \cite[\S4]{Fe19} for an explicit description of $e$-Jordan-cuspidal pairs of $\bG^F$, which give a parametrization of the blocks of $\tbG^F$ by \cite[Thm.~A]{KM15}.

\begin{lem}   \label{thm:block-eJGC}
 Let $s\in\tbG^F$ be a semisimple $\ell'$-element and
 $\tB\subseteq\cE_\ell(\tbG^F,s)$ be an $\ell$-block of $\tbG^F$.
 Then $\cW(\tB,\tbG)$ is non-empty if and only if one of the following holds. 
 \begin{enumerate}[\rm(1)]
  \item $\tB$ is of defect zero, in which case $\ell\nmid(q-\eps)$ and
   $\cW(\tB,\tbG)$ consists of an $e$-Jordan-cuspidal character of $\tbG^F$; or
  \item $\ell\mid(q-\eps)$, $\tB=\cE_\ell(\tbG^F,s)$ and $s$ has exactly one
   elementary divisor, denoted $\Ga$. Moreover, $m_\Ga(s)$ is an $\ell$-power.
 \end{enumerate}
\end{lem}

\begin{proof}
Suppose that $s=\prod_\Ga s_\Ga$ is the primary decomposition of $s$ so that
\[\C_{\tbG^F}(s)\cong\prod_{\Ga\in\cF}\GL_{m_\Ga(s)}((\eps q)^{\deg(\Ga)}).\]

First assume (1) holds. Then it follows by Lemma~\ref{lem:def-zero-pairs} that
$\cW(\tB,\tbG)$ is non-empty.
If (2) holds, then by Proposition~\ref{lem:e-GC-type-A},
$\C_{\tbG^F}(s)\cong\prod_{\Ga\in\cF}\GL_{m_\Ga(s)}((\eps q)^{\deg(\Ga)})$
possesses a unipotent $e$-GC character. 
Moreover, $\Z(\C_{\tbG}(s))_{\phi_e}\subseteq \Z(\tbG)_{\phi_e}$, and thus
$\cW(\tB,\tbG)$ is non-empty.

On the other hand, if $\cW(\tB,\tbG)$ contains an $e$-Jordan-cuspidal character
of $\tbG^F$, then $\tB$ is of defect zero, by the explicit description for
$e$-Jordan-cuspidal pairs of $\tbG^F$ in \cite[\S4]{Fe19}, as in~(1). Now
assume that
$\cW(\tB,\tbG)$ possesses an $e$-JGC character $\tchi^{s,\mu}$ which is not
$e$-Jordan-cuspidal. Then $\C_{\tbG^F}(s)$ possesses a unipotent $e$-cuspidal
character $\tchi^\mu$ which is not $e$-cuspidal.
By Proposition~\ref{lem:e-GC-type-A} and its proof, for every $\Ga\in\cF'$,
$m_\Ga(s)$ is an $\ell$-power and $\ell$ divides $(\eps q)^{\deg(\Ga)}-1$.
So $\Z(\C_{\tbG}(s))_{\phi_e}$ is non-trivial. On the other hand, we conclude
from $\Z(\C_{\tbG}(s))_{\phi_e}\subseteq \Z(\tbG)_{\phi_e}$ that $e=1$ or $2$
according as $\eps=1$ or $-1$ and $s$ has exactly one elementary divisor, as
in~(2). This completes the proof.
\end{proof}	

\begin{rmk}\label{rmk:eJGC-char-GL}
 In the situation of Lemma~\ref{thm:block-eJGC}(2), we let $\Ga$ be the (unique)
 elementary divisor of $s$ and let $m=m_\Ga(s)$. Then by
 Proposition~\ref{lem:e-GC-type-A}, $\cW(\tB,\tbG)$ consists of characters
 $\tchi^{s,\mu}$ where $\mu_\Ga$ is one of the hook partitions
 $(m),(m-1,1),(m-2,1^2),\ldots,(1^m)$. In particular, $|\cW(\tB,\tbG)|=m$.
\end{rmk}

%%%%%%%%%%%%%%%%%%%%%%%%%%%%%%%%%%%%%
\subsection{Radical subgroups of general linear and unitary groups}

Now we recall the classification of the radical subgroups and weights of
$\tG=\GL_n(\eps q)$ by Alperin--Fong \cite{AF90} and An \cite{An92,An93,An94}.

Let $d$ be the multiplicative order of $\eps q$ modulo $\ell$. 
If $\eps=1$, then $d=e$, and if $\eps=-1$, then $d=2e$, $e/2$ or $e$ if $e$ is
respectively odd, congruent to~2 modulo~4, or divisible by~4. In particular,
$\phi_d(\eps x)=\pm\phi_e(x)$. Let $a$ be the precise power of $\ell$ dividing
$(\eps q)^d-1$, that is, $a$ is the integer with $((\eps q)^d-1)_\ell=\ell^a$.

\iffalse\textcolor{red}{
Let $\al,\ga\ge0$ be integers, $Z_\al$ be the cyclic group of order
$\ell^{a+\al}$ and $E_\ga$ be an extra-special $\ell$-group of order
$\ell^{2\ga+1}$ and exponent $\ell$. The group $Z_\al E_\ga$, which denotes the
central product of $Z_\al$ and $E_\ga$ over $\Omega_1(Z_\al)=\Z(E_\ga)$, can be
embedded into $\GL_{\ell^\ga}((\eps q)^{d\ell^\al})$ uniquely up to conjugacy
in the sense that $Z_\al$ is identified with
$\Z(\GL_{\ell^\ga}((\eps q)^{d\ell^\al}))_\ell$. 
We let $\tR_{\al,\ga}$ denote the image of $Z_\al E_\ga$ under the composition 
\[ Z_\al E_\ga \hookrightarrow \GL_{\ell^\ga}((\eps q)^{d\ell^\al})
  \hookrightarrow \GL_{d\ell^{\al+\ga}}(\eps q).\]
By \cite[(1C)]{An94}, $\tR_{\al,\ga}$ is defined uniquely up to conjugacy in
$\GL_{d\ell^{\al+\ga}}(\eps q)$ in the sense that $\Z(\tR_{\al,\ga})$ is
primary. Let $m$ be a positive integer. We denote
$\tR_{m,\al,\ga}=\tR_{\al,\ga}\otimes I_m$ where $I_m$ stands for the identity
matrix of degree $m$. Then $\tR_{m,\al,\ga}$ is an $\ell$-subgroup of
$\GL_{md\ell^{\al+\ga}}(\eps q)$. }\fi
Recall the construction of the radical $\ell$-subgroups of general linear and
unitary groups from \cite{AF90,An92,An93,An94}; see also \cite[\S5.3]{FZ22}.
Let $\tR_{m,\al,\ga}$ be the $\ell$-subgroup of
$\GL_{md\ell^{\al+\ga}}(\eps q)$, defined as in \cite[\S4]{FLZ21a}.

For a positive integer $c$ we denote by $A_c$ the elementary abelian
$\ell$-group of order $\ell^c$ in its regular permutation representation.
The group $A_c$ can be embedded uniquely up to conjugacy as a transitive
subgroup of the symmetric group $\fS_{\ell^c}$. For a sequence
$\fc=(c_1,\ldots,c_t)$ of positive integers, we write $l(\fc):=t$ and
$|\fc|:=c_1+\cdots+c_t$. The group $A_{\fc}=A_{c_1}\wr\cdots\wr A_{c_t}$ is
embedded uniquely up to conjugacy as a transitive subgroup of
$\fS_{\ell^{|\fc|}}$ where $|\fc|:=c_1+\cdots+c_t$.
Let $\tR_{m,\al,\ga,\fc}=\tR_{m,\al,\ga}\wr A_{\fc}$.
For convenience, we also write $\tR_{m,\al,\ga,\fc}$ for $\tR_{m,\al,\ga}$ with
$\fc=(0)$, in which situation we set $|\fc|=l(\fc)=0$. Following
\cite{AF90,An92,An93,An94}, we call the groups $\tR_{m,\al,\ga,\fc}$ \emph{basic
subgroups} of $\tG_{m,\al,\ga,\fc}=\GL_{md\ell^{\al+\ga+|\fc|}}(\eps q)$.
Denote the centralizer and normalizer of $\tR_{m,\al,\ga,\fc}$ in
$\tG_{m,\al,\ga,\fc}$ by $\tC_{m,\al,\ga,\fc}$ and $\tN_{m,\al,\ga,\fc}$
respectively.

First let $|\fc|=0$ and set
$\tN_{m,\al,\ga}^0=\C_{\tN_{m,\al,\ga}}(\Z(\tR_{m,\al,\ga}))$.
According to \cite[\S 3.A]{FLZ21a}, 
\[\tC_{m,\al,\ga}\cong \GL_m((\eps q)^{d\ell^\al})\otimes I_{\ell^\ga}\ 
  \textrm{and} \ \tN_{m,\al,\ga}^0/\tR_{m,\al,\ga}
  \cong\Sp_{2\ga}(\ell)\ti (\tC_{m,\al,\ga}\tR_{m,\al,\ga}/\tR_{m,\al,\ga}).\]
Here we interpret $\Sp_0(\ell)$ as the trivial group.
In addition, $\tN_{m,\al,\ga}=\tN_{m,\al,\ga}^0\rtimes\langle v\rangle$ where
$v$ is a permutation matrix of order $d\ell^\al$.

Now let $\fc=(c_1,\ldots,c_t)$ with $|\fc|>0$. Then by \cite[\S4]{AF90} and
\cite[\S2]{An94}, $\tC_{m,\al,\ga,\fc}\cong \GL_m((\eps q)^{d\ell^\al})\otimes I_{\ell^{\ga+|\fc|}}$ and 
\[ \tN_{m,\al,\ga,\fc}/\tR_{m,\al,\ga,\fc}\cong
\tN_{m,\al,\ga}/\tR_{m,\al,\ga}\ti \GL_{c_1}(\ell)\ti\cdots\ti\GL_{c_t}(\ell).\]
See \cite{FLZ21a} for explicit sets of generators of the above groups in
matrix form.

Let $\tR$ be a radical $\ell$-subgroup of $\tG=\GL_n(\eps q)=\GL(V)$ or
$\GU(V)$ according as $\eps=1$ or $-1$, where $V$ is the underlying space
of~$\tG$. By \cite[(4A)]{AF90} and \cite[(2B)]{An94}, there exist
decompositions 
\begin{equation}\label{equ:decompo-radical-0}
V=V_0\oplus V_1\oplus\cdots\oplus V_s\quad
   \text{and}\quad \tR=\tR_0\ti \tR_1\ti\cdots\ti \tR_s
  \addtocounter{thm}{1}\tag{\thethm}
\end{equation}
such that $\tR_0$ is the trivial subgroup of $\GL(V_0)$ or $\GU(V_0)$, and
$\tR_i$ is a basic subgroup of $\GL(V_i)$ or $\GU(V_i)$ for $i\ge 1$.
Let $V_+=V_1\oplus \cdots \oplus V_s$, $\tR_+= \tR_1\ti\cdots \ti \tR_s$ and
$\tG_+=\GL(V_+)$ or $\GU(V_+)$.
Then $\C_{\tG}(\tR)=\tC_0\times \tC_+$, $\N_{\tG}(\tR)=\tN_0\times \tN_+$ where
$\tC_0=\tN_0=\GL(V_0)$ or $\GU(V_0)$, $\tC_+=\C_{\tG_+}(\tR_+)$ and 
$\tN_+=\N_{\tG_+}(\tR_+)$.

\begin{lem}\label{lem:radical-small-center}
 Let $\tR$ be a radical $\ell$-subgroup of $\tbG^F$ with decomposition
 (\ref{equ:decompo-radical-0}). If there exists a weight $(\tR,\tvhi)$ in
 $\Alp^0(\tbG^F,\tbG)$, then one of the following holds.
 \begin{enumerate}[\rm(1)]
  \item $\tR$ is the trivial subgroup, in which case $\ell\nmid(q-\eps)$ and
   $\tvhi\in\dz(\tbG^F)$; or
  \item $\ell\mid(q-\eps)$ and $\tR$ is a basic subgroup of $\tbG^F$.
 \end{enumerate}
\end{lem}

\begin{proof}
Suppose that $\tR$ has the decomposition~(\ref{equ:decompo-radical-0}).
Then $\Z(\tR)=\prod_{i=0}^s\Z(\tR_i)$ and from
\[\Z(\C_{\tbG}(\Z(\tR)))_{\phi_e}\subseteq\Z(\tbG)\]
we conclude that $\tR=\tR_i$ for some $0\le i\le s$. If $\tR$ is the trivial
subgroup of $\tbG^F$, then $\tvhi$ is of defect zero and thus
$\ell\nmid(q-\eps)$. Otherwise, $\tR=\tR_i$ for some $i\ge 1$ so that it
is a basic subgroup of $\tbG^F$.
Moreover, $\Z(\tbG)_{\phi_e}$ is non-trivial, which implies $\ell\mid(q-\eps)$.
\end{proof}

If $\tR_{m,\al,\ga,\fc}$ provides a weight of $\tG_{m,\al,\ga,\fc}$ belonging
to a block $\tB\subseteq\cE_\ell(\tG_{m,\al,\ga,\fc},s)$ for some semisimple
$\ell'$-element $s$, then the group 
\[\tR_{m,\al,\ga,\fc}\tC_{m,\al,\ga,\fc}/\tR_{m,\al,\ga,\fc}
  \cong \GL_m((\eps q)^{d\ell^\al})/\Z(\GL_m((\eps q)^{d\ell^\al}))_\ell\]
possesses an irreducible character of $\ell$-defect zero, and by
\cite[\S4]{FS82} (see also \cite[\S5.A]{FLZ21a}) it follows that $\ell\nmid m$
and $s$ has exactly one elementary divisor. In particular, if $\ga=|\fc|=0$,
then $\tR_{m,\al,0}$ is a defect group of $\tB$.

%%%%%%%%%%%%%%%%%%%%%%%%%%%%%%%%%%%%%
\subsection{Weights of general linear and unitary groups}

By \cite{FS82}, given a polynomial $\Ga\in\cF'$, there exists a unique block
$\tB_\Ga$ of $\tG_\Ga=\GL_{d_\Ga \deg(\Ga)}(\eps q)$ with defect group
$R_\Ga=R_{m_\Ga,\al_\Ga,0}$, where $m_\Ga\ge 1$ and $\al_\Ga\ge 0$ are integers
with $d_\Ga \deg(\Ga)=m_\Ga d\ell^{\al_\Ga}$ and $\ell\nmid m_\Ga$.
Let $\tC_\Ga=\C_{\tG_\Ga}(\tR_\Ga)$ and $\tN_\Ga=\N_{\tG_\Ga}(\tR_\Ga)$.
Then $\tC_\Ga\cong \GL_{m_\Ga}((\eps q)^{d \ell^{\al_\Ga}})$ and
$\tN_\Ga/\tC_\Ga$ is cyclic and of order $d\ell^{\al_\Ga}$.
Let $\fb_\Ga$ be a root block of $\tB_\Ga$, i.e., a block of $\tC_\Ga$ with
defect group $\tR_\Ga$ and ${\fb_\Ga}^{\tG_\Ga}=\tB_\Ga$, and let
$\ttheta_\Ga$ denote the canonical character of $\fb_\Ga$. Up to
$\tN_\Ga$-conjugacy, $\fb_\Ga$ and $\ttheta_\Ga$ are uniquely determined
by~$\Ga$. In addition, the group $(\tN_\Ga)_{\ttheta_\Ga}/\tC_\Ga$ is cyclic
and has order $d_\Ga$.

Let $\tR_{\Ga,\ga,\fc}=\tR_{m_\Ga,\al_\Ga,\ga,\fc}$. 
We denote by $\tC_{\Ga,\ga,\fc}$ and $\tN_{\Ga,\ga,\fc}$ the centralizer and
normalizer of $\tR_{\Ga,\ga,\fc}$ in
$\tG_{\Ga,\ga,\fc}:=\GL_{d_\Ga \deg(\Ga) \ell^{\ga+|\fc|}}(q)$, respectively.
We have $\tC_{\Ga,\ga,\fc}\tR_{\Ga,\ga,\fc}/\tR_{\Ga,\ga,\fc}\cong \tC_\Ga/\tR_\Ga$ since $\tC_{\Ga,\ga,\fc}=\tC_\Ga\otimes I_{\ell^{\ga+|\fc|}}$. 
From this we can define the character
$\ttheta_{\Ga,\ga,\fc}:=\ttheta_\Ga\otimes I_{\ell^{\ga+|\fc|}}$ of $\tC_{\Ga,\ga,\fc}$ by $\ttheta_{\Ga,\ga,\fc}(c\otimes I_{\ell^{\ga+|\fc|}}):=\ttheta_\Ga(c)$ for $c\in \tC_\Ga$.

First let $|\fc|=0$, in which case we abbreviate $\tR_{\Ga,\ga,\fc}$,
$\tC_{\Ga,\ga,\fc}$, $\tN_{\Ga,\ga,\fc}$, $\ttheta_{\Ga,\ga,\fc}$ to
$\tR_{\Ga,\ga}$, $\tC_{\Ga,\ga}$, $\tN_{\Ga,\ga}$, $\ttheta_{\Ga,\ga}$,
respectively.
Let $\tN^0_{\Ga,\ga}:=\C_{N_{\Ga,\ga}}(\Z(\tR_{\Ga,\ga}))$ so that 
\[\tN^0_{\Ga,\ga}/\tR_{\Ga,\ga}\cong (\tC_{\Ga,\ga}\tR_{\Ga,\ga}/\tR_{\Ga,\ga})\ti \Sp_{2\ga}(\ell).\]
Then $\ttheta_{\Ga,\ga}$ is $\tN^0_{\Ga,\ga}$-invariant and
$\dz(\tN^0_{\Ga,\ga}/\tR_{\Ga,\ga}\mid\ttheta_{\Ga,\ga})=\{ \ttheta_{\Ga,\ga} \ze_\ga \}$, where $\ze_\ga$ is the Steinberg character of $\Sp_{2\ga}(\ell)$. 
We interpret $\ze_0$ as the trivial character for $\ga=0$.
The quotient $(\tN_{\Ga,\ga})_{\ttheta_{\Ga,\ga}}/\tN^0_{\Ga,\ga}\cong (\tN_\Ga)_{\ttheta_\Ga}/\tC_\Ga$ is cyclic of order $d_\Ga$.
Thus $\dz((\tN_{\Ga,\ga})_{\ttheta_{\Ga,\ga}}/\tR_{\Ga,\ga}\mid \ttheta_{\Ga,\ga})$ has $d_\Ga$ elements; they are exactly the extensions of $\ttheta_{\Ga,\ga}\ze_\ga$ to $(\tN_{\Ga,\ga})_{\ttheta_{\Ga,\ga}}$. 

Now let $\fc=(c_1,\ldots,c_t)$ with $|\fc|>0$. Then as before,
\[(\tN_{\Ga,\ga,\fc})_{\ttheta_{\Ga,\ga,\fc}}/\tR_{\Ga,\ga,\fc}\cong (\tN_{\Ga,\ga})_{\ttheta_{\Ga,\ga}}/\tR_{\Ga,\ga}\ti \GL_{c_1}(\ell)\ti\cdots\ti \GL_{c_t}(\ell).\]
For $i=1,\ldots,t$, the group $\GL_{c_i}(\ell)$ has $\ell-1$ defect zero irreducible characters, the extensions of the Steinberg character of $\SL_{c_i}(\ell)$.
So $\dz((\tN_{\Ga,\ga,\fc})_{\ttheta_{\Ga,\ga,\fc}}/\tR_{\Ga,\ga,\fc}\mid \ttheta_{\Ga,\ga,\fc})$ has $d_\Ga (\ell-1)^{l(\fc)}$ elements.

Let $\delta$ be a non-negative integer.
Set \[\cC_{\Ga,\delta}:=\{\, (\tR,\tpsi)\mid\tR=\tR_{\Ga,\ga,\fc},\ \tpsi\in\dz((\tN_{\Ga,\ga,\fc})_{\ttheta_{\Ga,\ga,\fc}}/\tR_{\Ga,\ga,\fc}\mid \ttheta_{\Ga,\ga,\fc}), \ \ga+|\fc|=\delta \,\}.\]
Then by \cite[\S4]{AF90} and \cite[\S4]{An94}, $\cC_{\Ga,\delta}$ has
cardinality $d_\Ga \ell^\delta$.

Let $\tR=\tR_{m,\al,\ga,\fc}$ is a basic subgroup of $\tG=\tG_{m,\al,\ga,\fc}$.
Suppose that $(\tR,\tvhi)$ is a weight of $\tG$. Then by \cite{AF90,An94},
there exists a polynomial $\Ga\in\cF'$ such that $m=m_\Ga$, $\al=\al_\Ga$, and
up to conjugacy $\tR=\tR_{\Ga,\ga,\fc}$,
$\tvhi=\Ind_{(\tN_{\Ga,\ga,\fc})_{\ttheta_{\Ga,\ga,\fc}}}^{\tN_{\Ga,\ga,\fc}}(\tpsi)$ for some
$\tpsi\in\rdz((\tN_{\Ga,\ga,\fc})_{\ttheta_{\Ga,\ga,\fc}}/\tR_{\Ga,\ga,\fc}\mid\ttheta_{\Ga,\ga,\fc})$.
If such a weight $(\tR,\tvhi)$ is a $\tB$-weight where $\tB$ is a block of
$\tG$, then $\tB\subseteq\cE_\ell(\tbG^F,s)$ for some semisimple
$\ell'$-element $s\in\tbG^F$ which has exactly one elementary divisor $\Ga$
with $m_\Ga(s)=d_\Ga\ell^{\ga+|\fc|}$.

\begin{lem}   \label{thm:block-weights-nonemp}
 Let $s\in\tbG^F$ be a semisimple $\ell'$-element and
 $\tB\subseteq\cE_\ell(\tbG^F,s)$ be an $\ell$-block of $\tbG^F$.
 Then $\Alp(\tB,\tbG)$ is non-empty if and only if one of the following holds.
 \begin{enumerate}[\rm(1)]
  \item $\tB$ is of defect zero, in which case $\ell\nmid (q-\eps)$; or
  \item $\ell\mid(q-\eps)$, $\tB=\cE_\ell(\tbG^F,s)$ and $s$ has exactly one
   elementary divisor, denoted by $\Ga$. Moreover, $m_\Ga(s)$ is an
   $\ell$-power.
 \end{enumerate}
\end{lem}

\begin{proof}
Note that if $\tB$ is of defect zero, then $\Alp(\tB,\tbG)$ is non-empty. 
Now we assume (2). Then $d_\Ga=1$ for every $\Ga\in\cF$.
Let $\ga$ be the integer satisfying $\ell^\ga=m_\Ga(s)$.
Then $\C_{\tbG}(\Z(\tR))_{\phi_e}\subseteq\Z(\tbG)$ for $\tR=\tR_{\Ga,\ga}$.
Now $(\tR,\Ind^{\tN_{\Ga,\ga}}_{(\tN_{\Ga,\ga})_{\ttheta_{\Ga,\ga}}}(\psi))$ is
a $\tB$-weight of $\tbG^F$ for
$\psi\in\dz((\tN_{\Ga,\ga})_{\ttheta_{\Ga,\ga}}/\tR_{\Ga,\ga}\mid\ttheta_{\Ga,\ga})$,
then $\Alp(\tB,\tbG)$ is non-empty. 

On the other hand, we assume that $\Alp(\tB,\tbG)$ is non-empty. Let
$(\tR,\tvhi)$ be a $\tB$-weight of $\tbG^F$ with
$\Z(\C_{\tbG}(\Z(\tR)))_{\phi_e}\subseteq\Z(\tbG)$.
By Lemma~\ref{lem:radical-small-center}, we have that $\tR$ is the trivial
subgroup of $\tbG^F$, $\tvhi\in\dz(\tbG^F)$ (which implies that
$\ell\nmid(q-\eps)$ and $\tB$ is of defect zero) or $\ell\mid(q-\eps)$,
$\tR$ is a basic subgroup. In the latter situation, we write
$\tR=\tR_{\Ga,\ga,\fc}$ and let $\ttheta\in\Irr(\C_{\tbG^F}(\tR)\mid\tvhi)$,
then up to $\N_{\tbG^F}(\tR)$-conjugacy $\ttheta=\ttheta_{\Ga,\ga,\fc}$ for
some $\Ga\in\cF'$, and thus $s$ has exactly one elementary divisor $\Ga$ and
$m_\Ga(s)$ is an $\ell$-power. This completes the proof.
\end{proof}

\begin{rmk}\label{rmk:weights-blocks-number}
Suppose that we are in the situation of
Lemma~\ref{thm:block-weights-nonemp}(2). Let $\delta$ be the integer with
$\ell^\delta=m_\Ga(s)$.
Then $\tbG^F=\tG_{\Ga,\ga,\fc}$ and the $\tB$-weights of $\tbG^F$ are these
$(\tR,\tvhi)$ (up to conjugacy): $\tR=\tR_{\Ga,\ga,\fc}$ and
$\tvhi=\Ind_{(\tN_{\Ga,\ga,\fc})_{\ttheta_{\Ga,\ga,\fc}}}^{\tN_{\Ga,\ga,\fc}}(\tpsi)$ with $\ga+|\fc|=\delta$ and $\tpsi\in\dz((\tN_{\Ga,\ga,\fc})_{\ttheta_{\Ga,\ga,\fc}}/\tR_{\Ga,\ga,\fc}\mid\ttheta_{\Ga,\ga,\fc})$.
So there is a bijection between $\Alp(\tB,\tbG)$ and $\cC_{\Ga,\delta}$.
In particular,
$|\Alp(\tB,\tbG)|=|\cC_{\Ga,\delta}|=\ell^\delta=m_\Ga(s)$.
\end{rmk}	

By Lemma~\ref{thm:block-eJGC}, Remark~\ref{rmk:eJGC-char-GL},
Lemma~\ref{thm:block-weights-nonemp} and
Remark~\ref{rmk:weights-blocks-number}, we have the following.

\begin{cor}   \label{cor:num-coincide}
 Let $\tB$ be an $\ell$-block of $\tbG^F$.
 \begin{enumerate}[\rm(a)]
  \item $\Alp(\tB,\tbG)$ is non-empty if and only if $\cW(\tB,\tbG)$ is
   non-empty.
  \item $|\cW(\tB,\tbG)|=|\Alp(\tB,\tbG)|$.
 \end{enumerate}
\end{cor}	

%%%%%%%%%%%%%%%%%%%%%%%%%%%%%%%%%%%%%
\subsection{An equivariant bijection}

We set $\fZ=\{ z\in\barFF_q^\ti\mid z^{q-\eps}=1\}$ and identify $\fZ$ with
$\Z(\tbG^F)$. Recall from (\ref{equ:iso-cen-linear}) that
$\fZ\to\Lin(\tbG^F/\bG^F)$, $z\mapsto\widehat z$, is an isomorphism; in fact,
$\widehat z$ is the character in $\cE(\tbG^F,z)$ corresponding under Jordan
decomposition to $1_{\tbG^F}$.

Let $\Ga\in\cF$ and $\mathrm{RT}_\Ga$ be the set of roots of $\Ga$ in
$\barFF_q$. For $z\in \fZ$, we define $z\Ga$ to be the polynomial in $\cF$ with
$\mathrm{RT}_{z\Ga}=\{ zx\mid x\in \mathrm{RT}_\Ga\}$.
Let $\sigma\in\cB$, then we define $\sigma(\Ga)$ to be the elementary divisor
of $\sigma(s_\Ga)$ where $s_\Ga$ is a semisimple element of
$\GL_{\deg(\Ga)}(\eps q)$ which has a unique elementary divisor $\Ga$ and that
has multiplicity~1.

Let $\tchi^{s,\mu}$ be an irreducible character of $\tbG^F$.
We define $\sigma(\mu)$ and $z\mu$ to be the combinatorial objects with
$\sigma(\mu)_{\sigma(\Ga)}=\mu_{\Ga}$ and $(z\mu)_{z\Ga}=\mu_\Ga$.
Then $(\tchi^{s,\mu})^\sigma=\tchi^{\sigma(s),\sigma(\mu)}$ and
$\widehat z\tchi^{s,\mu}=\tchi^{zs,z\mu}$.

Let $\tB$ be a block of $\tbG^F$. Suppose that $\ell\mid(q-\eps)$,
$\tB=\cE_\ell(\tbG^F,s)$ and $s$ has a unique elementary divisor, and that one
has multiplicity $\ell^\delta$ for some positive integer $\delta$.
By Remark~\ref{rmk:weights-blocks-number}, there is a bijection between
$\Alp(\tB,\tbG)$ and $\cC_{\Ga,\delta}$.
So the elements of $\Alp(\tB,\tbG)$ can be parameterized by combinatorial
objects $(s,w)$, where $w\in\cC_{\Ga,\delta}$.
Since $\cC_{\Ga,\delta}$ has cardinality $\ell^\delta$, we may write
$\cC_{\Ga,\delta}=\{ (\tR_{\Ga,\delta,i},\psi_{\Ga,\delta,i})\mid
  1\le i\le\ell^\delta \}$.
Write $w_{\Ga,\delta,i}=(\tR_{\Ga,\delta,i},\psi_{\Ga,\delta,i})$.
By \cite[(1C)]{An94}, the group $\cB$ acts trivially on $\Re(\tG)$.
For $\sigma\in\cB$, there exists $g_i\in\tG$ such that $g_i\sigma$ fixes
$R_{\Ga,\delta,i}$ and thus as in \cite[\S5]{LZ18} and \cite[\S5]{Fe19}, we can
choose the labeling of $\cC_{\Ga,\delta}$ such that
$\psi_{\Ga,\delta,i}^{g_i\sigma}=\psi_{\sigma(\Ga),\delta,i}$ and
$\widehat z\psi_{\Ga,\delta,i}=\psi_{z\Ga,\delta,i}$.
Therefore, by \cite[Prop.~5.3]{LZ18} and \cite[Prop.~5.12]{Fe19}, if
$\overline{(\tR,\tvhi)}\in\Alp(\tB,\tbG)$ has label $(s,w_{\Ga,\delta,i})$,
then $\overline{(\tR,\tvhi)}^{\sigma}\in\Alp(\tB^\sigma,\tbG)$ has label
$(\sigma(s),w_{\sigma(\Ga),\delta,i})$ for $\sigma\in\cB$ and
$\overline{(\tR,\widehat z\tvhi)}\in\Alp(\widehat z\otimes\tB,\tbG)$ has label
$(zs,w_{z\Ga,\delta,i})$ for $z\in\fZ_{\ell'}$.

\begin{thm}   \label{thm:type-A-equ-bijection}
 There exists a blockwise bijection between $\cW(\tbG^F,\tbG)$ and
 $\Alp(\tbG^F,\tbG)$ which is compatible with the action of
 $\Lin_{\ell'}(\tG/G)\rtimes\cB$. 
\end{thm}
	
\begin{proof}
The assertion follows if we prove that for every $\ell$-block $\tB$ of $\tbG^F$
the groups $\Lin_{\ell'}(\tG/G)_\tB$ and $\cB_\tB$ act trivially on
$\cW(\tB,\tbG)$ and $\Alp(\tB,\tbG)$, since those two sets have the same
cardinality by Corollary~\ref{cor:num-coincide}. If $\ell\nmid(q-\eps)$, then
by Lemmas~\ref{thm:block-eJGC} and~\ref{thm:block-weights-nonemp}, either
$\cW(\tB,\tbG)$ is empty or $\tB$ is of defect zero.
Thus we may assume that $\ell\mid(q-\eps)$, $\tB=\cE_\ell(\tbG^F,s)$ for some
semisimple $\ell'$-element $s$ of $\tbG^F$ and $s$ has exactly one elementary
divisor $\Ga$ and $m_\Ga(s)$ is an $\ell$-power. 
By the above arguments, we see that $\Lin_{\ell'}(\tG/G)_{\tB}$ and $\cB_{\tB}$
act trivially on $\cW(\tB,\tbG)$ and $\Alp(\tB,\tbG)$ and this gives the
assertion.
\end{proof}

\begin{lem}   \label{lem:number-block-char}
 Let $\tB$ be an $\ell$-block of $\tbG^F$ such that $\cW(\tB,\tbG)$ is
 non-empty. Let $\tchi\in\cW(\tB,\tbG)$. Then the number of blocks of $\bG^F$
 covered by $\tB$ equals the cardinality of $\Irr(\bG^F\mid\tchi)$.
\end{lem}

\begin{proof}
If $\tB$ is of defect zero, then this lemma holds as is easy to check.	
Now we assume that $\ell\mid (q-\eps)$ and $s$ has exactly one elementary
divisor $\Ga$ and $m_\Ga(s)$ is an $\ell$-power. So by Clifford theory,
$|\Irr(\bG^F\mid\tchi)|$ equals the number of elements $z\in\fZ_{\ell'}$ with
$z\Ga=\Ga$.
On the other hand, the number of blocks of $\bG^F$ covered by $\tB$ also equals
the number of elements $z\in\fZ_{\ell'}$ with $z\Ga=\Ga$ by
\cite[Rem.~4.13]{Fe19} and \cite[Rem.~6.9]{FLZ23}. This completes the proof.
\end{proof}

Now we prove the Main Theorem~\ref{thm:mainthm-type-A} of this section.

\begin{proof}[Proof of Theorem \ref{thm:mainthm-type-A}]
We prove this assertion by applying \cite[Thm.~5.1]{Fe25} by taking
$A=\tbG^F\rtimes \cB$, $\tG=\tbG^F$, $G=\bG^F$, $E=\cB$,
$\widetilde\cI=\cW(\tbG^F,\tbG)$ and $\widetilde\cA=\Alp(\tbG^F,\tbG)$. Note
that $\tB$ is taken to be the union of the blocks $\widetilde b$ of $\tbG^F$
with non-empty $\cW(\widetilde b,\tbG)$ here. The condition (ii) and (iii) of
\cite[Thm.~5.1]{Fe25} follows by \cite[Thm.~4.1]{CS17a} and
\cite[Thm.~7.1]{FLZ21a} respectively. Moreover, since $\tG/G$ is cyclic,
condition (i) of \cite[Thm.~5.1]{Fe25} holds automatically.
By Theorem~\ref{thm:type-A-equ-bijection}, it remains to verify condition
(iv.b) of \cite[Thm.~5.1]{Fe25}, which follows by
\cite[Prop.~5.6]{Fe25} and Lemma~\ref{lem:number-block-char}.
\end{proof}

%%%%%%%%%%%%%%%%%%%%%%%%%%%%%%%%%%%%%%%%%%%%%%%%%%%%%%%%%%%%%%%%%%%%%%%%%
%%%%%%%%%%%%%%%%%%%%%%%%%%%%%%%%%%%%%%%%%%%%%%%%%%%%%%%%%%%%%%%%%%%%%%%%%
\section{The inductive conditions for Alperin's weight conjecture}   \label{sec:IBAW}

In this section, we reformulate the inductive conditions for Alperin's weight
conjecture for groups of Lie type in terms of generic weights. 
Throughout the section $\bG$ is a simple algebraic group of simply connected
type over $\barFF_p$. Let $\Phi$ and $\Delta$ denote respectively the set of
roots and simple roots of $\bG$ determined by the choice of a maximal torus and
a Borel subgroup containing it. To describe Frobenius endomorphisms of~$\bG$,
we use the Chevalley generators $x_\alpha(t)$
($t\in\barFF_q$, $\alpha\in\Phi$) as in~\cite[Thm.~1.12.1]{GLS98}.

Recall the endomorphisms of $\bG$ described as in~\cite[\S2]{MS16}.
Let $F_0:\bG\to\bG$ denote the field endomorphism of $\bG$ given by
$F_0(x_\alpha(t))=x_\alpha(t^p)$ for $t\in\barFF_p$ and $\al\in\Phi$.
A (length-preserving) automorphism $\tau$ of the Dynkin diagram associated to
$\Delta$ (and hence an automorphism of $\Phi$) determines a graph automorphism
$\gamma$ of $\bG$ given by $\gamma(x_\alpha(t)):=x_{\tau(\alpha)}(t)$ for
$t\in\barFF_p$ and $\al\in \pm\Delta$. Any such $\gamma$ commutes with $F_0$.

Suppose that $\Z(\bG)$ has rank $r$ as finite abelian group.
Let $\bZ$ be a torus of rank~$r$ with an embedding of
$\Z(\bG)$. Let us set $\tbG:=\bG\times_{\Z(\bG)}\bZ$ the central product of
$\bG$ and $\bZ$ over $\Z(\bG)$. Then $\tbG$ is a connected reductive group such
that the natural map $\bG\hookrightarrow\tbG$ is a regular embedding.
As in \cite[p.~874]{MS16}, we can extend $F_0$ to a Frobenius endomorphism of
$\tbG$ and $\gamma$ to an automorphism of $\tbG$.

Consider a Frobenius endomorphism $F:=F_0^f\gamma$, with $f$ a positive integer
and $\gamma$ a (possibly trivial) graph automorphism of $\bG$.
%, leaving aside the types $^2{\ty{B}}_2$, $^2{\ty{G}}_2$ and $^2{\ty{F}}_4$.
Then $F$ defines an $\FF_q$-structure on $\tbG$, where $q=p^f$.
The groups of rational points $G=\bG^F$ and $\tG=\widetilde\bG^F$ are
finite. %~(see~\cite[Thm.~21.5]{MT11}).
Let $\cB$ be the subgroup of $\Aut(\bG^F)$ generated by $F_0$ (here we identify
$F_0$ with $F_0\rceil_G$) and the graph automorphisms commuting with $F$.
Then $\tbG^F\rtimes \cB$ is well defined and induces all automorphisms of
$\bG^F$ (see~\cite[Thm.~2.5.1]{GLS98}). 
Let $\Diag(\bG^F)$ be the subgroup of $\Aut(\bG^F)$ induced by $\tbG^F$ and let
$\Diag_\ell(\bG^F)$ be the subgroup of $\Diag(\bG^F)$ induced by $\tG'$ where
$\bG^F\le\tG'\le\tbG^F$ with $\tG'/\bG^F=(\tbG^F/\bG^F)_\ell$.
Note that if $R$ is a radical $\ell$-subgroup of $\bG^F$ and
$\bT:=\Z^\circ(\C^\circ_\bG(\Z(R)))_{\phi_e}$, then
$\N_{\tbG^F\rtimes\cB}(R)\subseteq \N_{\tbG^F\rtimes\cB}(\bT)$.

\begin{condition}   \label{condition-simply-connnected}
	Suppose that $\bG$ is simple and simply connected and $F\colon\bG\to\bG$ is
	a Frobenius endomorphism with respect to an $\FF_q$-structure. 
	Let $\ell$ be an odd prime not dividing $q$ such that $\ell$ is good for
	$\bG$ and $\ell\nmid |\Z(\bG)^F|$. Assume that $\ell>3$ if
	$\bG^F={}^3\ty{D}_4(q)$. Let $e:=e_\ell(q)$.
	Let $\tbG$ and $\cB$ be defined as above. 
\end{condition}

In Condition~\ref{condition-simply-connnected}, from $\ell\nmid |\Z(\bG)^F|$ we
deduce that $\ell$ divides none of $|\cZ(\bG)^F|$, $|\cZ(\bG)_F|$ or
$|\cZ(\bG^*)^F|$ (in fact, $\Z(\bG^*)^F=1$).

We will prove the following theorem in this section.

\begin{thm}   \label{main-thm-weights}
 Keep Condition \ref{condition-simply-connnected}. Let $B$ be an $\ell$-block
 of $\bG^F$. Then there is a $(\tbG^F\rtimes \cB)_B$-equivariant bijection
 $\vOm\colon \cW(B)\to\Alp(B)$ such that for every
 $\overline{(\bT,\eta)}\in\cW(B)$, there exists a $B$-weight $(R,\vhi)$ of
 $\bG^F$ with $\overline{(R,\vhi)}=\vOm(\overline{(\bT,\eta)})$ satisfying
 \begin{enumerate}[\rm(a)]
  \item $\bT=\Z^\circ(\C^\circ_\bG(\Z(R)))_{\phi_e}$,
   $\bl(\vhi)^{\N_{\bG^F}(\bT)}=\bl(\eta)$ and
  \item $((\tbG^F\rtimes\cB)_{\bT,\eta},\N_{\bG^F}(\bT),\eta)
   \geqslant_b((\tbG^F\rtimes \cB)_{R,\vhi},\N_{\bG^F}(R),\vhi)$.
 \end{enumerate}
\end{thm}

%%%%%%%%%%%%%%%%%%%%%%%%%%%%%%%%%%%%%
\subsection{The proof of Theorem~\ref{main-thm-weights}}
Theorem~\ref{main-thm-weights} follows from the following theorem by letting
$\bT$ run through a representative set of the $\bG^F$-conjugacy classes of
$e$-tori of $\bG$.

\begin{thm}   \label{thm:main-theorem-single-block}
 Keep Condition \ref{condition-simply-connnected}. Let $\bT$ be an $e$-torus
 of~$\bG$ with $\bT=\Z^\circ(\C_\bG(\bT))_{\phi_e}$. Then there is a
 $(\tbG^F\rtimes \cB)_\bT$-equivariant bijection
 \[\vOm_\bT\colon \cW^0(\bG^F,\bT)
   \to\Alp^0(\bG^F,\bT)/\sim_{\N_{\bG^F}(\bT)}\]
 such that for every $\eta\in\cW^0(\bG^F,\bT)$, there exists a weight
 $(R,\vhi)$ in $\Alp^0(\bG^F,\bT)$ whose $\N_{\bG^F}(\bT)$-orbit corresponds to
 $\eta$ via $\vOm_\bT$ satisfying $\bl(\vhi)^{\N_{\bG^F}(\bT)}=\bl(\eta)$
 and
 \[((\tbG^F\rtimes \cB)_{\bT,\eta},\N_{\bG^F}(\bT),\eta)
  \geqslant_b((\tbG^F\rtimes \cB)_{R,\vhi},\N_{\bG^F}(R),\vhi).\]
\end{thm}	

\begin{proof}
Let $\bL:=\C_\bG(\bT)$. By \cite[Prop.~2.4]{GH91} and \cite[Prop.~13.12]{CE04},
$\ell$ is good for $\bL$, and does not divide the orders of $\cZ(\bL)_F$,
$\cZ(\bL)^F$ and $\cZ(\bL^*)^F$. Since $\bL$ is a Levi subgroup of $\bG$, by
\cite[Prop.~12.14]{MT11}, $[\bL,\bL]$ is a simply connected semisimple
algebraic group, that is,
\[ [\bL,\bL]=\bH_1\ti\cdots\ti\bH_s,\]
where for every $1\le i\le s$, $\bH_i$ is a simply connected simple algebraic
group. Also, the rational types of $(\bL,F)$ do not include type $^3\ty{D}_4$
when $\ell=3$. Indeed, for $\bL\ne\bG$ to have a component of type $^3\ty{D}_4$
the group $\bG$ has to be of exceptional type, but then $\ell=3$ is bad for
$\bG$, contrary to assumption. Moreover, $[\bL,\bL]\embed \bL$ is an
$\ell$-regular embedding as $|\cZ(\bL)_F|$ is prime to $\ell$, and so
$\Alp^0(\bL^F,\bL)/\sim_{\bL^F}=\Alp(\bL^F\mid\Alp_0([\bL,\bL]^F))$ (as defined
in \S\ref{subsec:preli-weights}).
By a similar proof as for \cite[Prop.~6.3]{FS21}, one shows
\[\cE(\bL^F,\ell')=
 \Irr(\bL^F\mid\cE([\bL,\bL]^F,\ell'))\cap \Irr(\bL^F\mid 1_{\Z(\bL^F)_\ell}).\]
Thus according to Lemma~\ref{lem:Clifford-eJGC},
\[\cW(\bL^F,\bL)=\Irr(\bL^F\mid\cW([\bL,\bL]^F,[\bL,\bL]))
  \cap\Irr(\bL^F\mid 1_{\Z(\bL^F)_\ell}).\]
Here $\cW(\bL^F,\bL)$ is defined as in \S\ref{subsec:generic-weights}.

Since $\bL=[\bL,\bL]\Z(\bL)$ and $\C_\bH(\bH^F)=\C_\bH(\bH)=\Z(\bH)$ for any
semisimple group $\bH$, we have
$$C_\bL([\bL,\bL]^F)
  =\C_{[\bL,\bL]}([\bL,\bL]^F)\Z(\bL)=\Z([\bL,\bL])\Z(\bL)=\Z(\bL)$$
and hence
$$\Z(\bL^F)\le\C_{\bL^F}([\bL,\bL]^F)=\Z(\bL)^F$$
are equal.

The action of the Frobenius endomorphism $F$ induces a permutation $\sigma$ on
the set $\{\bH_1\,\ldots,\bH_s\}$ and we decompose
$\sigma=\sigma_1\cdots\sigma_t$ into disjoint cycles.
For $1\le i\le t$, let $\Sigma_i$ be the support of the permutation $\sigma_i$
and let $n_i=|\Sigma_i|$. Then the inclusion map
$\bH_{k_i}\embed \prod_{j\in\Sigma_i}\bH_j$ induces an isomorphism
$\bH_{k_i}^{F^{n_i}}\cong (\prod_{j\in\Sigma_i}\bH_j)^F$ for any
$k_i\in\Sigma_i$ (in the following we fix one $k_i$ in every $\Sigma_i$) for
every $1\le i\le t$. Thus we have
\[[\bL,\bL]^F= \prod_{i=1}^t\bigg(\prod_{j\in\Sigma_i}\bH_j\bigg)^F\cong
  \prod_{i=1}^t \bH_{k_i}^{F^{n_i}}.\]
Write $H_{k_i}:=(\prod_{j\in\Sigma_i}\bH_j)^F$.

Let $(R_0,\vhi_0)$ be a weight of $[\bL,\bL]^F$. Then
$R_0=R_{0,1}\ti \cdots\ti R_{0,t}$ and
$\vhi_0=\vhi_{0,1}\ti\cdots\ti \vhi_{0,t}$ where $(R_{0,i},\vhi_{0,i})$ is a
weight of $H_{k_i}$ for $1\le i\le t$. %We write
%\[(R_0,\vhi_0)=:(R_{0,1},\vhi_{0,1})\ti\cdots\ti (R_{0,t},\vhi_{0,t}).\]

Note that $[\bL,\bL]=[\bL,\bL]_\ba\ti\bL_\bb$ (in the sense of
Notation~\ref{notation:ab}), $[\bL,\bL]^F=[\bL,\bL]_\ba^F\ti \bL_\bb^F$ and
$\ell\nmid |\Z(\bL_\bb)^F|$ (by Condition~\ref{condition-simply-connnected}).
Let $[\bL,\bL]_\ba\to\tbL_\ba$ be an $\ell$-regular embedding.
Then $[\bL,\bL]\to\tbL_\ba\ti\bL_\bb$ is also $\ell$-regular, and
$\Alp_0([\bL,\bL]^F)$ consists of the conjugacy classes of weights of
$[\bL,\bL]^F$ covered by the elements in
\[\Alp(\tbL_{\ba}^F\ti\bL_\bb^F,\tbL_\ba\ti\bL_\bb)
  =\Alp(\tbL_\ba^F,\tbL_\ba)\ti\Alp(\bL_\bb^F,\bL_\bb).\]
So $\Alp_0([\bL,\bL]^F)=\Alp_0([\bL,\bL]_\ba^F)\ti \Alp(\bL_\bb^F,\bL_\bb)$,
and by Lemma~\ref{lem:weights-ell'center},
$\Alp(\bL_\bb^F,\bL_\bb)=\Alp_0(\bL_\bb^F)$.
By Corollary~\ref{cor-partition-wei}, if
$(R,\vhi)\in\Alp^0(\bL_\bb^F,\bL_\bb)$, then $R=1$ and thus
$\vhi\in\dz(\bL_\bb^F)$. 
By Corollary~\ref{cor-par-gen-wei}, $\cW(\bL_\bb^F,\bL_\bb)=\dz(\bL_\bb^F)$.
Therefore, for $1\le i\le t$, according to \cite[Thm.~4.3]{Fe25} and
Theorem~\ref{thm:mainthm-type-A}, there exists a blockwise
$\Aut(\bH_{k_i}^{F^{n_i}})$-equivariant bijection 
\[f_{k_i}\colon \cW(\bH_{k_i}^{F^{n_i}},\bH_{k_i})\to
  \Alp_0(\bH_{k_i}^{F^{n_i}})/\!\sim_{\Diag_\ell(\bH_{k_i}^{F^{n_i}})}\]
such that for any character $\ze_{0,k_i}\in \cW(\bH_{k_i}^{F^{n_i}},\bH_{k_i})$,
there exists a weight $(R_{0,k_i},\vhi_{0,k_i})$ of $H_{k_i}$ whose
$\Diag_\ell(H_{k_i})$-orbit corresponds to $\ze_{0,k_i}$ via $f_{k_i}$ and
satisfies that
\[ (H_{k_i}\rtimes \Aut(H_{k_i})_{\ze_{0,k_i}},H_{k_i},\ze_{0,k_i})
  \geqslant_{(g),b}( (H_i\rtimes \Aut(H_{k_i}))_{R_{0,k_i},\vhi_{0,k_i}}, \N_{H_{k_i}}(R_{0,k_i})_{\vhi_{0,k_i}},\vhi_{0,k_i} ) \]
is normal with respect to
$\N_{H_{k_i}\rtimes \Diag_\ell(H_{k_i})}(R_{0,i})_{\vhi_{0,i}}$.

Let $L'$ be the subgroup of $\bL^F$ containing $[\bL,\bL]^F$ such that
$L'/[\bL,\bL]^F=(\bL^F/[\bL,\bL]^F)_\ell$.
According to Lemma~\ref{lem:diga-ell'}, the automorphisms induced by $L'$ on
$[\bL,\bL]^F$ form $\Diag_{\ell}([\bL,\bL]^F)$. Define the map
\[\vOm_{\bT,0}:\cW([\bL,\bL]^F,[\bL,\bL])\to\Alp_0([\bL,\bL]^F)/\!\sim_{L'}\]
by 
\[\ze_{0,k_1}\ti\cdots\ti\ze_{0,k_t}\mapsto
  f_{k_1}(\ze_{0,k_1})\ti\cdots\ti f_{k_t}(\ze_{0,k_t}),\]
where $\ze_{0,k_i}\in \cW(\bH_{k_i}^{F^{n_i}},\bH_{k_i})$ for $1\le i\le t$.

Let $\{ k_1,\ldots,k_t \}=A_1\cup\cdots\cup A_u$ be the partition such that
$k_j,k_l\in A_i$ if and only if $n_j=n_l$ and there exists an isomorphism of
algebraic groups $\bH_{k_j}\to\bH_{k_l}$ commuting with the action of $F^{n_j}$.
For each $1\le i\le u$ we fix a representative $y_{i}\in A_i$.
Thus we can identify $[\bL,\bL]^F$ with $\prod_{i=1}^u H_{y_i}^{|A_i|}$, and so
\[\Aut([\bL,\bL]^F)\cong\prod_{i=1}^u \Aut(H_{y_i})\wr\fS_{|A_i|}.\]

Let $c_0:=c_{0,k_1}\ti\cdots\ti c_{0,k_t}$ be a block of $[\bL,\bL]^F$, where $c_{0,k_i}$ is a block of $H_{k_i}$. 
For any $1\le i\le u$ we define a partition $A_i=I_{i,1}\cup\cdots\cup I_{i,w_i}$ such that for $k_j,k_l\in A_i$ we have $k_j,k_l\in I_{i,j}$ if and
only if $c_{0,k_j}=c_{0,k_l}$, under the induced isomorphism
$H_{k_j}\cong H_{k_l}$.
For each $1\le i\le u$, $1\le j\le w_i$ we fix a representative $z_{i,j}\in I_{i,j}$.
Without loss of generality, we may assume that $c_0=\otimes_{i=1}^{u}\otimes_{j=1}^{w_i}c_{0,z_{i,j}}$, and thus \[\Aut([\bL,\bL]^F)_{c_0}=\prod_{i=1}^{u}\prod_{j=1}^{w_i}\Aut(H_{z_{i,j}})_{c_{0,z_{i,j}}}\wr \fS_{|I_{i,j}|}.\]
It can be checked directly that the stabilizer of an irreducible character in $c_0$ (resp.\ of a $c_0$-weight of $[\bL,\bL]^F$) is also a direct product of wreath products, and from this one checks that the bijection $\vOm_{\bT,0}$ is blockwise and $\Aut([\bL,\bL]^F)$-equivariant.

According to Theorems~4.1 and~4.2 of \cite{Fe25}, for every
$\ze_0\in\cW([\bL,\bL]^F,[\bL,\bL])$, there is a weight $(R_0,\vhi_0)$ of
$[\bL,\bL]^F$ whose $L'$-orbit corresponds to $\ze_0$ via
$\vOm_{\bT,0}$ satisfying
$(\tbG^F\rtimes\cB)_{\bT,R_0,\vhi_0}\subseteq (\tbG^F\rtimes\cB)_{\bT,\ze_0}$
and 
\[((\tbG^F\rtimes\cB)_{\bT,\ze_0},[\bL,\bL]^F,\ze_0)\geqslant_{(g),b}
  ((\tbG^F\rtimes\cB)_{\bT,R_0,\vhi_0},\N_{[\bL,\bL]^F}(R_0)_{\vhi_0},\vhi_0)\]
is normal with respect to $\N_{L'}(R_0)_{\vhi_0}$.

By \cite[Thm.~5.8]{Fe25}, we can obtain a blockwise
$(\tbG^F\rtimes\cB)_\bT$-equivariant bijection 
\[\vOm_\bT\colon\cW^0(\bL^F,\bT)\to\Alp^0(\bL^F,\bT)/\sim_{\bL^F}\]
such that for every $\ze\in\cW^0(\bL^F,\bT)$, there is a weight $(R,\vhi)$ of
$\bL^F$ with $\vOm_\bT(\ze)=\overline{(R,\vhi)}$ satisfying
$(\tbG^F\rtimes\cB)_{R,\vhi}\subseteq (\tbG^F\rtimes\cB)_{\bT,\ze}$ and
\[((\tbG^F\rtimes\cB)_{\bT,\ze},\bL^F,\ze)
  \geqslant_{b}((\tbG^F\rtimes\cB)_{R,\vhi},\N_{\bL^F}(R)_{\vhi},\vhi).\]

Recall that
\[\cW^0(\bG^F,\bT)
 =\{ \eta\in\rdz(\N_{\bG^F}(\bT)\mid \ze)\mid\ze\in\cW^0(\bL^F,\bT) \}\]
and by Lemma~\ref{lem:cen-weights}, 
\[\Alp^0(\bG^F,\bT)=\{ (R,\vhi')\mid (R,\vhi)\in\Alp^0(\bL^F,\bT),
   \vhi'\in\rdz(\N_{\bG^F}(R)\mid\vhi) \}.\] 
Therefore, using the arguments in the proof of \cite[Prop.~4.7]{NS14} we can
conclude.
\end{proof}

%%%%%%%%%%%%%%%%%%%%%%%%%%%%%%%%%%%%%
\subsection{Criteria for the inductive conditions}
We can now reformulate the inductive Alperin weight (AW) condition from
\cite{NT11} and the inductive blockwise Alperin weight (BAW) condition from
\cite{Sp13} in terms of generic weights.

\begin{lem}   \label{lem:bijection-tbG}
 Keep Condition \ref{condition-simply-connnected}. Then there exists a
 $(\Lin_{\ell'}(\tbG^F/\bG^F)\rtimes\cB)$-equivariant bijection
 \[\widetilde\vOm\colon \cW(\tbG^F)\to\Alp(\tbG^F)\] such that 
 \begin{enumerate}[\rm(a)]
  \item 
  \begin{enumerate}[\rm(1)]
   \item $\widetilde\vOm(\cW(\tB))=\Alp(\tB)$ for every $\ell$-block $\tB$
    of $\tbG^F$, 
   \item for $(\bT,\eta)\in\cW^0(\bG^F)$, one has
    $\widetilde\vOm(\cW(\tbG^F\mid\overline{(\bT,\eta)}))
    =\Alp(\bG^F\mid\vOm(\overline{(\bT,\eta)}))$, where $\vOm$ is the bijection
    from Theorem~\ref{main-thm-weights}, and
  \end{enumerate}
  \item for every $\overline{(\widetilde\bT,\widetilde\eta)}\in\cW(\tbG^F)$
   there exists a weight $(\tR,\tvhi)$ of $\tbG^F$ with
   $\overline{(\tR,\tvhi)}=\vOm(\overline{(\widetilde\bT,\widetilde\eta)})$
   satisfying
  \begin{enumerate}[\rm(1)]
   \item $\widetilde\bT=\Z^\circ(\C^\circ_{\tbG}(\Z(\tR)))_{\phi_e}$, and
   \item $((\tbG^F\rtimes \cB)_{\widetilde\bT,\widetilde\eta},\N_{\tbG^F}(\widetilde\bT),\widetilde\eta)\geqslant_b((\tbG^F\rtimes\cB)_{\tR,\tvhi},\N_{\tbG^F}(\tR),\tvhi)$.
  \end{enumerate}
 \end{enumerate}		 
\end{lem}

\begin{proof}
Since $\cZ(\bG)_F\cong\tbG^F/\bG^F\Z(\tbG)^F$ is an $\ell'$-group, the proof of
\cite[Prop.~5.2]{Fe19} shows that the map
$\Re(\bG^F)\to\Re(\tbG^F)$, $R\mapsto R\Z(\tbG)^F_\ell$, is bijective.
Let $\bT$ be an $e$-torus of $\bG$, $R$ a radical $\ell$-subgroup of
$\bG^F$, $\widetilde\bT:=\Z(\tbG)\bT$ and $\tR:=R\Z(\tbG)^F_\ell$.
Then $\widetilde\bT=\Z^\circ(\C^\circ_{\tbG}(\Z(\tR)))_{\phi_e}$ if and only if
$\bT=\Z^\circ(\C^\circ_\bG(\Z(R)))_{\phi_e}$.
So this lemma follows from Theorem~\ref{main-thm-weights} and
\cite[Thm.~5.8]{Fe25}.
\end{proof}

\begin{thm}   \label{thm:criterion-AW}
 Keep Condition~\ref{condition-simply-connnected}. Assume that
 $S:=\bG^F/\Z(\bG^F)$ is simple and does not have an exceptional covering group.
 Suppose that $\cE(\bG^F,\ell')$ is a uni-triangular basic set for $\bG^F$ and
 the following conditions are satisfied:	
 \begin{enumerate}[\rm(1)]
  \item \begin{enumerate}[\rm(a)]
	\item $\tbG^F/\bG^F$ is abelian and
     $\C_{\tbG^F\rtimes\cB}(\bG^F)=\Z(\tbG)^F$,
	\item $\cB$ is abelian or isomorphic to the direct product of a cyclic
     group with the symmetric group $\fS_3$, 
	\item every character in $\cE(\bG^F,\ell')$ extends to its stabilizer in
     $\tbG^F$, and
	\item for every $(\bT,\eta)\in\cW^0(\bG^F)$, the character $\eta$ extends
     to its stabilizer in $\N_{\tbG^F}(\bT)$.
  \end{enumerate}
  \item For every $\tchi\in\cE(\tbG^F,\ell')$, there exists some character
   $\chi_0\in\Irr(\bG^F\mid\tchi)$ such that $(\tbG^F\rtimes \cB)_{\chi_0}
    =\tbG^F_{\chi_0}\rtimes \cB_{\chi_0}$ and $\chi_0$ extends to
   $\bG^F\rtimes \cB_{\chi_0}$.
  \item For every $(\widetilde\bT,\widetilde\eta)\in\cW^0(\tbG^F)$, there
   exists $\eta_0\in\Irr(\N_{\bG^F}(\widetilde\bT)\mid\widetilde\eta)$ such that
   $(\tbG^F\rtimes \cB)_{\widetilde\bT,\eta_0} = \tbG^F_{\widetilde\bT,\eta_0} (\bG^F\rtimes \cB)_{\widetilde\bT,\eta_0}$ and	$\eta_0$ extends to $(\bG^F\rtimes \cB)_{\widetilde\bT,\eta_0}$.
  \item There exists a $(\Lin_{\ell'}(\tbG^F/\bG^F)\rtimes\cB)$-equivariant
   bijection
   \[\widetilde\vOm\colon \cE(\tbG^F,\ell') \to \cW(\tbG^F)\]
   such that for every $\widetilde\nu \in \Lin_{\ell'}(\Z(\tbG)^F)$ and every
   $\tchi\in\Irr(\tbG^F\mid\widetilde\nu)\cap\cE(\tbG^F,\ell')$, we have
   $\widetilde\vOm(\tchi)=\overline{(\tbT,\widetilde\eta)}$ with
   $\widetilde\eta\in\Irr(\N_{\tbG^F}(\tbT)\mid\widetilde\nu)$.
 \end{enumerate}
 Then the inductive AW condition holds for the simple group $S$ and the prime
 $\ell$.
\end{thm}

\begin{proof}
Here we use \cite[Thm.~3.3]{BS22}.
By the fact that $\cZ(\bG)_F\cong\tbG^F/\bG^F\Z(\tbG)^F$ is an $\ell'$-group,
we know condition (ii.2) of \cite[Thm.~3.3]{BS22} is satisfied.
Note that $\cE(\tbG^F,\ell')\subseteq\Irr(\tbG^F\mid 1_{\Z(\tbG)^F_\ell})$.
Since $\cE(\bG^F,\ell')$ is an $\Aut(\bG^F)$-stable uni-triangular basic set for
$\bG^F$, there exists a $(\tbG^F\rtimes\cB)$-equivariant bijection
$\varrho\colon\cE(\bG^F,\ell')\to\IBr(\bG^F)$ such that for every
$\chi\in\cE(\bG^F,\ell')$ and every $\bG^F\le H\le (\tbG^F\rtimes\cB)_\chi$,
if $\chi$ extends to $H$, then $\varrho(\chi)$ extends to $H$ (see, e.g.,
\cite[Lemma~2.9]{FLZ19}).
By \cite[Prop.~6.3]{FS21}, $\cE(\tbG^F,\ell')$ is a uni-triangular basic set for
$\tbG^F$ (of course it is $(\Lin_{\ell'}(\tbG^F/\bG^F)\rtimes\cB)$-stable) and
thus there exists a $(\Lin_{\ell'}(\tbG^F/\bG^F)\rtimes\cB)$-equivariant
bijection $\widetilde\varrho\colon\cE(\tbG^F,\ell')\to\IBr(\tbG^F)$.
In addition, for all $\chi\in\cE(\bG^F,\ell')$, we have
$\widetilde\varrho (\Irr(\tbG^F\mid\chi)\cap \cE(\tbG^F,\ell'))
 =\IBr(\tbG^F\mid\varrho(\chi))$.

Therefore, (1) (resp. (2), (3), (4)) implies condition (i) (resp. (iii), (iv),
(ii)) of \cite[Thm.~3.3]{BS22} by the above paragraph,
Lemma~\ref{lem:bijection-tbG} and its proof, and Theorem~\ref{main-thm-weights}.
So the inductive AW condition holds for the simple group $S$ and the
prime~$\ell$.
\end{proof}

\begin{rmk}   \label{rmk:condition-i-ii}
We notice that conditions (1.a) and (1.b) of Theorem~\ref{thm:criterion-AW}
hold. Moreover, (1.c) of Theorem~\ref{thm:criterion-AW} follows by \cite{Lu88},
while in \cite[Thm.~4.1]{CS17a}, \cite[Thm.~3.1]{CS17b}, \cite[Thm.~B]{CS19}
and \cite[Thm.~A]{Sp23}, condition (2) of Theorem~\ref{thm:criterion-AW} is
proved.
\end{rmk}	

\begin{cor}   \label{thm-iAW-geqslant_c}
 Keep Condition \ref{condition-simply-connnected}.  Assume that
 $\cE(\bG^F,\ell')$ is a uni-triangular basic set for $\bG^F$ and that
 $S=\bG^F/\Z(\bG^F)$ is simple and does not have an exceptional covering group.
 Suppose that there exists a $(\tbG^F\rtimes\cB)$-equivariant bijection
 \[\vOm\colon\cE(\bG^F,\ell')\to\cW(\bG^F)\]
 such that for every $\chi\in\cE(\bG^F,\ell')$ and
  $\vOm(\chi)=\overline{(\bT,\eta)}$, one has
 \[((\tbG^F\rtimes\cB)_\chi,\bG^F,\chi)
   \geqslant_c((\tbG^F\rtimes\cB)_{\bT,\eta},\N_{\bG^F}(\bT),\eta).\]
 Then the inductive AW condition holds for the simple group $S$ and the
 prime~$\ell$.
\end{cor}	
		
\begin{proof}
By Remark~\ref{rmk:condition-i-ii}, (1.a)--(1.c) and (2) of
Theorem~\ref{thm:criterion-AW} are satisfied. Moreover, by our assumption and
Theorem~\ref{main-thm-weights}, conditions~(1.d) and~(3) hold. In addition,
condition~(4) also follows from the assumptions, so
Theorem~\ref{thm:criterion-AW} gives the result.
\end{proof}

\begin{thm}   \label{thm:criterion-BAW}
 Keep Condition \ref{condition-simply-connnected}. Let $B$ be a union of
 $\ell$-blocks of $\bG^F$ which is a $\tbG^F$-orbit and $\tB$ be the union of
 blocks of $\tG$ covering $B$. Assume that $\Irr(B)\cap\cE(\bG^F,\ell')$ is a
 uni-triangular basic set of $B$ and the following hold.
 \begin{enumerate}[\rm(1)]
	\item \begin{enumerate}[\rm(a)]
	\item $\tbG^F/\bG^F$ is abelian and
     $\C_{\tbG^F\rtimes\cB}(\bG^F)=\Z(\tbG)^F$,
	\item $\cB$ is abelian or isomorphic to the direct product of a cyclic
     group with $\fS_3$, 
	\item every character in $\Irr(B)\cap\cE(\bG^F,\ell')$ extends to its
     stabilizer in $\tbG^F$, and
	\item for every $(\bT,\eta)\in\cW^0(B)$, the character $\eta$ extends to
     its stabilizer in $\N_{\tbG^F}(\bT)$.
  \end{enumerate}
  \item For every $\tchi\in\Irr(\tB)\cap\cE(\bG^F,\ell')$, there exists some
     character $\chi_0\in\Irr(\bG^F\mid\tchi)$ such that
     $(\tbG^F\rtimes \cB)_{\chi_0}=\tbG^F_{\chi_0}\rtimes \cB_{\chi_0}$ and
     $\chi_0$ extends to $\bG^F\rtimes \cB_{\chi_0}$.
  \item For every $(\widetilde\bT,\widetilde\eta)\in\cW^0(\tB)$, there exists
   $\eta_0\in\Irr(\N_{\bG^F}(\widetilde\bT)\mid\widetilde\eta)$ such that
   $(\tbG^F\rtimes \cB)_{\widetilde\bT,\eta_0}
      = \tbG^F_{\widetilde\bT,\eta_0} (\bG^F\rtimes \cB)_{\widetilde\bT,\eta_0}$
   and $\eta_0$ extends to $(\bG^F\rtimes \cB)_{\widetilde\bT,\eta_0}$.
  \item There exists a $(\Lin_{\ell'}(\tbG^F/\bG^F)\rtimes \cB_{\tB})$-equivariant bijection
   \[\widetilde\vOm\colon \Irr(\tB)\cap\cE(\tbG^F,\ell') \to \cW(\tB)\]
   such that
   \begin{enumerate}[\rm(a)]
	 \item $\widetilde\vOm$ preserves blocks, and
	 \item if the character $\tchi$ in (2) and
	  $(\tbT,\widetilde\eta)\in\cW^0(\tB)$ in (3) satisfy
	  $\overline{(\tbT,\widetilde\eta)}=\widetilde\vOm(\tchi)$, then the
	  character $\chi_0$ in (2) and $(\bT,\eta_0)\in\cW^0(B)$ in (3) can be
	  chosen in the same block of $\bG^F$, and to satisfy that
	  $\bl(\hchi)=\bl(\widehat\eta)^{\tbG^F_\chi}$, where
	  $\hchi\in\Irr(\tbG^F_{\chi_0}\mid\chi_0)$ is the Clifford correspondent
	  of $\tchi$ and $\widehat\eta\in\Irr(\N_{\tbG^F}(\bT)_{\eta_0}\mid\eta_0)$
	  is the Clifford correspondent of $\widetilde\eta$.
	\end{enumerate}		
 \end{enumerate}
 Then the inductive BAW condition holds for every block in $B$.
\end{thm}
	
\begin{proof}		
We use \cite[Thm.~2.2]{FLZ23} and the proof is similar to the one of
Theorem~\ref{thm:criterion-AW}.
\end{proof}

\begin{rmk}
If $\Out(\bG^F)$ is abelian and conditions (1), (2), (3) and (4.a) in
Theorem~\ref{thm:criterion-BAW} hold, then by a similar argument as above and
using \cite[Thm.~4.5]{BS22}, we can prove that the inductive BAW
condition holds for every block in $B$.
\end{rmk}	

\begin{cor}   \label{cor:new-criterion-BAW}
 Keep Condition \ref{condition-simply-connnected}. 
 Let $B$ be an $\ell$-block of $\bG^F$. Assume that
 $\Irr(B)\cap\cE(\bG^F,\ell')$ is a uni-triangular basic set of $B$ and there
 exists a $(\tbG^F\rtimes\cB)_B$-equivariant bijection
 \[\vOm\colon\Irr(B)\cap\cE(\bG^F,\ell')\to\cW(B)\]
 such that for every $\chi\in\Irr(B)\cap\cE(\bG^F,\ell')$ and
 $\vOm(\chi)=\overline{(\bT,\eta)}$, one has
 \[((\tbG^F\rtimes\cB)_\chi,\bG^F,\chi)
   \geqslant_b((\tbG^F\rtimes\cB)_{\bT,\eta},\N_{\bG^F}(\bT),\eta).\]
 Then the inductive BAW condition holds for $B$.
\end{cor}	
		
\begin{proof}
We use Theorem~\ref{thm:criterion-BAW}, and the proof is just similar to the
one for Theorem~\ref{thm-iAW-geqslant_c}.		
\end{proof}

%%%%%%%%%%%%%%%%%%%%%%%%%%%%%%%%%%%%%%%%%%%%%%%%%%%%%%%%%%%%%%%%%%%%%%%%%
%%%%%%%%%%%%%%%%%%%%%%%%%%%%%%%%%%%%%%%%%%%%%%%%%%%%%%%%%%%%%%%%%%%%%%%%%
\section{Character correspondences for relative Weyl groups}   \label{sec:open-problem}

According to Lemma~\ref{lem-genwei-dz}, we may hope to use the defect zero
characters of relative Weyl groups to describe the generic weights, as we
expect Assumption~\ref{ext-char-e-cuspidal} to often hold.
On the other hand, under Condition~\ref{condition-simply-connnected}, the
Brauer characters in a block are in bijection with the characters of a relative
Weyl group, if the generalized $e$-Harish-Chandra theory is known to hold.
From this, we propose a question for character correspondences at the level of
relative Weyl groups.

Throughout this section, we assume that $\bG$ is connected reductive, $F$ a
Frobenius map with respect to an $\FF_q$-structure, and $\ell$ is odd, good for
$\bG$ and does not divide $q|\cZ(\bG)^F|$. Let $e=e_\ell(q)$.

%%%%%%%%%%%%%%%%%%%%%%%%%%%%%%%%%%%%%%%%%%
\subsection{A question}

Let $B$ be an $\ell$-block of $\bG^F$. Under
Assumption~\ref{ext-char-e-cuspidal},
\[|\cW(B)|=\sum\limits_{(\bL,\ze)}|\dz(W_{\bG^F}(\bL,\ze))|,\]
where $(\bL,\ze)$ runs through the $\bG^F$-conjugacy classes of $e$-JGC pairs
of $\bG$ with $\ze\in\cE(\bL^F,\ell')$ and $B=\RLG(\bl(\ze))$.

If generalized $e$-Harish-Chandra theory holds, we have a bijection
\[\Irr(B)\cap\cE(\bG^F,\ell')\to \Irr(W_{\bG^F}(\bL_0,\ze_0)),\]
where $(\bL_0,\ze_0)$ is an $e$-Jordan-cuspidal pair corresponding to $B$ by
\cite[Thm.~A(e)]{KM15} so that $\ze_0\in\cE(\bL_0^F,\ell')$ and
$\R^\bG_{\bL_0}(\bl(\ze_0))=B$. 
By \cite[Thm.~14.4]{CE04}, $\Irr(B)\cap\cE(\bG^F,\ell')$ is a basic set for $B$.

In the spirit of Theorem~\ref{main-thm-weights} and the Alperin weight
conjecture, we propose the following question:

\begin{question}\label{ques:partition-local}
 Let $(\bL_0,\ze_0)$ be an $e$-Jordan-cuspidal pair of $\bG$ with
 $\ze_0\in\cE(\bL_0^F,\ell')$. Is there a bijection
  \begin{equation}\label{equ:partition-local}
    \Irr(W_{\bG^F}(\bL_0,\ze_0))\to \coprod\limits_{(\bL,\ze)}\dz(W_{\bG^F}(\bL,\ze))
    \addtocounter{thm}{1}\tag{\thethm}
  \end{equation}
  where $(\bL,\ze)$ runs through the $\bG^F$-conjugacy classes of $e$-JGC pairs
  of $\bG$ with $\ze\in\cE(\bL^F,\ell')$ and
  $\R^\bG_{\bL_0}(\bl(\ze_0))=\RLG(\bl(\ze))$?
\end{question}

In what follows, we will prove that bijections (\ref{equ:partition-local})
exist for all unipotent blocks, as well as for quasi-isolated blocks of
exceptional groups.

%%%%%%%%%%%%%%%%%%%%%%%%%%%%%%%%%%%%%%%%%%%%%%
\subsection{Bijections (\ref{equ:partition-local}) for unipotent blocks}

First, for groups with abelian Sylow $\ell$-subgroups and blocks with abelian
defect groups, the correspondence (\ref{equ:partition-local}) is just the
identity map, so we have:

\begin{lem}\label{lem:biject-abel-sylow}
 Suppose that $\bG$ is simple and $\bG^F$ has abelian Sylow $\ell$-subgroups.
 Then there is a bijection as in~(\ref{equ:partition-local}) for all
 $\ell$-blocks of $\bG^F$.
\end{lem}

\begin{proof}
This follows from Lemma~\ref{cor:abel-sylow} and the fact that
$W_{\bG^F}(\bL_0,\ze_0)$ is an $\ell'$-group (cf. \cite[Prop.~2.4]{Ma14}).
\end{proof}

\begin{lem}
 Assume that $\bG$ is an $F$-stable Levi subgroup of a simple algebraic group
 $\bH$ of simply connected type with a Frobenius endomorphism extending $F$.
 Let $\ell$ be a prime not dividing $q$ such that $\ell$ is odd and good
 for~$\bG$, and with $\ell>3$ if $\bG^F={}^3\ty{D}_4(q)$. If $B$ is an
 $\ell$-block of $\bG^F$ with abelian defect groups, then there is a bijection
 as in~(\ref{equ:partition-local}) for $B$.
\end{lem}

\begin{proof}
This follows from (a) and (b) of Proposition~\ref{prop:generic-abeldef}.
\end{proof}

In the following, we will show:

\begin{thm}   \label{thm:unipotent}
 Assume that $\ell$ is odd, good for $\bG$ and does not divide $|\cZ(\bG)^F|$. 
 Then there is a bijection (\ref{equ:partition-local}) for all unipotent
 $\ell$-blocks of $\bG^F$.
\end{thm}

We start with exceptional groups.

\begin{prop}   \label{prop:bijection exc}
 Let $\bG$ be of exceptional type and $\ell$ a good prime for $\bG$. Then there
 is a bijection as in~(\ref{equ:partition-local}) for all unipotent
 $\ell$-blocks of $\bG^F$.
\end{prop}

\begin{table}[htb]
\caption{Principal blocks for good primes $\ell$ at $e=1$}   \label{tab:exc}
$\begin{array}{c|cccc}
  \bG^F& \ell& |W_{\bG^F}(\bL_0,1)|& \sum|\dz(W_{\bG^F}(\bL,\ze))|\\
\hline
 \ty{E}_6(q)& 5&  25& 15+2\cdot5\\
 \ty{E}_7(q)& 5&  60& 30+6\cdot5\\
 \ty{E}_7(q)& 7&  60& 46+2\cdot7\\
 \ty{E}_8(q)& 7& 112& 84+4\cdot7\\
\end{array}$
\end{table}

\begin{proof}
If the Sylow $\ell$-subgroups of $\bG^F$ are abelian, this was seen in
Lemma~\ref{lem:biject-abel-sylow}.
Thus (e.g. by \cite[Thm.~25.14]{MT11}) we are left with the case that $\ell=5$
and $\bG$ is of type $\ty{E}_6$ or $\ty{E}_7$, or $\ell=7$ and $\bG$ is of type
$\ty{E}_7$ or $\ty{E}_8$, and moreover $e:=e_\ell(q)\in\{1,2\}$ in all cases.
Then the only unipotent block $B$ of $\bG^F$ with non-abelian defect is the
principal block
(by the description of defect groups in \cite{CE94}). First assume $e=1$.
Then $B$ is labelled by the $1$-cuspidal pair $(\bL_0,1)$ with $\bL_0\le\bG$
the centralizer of a Sylow 1-torus. By Lemma~\ref{lem:noncus-GC}, if
$(\bL,\ze)\ne(\bL_0,1)$ is $1$-GC with $\RLG(\bl(\ze))=B$ then $\bL$ has
a single component of type $\ty{A}_{\ell-1}(q)$ (by rank considerations), and
each such has $\ell$ unipotent $e$-GC characters. The corresponding
cardinalities are listed in Table~\ref{tab:exc} by which our claim follows.
Note that the Sylow 5-subgroups of $^2\ty{E}_6(q)$ with $e=1$ are abelian.
The case $e=2$ is entirely analogous.
\end{proof}

Curiously, the very same numbers as in column~4 of Table~\ref{tab:exc} appeared
in \cite[Tab.~3]{KMS}, originating in the associated fusion systems.

Observe that in the situation of Proposition~\ref{prop:bijection exc} the
characters in $W_{\bG^F}(\bL_0,1)$ not of defect zero are of height~0, since
$W_{\bG^F}(\bL_0,1)$ has cyclic Sylow $\ell$-subgroups, and thus
are in relation with the unipotent height zero characters of
$\ty{A}_{\ell-1}(q)$ by a McKay bijection. We also have the analogue for
quasi-isolated blocks:

\begin{prop}   \label{prop:bijection exc qi}
 Let $\bG$ be of exceptional type and $\ell$ a good prime for $\bG$. Then there
 is a bijection as in~(\ref{equ:partition-local}) for all quasi-isolated
 $\ell$-blocks of $\bG^F$.
\end{prop}

\begin{table}[htb]
\caption{Quasi-isolated $\ell$-blocks for good primes $\ell$ at $e=1$}   \label{tab:exc qi}
$\begin{array}{c|cccc}
 \bG^F& \ell& \C_{\bG^{*F}}(s)& |W_{\bG^F}(\bL_0,1)|& \sum|\dz(W_{\bG^F}(\bL,\ze))|\\
\hline\hline
 \ty{E}_6(q)& 5& \ty{A}_5(q)\ty{A}_1(q)& 22& 2\cdot6+2\cdot5\\
\hline
 \ty{E}_7(q)& 5& \ty{A}_5(q)\ty{A}_2(q)& 33& 3\cdot6+3\cdot5\\
 \ty{E}_7(q)& 5& \ty{D}_6(q)\ty{A}_1(q)& 74& 2\cdot27+2\cdot10\\
 \ty{E}_7(q)& 5& \ty{A}_7(q).2& 44& 14+6\cdot5\\
 \ty{E}_7(q)& 5& \Phi_1.\ty{E}_6(q).2& 50& 30+4\cdot5\\
 \ty{E}_7(q)& 7& \ty{A}_7(q).2& 44& 30+2\cdot7\\
\hline
 \ty{E}_8(q)& 7& \ty{A}_8(q)& 30& 16+2\cdot7\\
 \ty{E}_8(q)& 7& \ty{D}_8(q)& 100& 86+2\cdot7\\
 \ty{E}_8(q)& 7& \ty{E}_7(q)\ty{A}_1(q)& 120& 92+4\cdot7\\
 \ty{E}_8(q)& 7& \ty{A}_7(q)\ty{A}_1(q)& 44& 30+2\cdot7\\
\end{array}$
\end{table}

\begin{proof}
Let $B$ be a quasi-isolated $\ell$-block of $\bG^F$. If the defect groups of
$B$ are abelian, this was seen in Lemma~\ref{lem:biject-abel-sylow}. If $B$ is
unipotent, see Proposition~\ref{prop:bijection exc}. Arguing as there, we have
$\ell\in\{5,7\}$ and $\bG$ of type $\ty{E}_n$. Moreover, if $B$ is labeled by
the $e$-cuspidal pair $(\bL,\la)$ then $|W_{\bG^F}(\bL,\la)|$ is divisible
by~$\ell$. According to the tables in \cite{KM13} the only cases coming up
are those given in Table~\ref{tab:exc qi}. Here, $s\in G^*$ is a semisimple
$\ell'$-element such that $\Irr(B)\subseteq\cE_\ell(G,s)$. We can now argue as
before.
\end{proof}

Observe that by \cite[Prop.~6.5]{KMS} all characters in
$\Irr(B)\cap\cE(\bG^F,\ell')$ corresponding to defect zero characters of one
fixed relative Weyl group (so to one $e$-Harish-Chandra series) have the same
height.

Now we consider correspondences (\ref{equ:partition-local}) for unipotent
blocks of classical groups. That is, we need to deal with the groups $\fS_n$,
$C_e\wr\fS_n$ (with $\ell\nmid e$) and $G(2e,2,n)$ (with $\ell\nmid 2e$).

We first look at $\fS_n$, the symmetric group on $\{1,2,\ldots,n\}$. Let $\ell$
be a prime. For a non-negative integer~$n$, if 
\begin{equation}\label{equ-ell-expansion}
	n=\sum_{i\ge 0} \beta_i \ell^i
	\addtocounter{thm}{1}\tag{\thethm}
\end{equation}
for some non-negative integers $\beta_i$, then we say that
(\ref{equ-ell-expansion}) is an \emph{$\ell$-expansion} of $n$. 
If $\beta_i<\ell$ for every $i\ge 0$, then it is called the \emph{$\ell$-adic
expansion} of $n$.

Let $\nu_\ell$ be the exponential valuation associated to the prime $\ell$,
normalized so that $\nu_\ell(\ell)=1$. For finite groups $H\le G$ we abbreviate
$\nu_\ell(|G:H|)$ to $\nu_\ell(G:H)$. In particular, $\nu_\ell(G)$ stands for
$\nu_\ell(|G|)$. If $\chi\in\Irr(G)$, then we denote by $\de(\chi)$ the defect
of~$\chi$, that is, $\de(\chi)=\nu_\ell(G)-\nu_\ell(\chi(1))$. In addition, we
denote by $\Irr_{\ell'}(G)$ the set of irreducible characters of $G$ of degree
prime to $\ell$.

Recall that the partitions of $n$ are in bijection with the conjugacy classes
of Young subgroups of $\fS_n$. Let $\mu=(\mu_1,\ldots,\mu_k)\vdash n$.
Then the corresponding \emph{Young subgroup} of $\fS_n$ is
\[\fS_{\mu}=\fS_{\{1,2,\ldots,\mu_1\}}\ti\fS_{\{\mu_1+1,\mu_1+2,\ldots,\mu_1+\mu_2\}}\ti\cdots\ti\fS_{\{n-\mu_k+1,n-\mu_k+2,\ldots,n\}}.\]
We have $\fS_{\mu}\cong\fS_{\mu_1}\ti\cdots\ti\fS_{\mu_k}$.

\begin{lem}   \label{lem:exten-linear}
 Let $G=N\rtimes H$ and let $\la$ be a $G$-invariant linear character of $N$.
 Then $\hat\la$ defined by $\widehat\la(nh):=\la(n)$, for $n\in N$, $h\in H$,
 is an extension of $\la$ to $G$.
\end{lem} 

\begin{proof}
This can be checked directly.
\end{proof}

\begin{lem}  \label{lem-ell'char}
 Assume that $n=\ell^m$ with $m\ge 0$ and $e$ is an integer with $\ell\nmid e$.
 Then we have the following.
 \begin{enumerate}[\rm(a)]
  \item The characters of $\Irr_{\ell'}(\fS_{n})$ are labeled by hook
   partitions.
  \item The set $\Irr_{\ell'}(C_e\wr\fS_{n})$ consists of $en$ characters.
 \end{enumerate}
\end{lem}	
 
\begin{proof}
(a) follows from \cite[\S3 and \S4]{Ma71}, and can also be proved by direct
calculation using the hook formula.
 
For (b), let $G=C_e\wr\fS_{n}$, $\chi\in\Irr_{\ell'}(G)$ and
$\theta\in\Irr(C_e^n\mid\chi)$. Then $\chi=\Ind^{G}_{H}(\hat\theta)$ where
$H=G_\theta$ and $\hat\theta\in\Irr(H\mid\theta)$.
Since $\chi$ is of $\ell'$-degree, we see that $\ell\nmid[G:H]$.
Also note that $H=C_e^n\rtimes H'$ where $H'$ is a Young subgroup of $\fS_n$
and $n$ is an $\ell$-power, which force that $H'=\fS_n$.
Moreover, $\theta=\tau^{\boxtimes n}$ where $\tau\in\Irr(C_e)$.

Conversely, if $\theta=\tau^{\boxtimes n}\in\Irr(C_e^n)$ where
$\tau\in\Irr(C_e)$, then $\theta$ extends to $G$ by
Lemma~\ref{lem:exten-linear}.
So Gallagher’s theorem implies a bijection between $\Irr_{\ell'}(\fS_n)$ and
$\Irr_{\ell'}(G\mid\theta)$, and thus (b) holds. 
\end{proof}	

An \emph{$\ell$-Young subgroup} of $\fS_n$ is a conjugate of some Young
subgroup $\fS_{\mu}$ such that $\mu=(\mu_1,\ldots,\mu_k)\vdash n$
where $\mu_i$ is an $\ell$-power for any $1\le i\le k$.
If $Y=C_e\wr Y'$ where $Y'$ is an $\ell$-Young subgroup of $\fS_n$, then we say
that $Y$ is an \emph{$\ell$-Young subgroup} of $C_e\wr\fS_n$.
Now the complex reflection group $G(2e,2,n)$ is a semidirect product
$N\rtimes\fS_n$ where $N$ is a certain
subgroup of $C_{2e}^n$ of index 2. If $Y=N\rtimes Y'$ where $Y'$ is an
$\ell$-Young subgroup of $\fS_n$, then we say that $Y$ is an \emph{$\ell$-Young
subgroup} of $G(2e,2,n)$.
For any $\ell$-expansion $n=\sum\limits_{i\ge 0}\beta_i \ell^i$ we can define a
partition $\mu\vdash n$ such that $\ell^i$ appears $\beta_i$ times for every
$i\ge 0$, so that
$\fS_\mu\cong\prod\limits_{i\ge 0}(\fS_{\ell^i})^{\beta_i}$.
This induces a bijection between the $\ell$-expansions of $n$ and the conjugacy
classes of $\ell$-Young subgroups of $\fS_n$, $C_e\wr\fS_n$ or $G(2e,2,n)$.

Let $G=\fS_n$, $C_e\wr\fS_n$ or $G(2e,2,n)$.
If $Y$ is an $\ell$-Young subgroup of $G$ and $\ze\in\Irr_{\ell'}(Y)$, then we
call $(Y,\ze)$ an \emph{$\ell$-Young pair} of $G$.
Let $\cY_G$ be the set of $\ell$-Young pairs of $G$.

\begin{prop}   \label{prop-bijection-sym}
 Let $G=\fS_n$, $C_e\wr\fS_n$ (with $\ell\nmid e$) or $G(2e,2,n)$ (with
 $\ell\nmid 2e$). There is a bijection between $\Irr(G)$ and the $G$-conjugacy
 classes of triples $(Y,\ze,\la)$ where $(Y,\ze)\in \cY_G$ and
 $\la\in\dz(\N_G(Y)_\ze/Y)$ such that $\de(\chi)=\nu_\ell(Y)$ whenever
 $\chi\in\Irr(G)$ corresponds to $(Y,\ze,\la)$.
\end{prop}

\begin{proof}
First we let $G=\fS_n$, in which situation the proof can be found in the proof
of \cite[Prop.~(4.9)]{MO83} and we recall it as follows. Recall (from, for
example, \cite[p.~29]{Ol76}) that the partitions of $n$ are in natural
bijection with the $\ell$-core towers $(\ka_{i,j})_{i\ge 0,\ 1\le j\le\ell^i}$
with $n=\sum\limits_{i\ge 0}\sum\limits_{j=1}^{\ell^i}|\ka_{i,j}|\ell^i$.
Let $\chi\in\Irr(G)$ be a character labeled by $\mu\vdash n$.
Let $(\ka_{i,j})_{i\ge 0,\ 1\le j\le\ell^i}$ be the $\ell$-core tower of $\mu$.
Write $\beta_i(\ell,\mu)=\sum\limits_{j=1}^{\ell^i} |\ka_{i,j}|$ so that
$n=\sum\limits_{i\ge 0}\beta_i(\ell,\mu)\ell^i$.
By \cite[(3.3)]{Ma71}, we have
\[ \de(\chi)=\frac{n-\sum\limits_{i\ge 0}\beta_i(\ell,\mu)}{\ell-1}.\]

For every $i\ge0$, we write $\Irr_{\ell'}(\fS_{\ell^i})=\{\xi_{i,j}\mid 1\le j\le \ell^i\}$ by Lemma~\ref{lem-ell'char}, i.e., we fix an order for the elements of $\Irr_{\ell'}(\fS_{\ell^i})$.

Let $n=\sum\limits_{i\ge 0}\beta_i \ell^i$ be an $\ell$-expansion of $n$ and let $Y$ be the corresponding $\ell$-Young subgroup of $G$ so that $Y\cong\prod\limits_{i\ge 0}(\fS_{\ell^i})^{\beta_i}$.  Then
\[\nu_\ell(Y)=\sum_{i\ge 0}\beta_i\nu_\ell(\ell^i!)=\frac{n-\sum\limits_{i\ge 0}\beta_i}{\ell-1}.\]
Let $(Y,\ze)$ be an $\ell$-Young pair of $G$ and $\la\in\dz(\N_G(Y)_\ze/Y)$.
We write $\ze=\prod\limits_{i\ge0}\prod\limits_{t=1}^{u_i}\ze_{i,t}^{\beta_{i,t}}$
where $\ze_{i,t}\in\Irr_{\ell'}(\fS_{\ell^i})$ with $\ze_{i,t}\ne\ze_{i,t'}$ if
$t\ne t'$ and $\beta_i=\sum\limits_{t=1}^{u_i}\beta_{i,t}$.
From this, $\N_G(Y)_\ze/Y\cong \prod\limits_{i\ge 0}\prod\limits_{t=1}^{u_i}\fS_{\beta_{i,t}}$.
Therefore, we may write $\la=\prod\limits_{i\ge 0}\prod\limits_{t=1}^{u_i}\la_{i,t}$ with $\la_{i,t}\in\dz(\fS_{\beta_{i,t}})$.
Define an $\ell$-core tower $(\ka_{i,j})_{i\ge 0,\ 1\le j\le\ell^i}$ corresponding to $(Y,\ze,\la)$: we let $\ka_{i,j}=\emptyset$ if none of $\ze_{i,t}$ ($1\le t\le u_i$) equals $\xi_{i,j}$, while we let $\ka_{i,j}$ be the $\ell$-core partition of $\beta_{i,k}$ corresponding to the irreducible character $\la_{i,t}$ if $\ze_{i,t}=\xi_{i,j}$.

This defines a bijection between the $G$-conjugacy classes of triples $(Y,\ze,\la)$ with $(Y,\ze)\in \cY_G$ and $\la\in\dz(\N_G(Y)_\ze/Y)$ and the $\ell$-core towers $(\ka_{i,j})_{i\ge 0,\ 1\le j\le\ell^i}$ with $n=\sum\limits_{i\ge 0}\sum\limits_{j=1}^{\ell^i}|\ka_{i,j}|\ell^i$, and thus completes the proof for the case $G=\fS_n$.

Next, we let $G=C_e\wr\fS_n$ (with $\ell\nmid e$). The irreducible characters of $G$ are in bijection with $e$-tuples $(\mu_1,\ldots,\mu_e)$ with $|\mu_1|+\cdots+|\mu_e|=n$. If $\chi\in\Irr(G)$ is labeled by $(\mu_1,\ldots,\mu_e)$, then it can be shown that
\[ \de(\chi)=\frac{n-\sum\limits_{k=1}^e\sum\limits_{i\ge 0}\beta_i(\ell,\mu_k)}{\ell-1}.\]

For every $i\ge 0$, we write $\Irr_{\ell'}(C_e\wr\fS_{\ell^i})=\{\, \xi_{k,i,j}\mid 1\le k\le e, 1\le j\le \ell^i \,\}$ by Lemma~\ref{lem-ell'char}, i.e., we fix an order for the elements of $\Irr_{\ell'}(C_e\wr\fS_{\ell^i})$. Here we assume further that $\xi_{k,i,j}$ and $\xi_{k',i,j'}$ cover the same character in $\Irr(C_e^{\ell^i})$ if and only $k=k'$.

Let $n=\sum\limits_{i\ge 0}\beta_i \ell^i$ be an $\ell$-expansion of $n$ and let $Y$ be the corresponding $\ell$-Young subgroup of $G$ so that $Y\cong \prod\limits_{i\ge 0}(C_e\wr\fS_{\ell^i})^{\beta_i}$.  Then
$\nu_\ell(Y)=(n-\sum\limits_{i\ge 0}\beta_i)/(\ell-1)$.
Let $(Y,\ze)$ be an $\ell$-Young pair of $G$ and $\la\in\dz(\N_G(Y)_\ze/Y)$.
We write $\ze=\prod\limits_{i\ge 0}\prod\limits_{t=1}^{u_i}\ze_{i,t}^{\beta_{i,t}}$ where
$\ze_{i,t}\in\Irr_{\ell'}(C_e\wr\fS_{\ell^i})$ with $\ze_{i,t}\ne\ze_{i,t'}$ if
$t\ne t'$ and $\beta_i=\sum\limits_{t=1}^{u_i}\beta_{i,t}$.
Moreover, $\N_G(Y)_\ze/Y\cong \prod\limits_{i\ge 0}\prod\limits_{t=1}^{u_i}\fS_{\beta_{i,t}}$.
Therefore, we may write $\la=\prod\limits_{i\ge 0}\prod\limits_{t=1}^{u_i}\la_{i,t}$ with $\la_{i,t}\in\dz(\fS_{\beta_{i,t}})$.
For each $1\le k\le e$, we define an $\ell$-core tower $(\ka_{k,i,j})_{i\ge 0,\ 1\le j\le\ell^i}$ corresponding to $(Y,\ze,\la)$: we let $\ka_{k,i,j}=\emptyset$ if none of $\ze_{i,t}$ ($1\le t\le u_i$) equals $\xi_{k,i,j}$, while we let $\ka_{k,i,j}$ be the $\ell$-core partition of $\beta_{i,t}$ corresponding to $\la_{i,t}$ if $\ze_{i,t}=\xi_{k,i,j}$.

This defines a bijection between the $G$-conjugacy classes of triples $(Y,\ze,\la)$ with $(Y,\ze)\in \cY_G$ and $\la\in\dz(\N_G(Y)_\ze/Y)$ and the $\ell$-core towers $((\ka_{k,i,j})_{i\ge 0,\ 1\le j\le\ell^i})_{1\le k\le e}$ with $n=\sum\limits_{k=1}^e\sum\limits_{i\ge 0}\sum\limits_{j=1}^{\ell^i}|\ka_{k,i,j}|\ell^i$, and thus completes the proof for the case $G=C_e\wr\fS_n$.

Finally, let $G=G(2e,2,n)$ (with $\ell\nmid 2e$).
Note that $G$ is of index 2 in $\tG:=G(2e,1,n)=C_{2e}\wr\fS_n$.
We identify $G$ with $N\rtimes \fS_n$ where $N$ is a subgroup of $C_{2e}^n$ of index 2.
Let $\tchi\in\Irr(\tG)$. Then $\tchi$ covers one or two irreducible characters of $G$. Let $\theta\in\Irr(C_{2e}^n)$, and let $\widehat\theta$ be the extension of $\theta$ to $\tG_{\theta}$ as in Lemma~\ref{lem:exten-linear}.
Then $\tchi=\Ind_{\tG_{\theta}}^{\tG}(\widehat\theta\eta)$ where $\eta\in\Irr(\tG_{\theta}/C_{2e}^n)$.
Let $\tau$ be the non-trivial linear character of $\tG/G$ which by restriction can be also regarded as a character of $C_{2e}^n/N$.
Then \[\tau\tchi=\tau\Ind_{\tG_{\theta}}^{\tG}(\widehat\theta\eta)=\Ind_{\tG_{\theta}}^{\tG}(\tau\widehat\theta\eta)=\Ind_{\tG_{\theta}}^{\tG}(\widehat{\tau\theta}\eta),\]
and thus $\tau\tchi=\tchi$ if and only if $\theta$ and $\tau\theta$ are $\tG$-conjugate by construction.
So by Clifford theory, $\tchi$ covers two irreducible characters of $G$ if and only if $\tchi$ corresponds to $(\mu_1,\mu_2,\ldots,\mu_{2e})$ such that $\mu_i=\mu_{e+i}$ for $1\le i\le e$.

Let $Y'$ be an $\ell$-Young subgroup of $\fS_n$, $\widetilde Y=C_{2e}\wr Y'$ and $Y=N\rtimes Y'$.
For $\widetilde \ze\in\Irr(\widetilde Y)$ and $\ze\in\Irr(Y\mid\widetilde \ze)$, $(\widetilde Y,\widetilde\ze)\in\cY_{\tG}$ if and only if $(Y,\ze)\in\cY_{G}$.
By the proof of Lemma~\ref{lem-ell'char}(b), if $\widetilde \ze\in\Irr_{\ell'}(\widetilde Y)$ and $\widetilde\eta\in\Irr((C_{2e})^n\mid\widetilde \ze)$, then $\widetilde\ze=\widehat{\widetilde\eta}\ze'$ where $\ze'\in\Irr_{\ell'}(Y')$, and $\widehat{\widetilde\eta}\in\Lin(\widetilde Y)$ is defined as in Lemma~\ref{lem:exten-linear}.
As $\tau\widetilde\ze=\tau\widehat{\widetilde\eta}\ze'=\widehat{\tau\widetilde\eta}\ze'$, we see that $\ze:=\Res^{\widetilde Y}_Y(\widetilde\ze)$ is irreducible.
Thus $\N_{\tG}(\widetilde Y)_{\widetilde\ze}/\widetilde Y$ can be regarded as a subgroup of $\N_{G}(Y)_\ze/Y$.
More precisely, $\N_{G}(Y)_\ze/Y\cong \N_{\tG}(\widetilde Y)_{\widetilde\ze}/\widetilde Y$ if and only if $\tau\widetilde\ze$ and $\widetilde\ze$ are not $\N_{\tG}(\widetilde Y)$-conjugate, and when $\tau\widetilde\ze$ and $\widetilde\ze$ are $\N_{\tG}(\widetilde Y)$-conjugate, we have $\N_{G}(Y)_\ze/Y\cong (\N_{\tG}(\widetilde Y)_{\widetilde\ze}/\widetilde Y)\ti C_2$.
For $\la\in\dz(\N_G(Y)_\ze/Y)$ and $\widetilde\la\in\dz(\N_{\tG}(\widetilde Y)_{\widetilde\ze}/\widetilde Y)$, if moreover $\widetilde\la\in\Irr(\N_G(Y)_\ze/Y\mid\la)$, then we say $(\widetilde Y,\widetilde \ze,\widetilde\la)$ covers $(Y,\ze,\la)$.
Therefore, it suffices to show that if $\tchi\in\Irr(\tG)$ corresponds to the triple $(\widetilde Y,\widetilde\ze, \widetilde\la)$, then the number of irreducible characters of $G$ covered by $\tchi$ is equal to the number of $G$-conjugacy classes of triples $(Y,\ze,\la)$ of $G$ covered by $(\widetilde Y,\widetilde\ze, \widetilde\la)$.
By the above arguments, $(\widetilde Y,\widetilde\ze, \widetilde\la)$ covers one or two triples $(Y,\ze,\la)$ of $G$, and the number is two if and only if $\tau\widetilde\ze$ and $\widetilde\ze$ are $\N_{\tG}(\widetilde Y)$-conjugate, and hence if and only if $(\widetilde Y,\widetilde\ze, \widetilde\la)$ corresponds to $((\ka_{k,i,j})_{i\ge 0,\ 1\le j\le\ell^i})_{1\le k\le 2e}$ such that $\ka_{k,i,j}=\ka_{k+e,i,j}$ for $1\le i\le e$ and any $k,j$.
This completes the proof.
\end{proof}

Therefore, we have complete the proof of Theorem~\ref{thm:unipotent} by combing
Proposition~\ref{prop:bijection exc} and \ref{prop-bijection-sym}.

%%%%%%%%%%%%%%%%%%%%%%%%%%%%%%%%%%%%%%%%%%%%%%%%%%%%%%%%%%%%%%%%%%%%%%%%%
%%%%%%%%%%%%%%%%%%%%%%%%%%%%%%%%%%%%%%%%%%%%%%%%%%%%%%%%%%%%%%%%%%%%%%%%%

\end{document}